\documentclass[a4paper,10pt]{amsart}
\usepackage{epsf}  
\usepackage{amsfonts, amssymb, amsthm, amsmath}  
\usepackage{hyperref}
\usepackage{tikz}
\usetikzlibrary{matrix,arrows,decorations.pathmorphing}
\usepackage{tikz-cd}
\usepackage{graphicx}
\usepackage{psfrag}
\usepackage{amsthm}
\newtheorem{thm}{Theorem}  
\newtheorem*{*thm}{Theorem}

\newtheorem{cor}[thm]{Corollary}  
\newtheorem{lemma}[thm]{Lemma}  
\newtheorem{remark}[thm]{Remark}  
\newtheorem{defn}[thm]{Definition}  
\newtheorem{prop}[thm]{Proposition}  
\newtheorem{claim}[thm]{Claim}  
\newtheorem{example}[thm]{Example}  
\numberwithin{thm}{section}  
\def\pf{\noindent\emph{Proof: }}  
\def\stop{\hfill$\square$}  
\def\co{\colon\thinspace}

\providecommand{\ov}[1]{\hspace{-.1cm}\downarrow_{#1}}

\providecommand{\totl}[1]{\ensuremath{\lceil #1\rceil }}
\providecommand{\totb}[1]{\ensuremath{\underline{ #1}}}
\DeclareMathOperator{\End}{End}

\newcommand{\ro}{{}^{r}\Omega}
\providecommand{\rof}{{}^{\phantom{f}r}_{fg}\Omega}
\newcommand{\rhf} {{}^{\phantom{f}r}_{fg}H}
\newcommand{\rh}{{}^{r}H}
\newcommand{\ex}{\bold}

\providecommand{\C}[2]{\ensuremath {C^{#1,\underline{#2}}}}

\newcommand{\tc}[1]{\check\rvert_{#1}}

\newcommand{\dmsw}{\mathcal M^{st}_{\bullet}}
\newcommand{\dmod}{\mathcal M_{\bullet}}

\newcommand{\Msw}{\mathcal M^{st}}

\providecommand{\fp}[2]{{}_{\hspace{3pt}#1\hspace{-2pt}}\times_{#2}}
\newcommand{\exte}{\subset_{e}}
\newcommand{\rexte}{\ {}_{e}\!\!\supset}

\DeclareMathOperator{\expl}{Expl}

\newcommand{\dbar}{\bar{\partial}}

\providecommand{\et}[2]{\ensuremath{\bold T^{#1}_{#2}}}
\providecommand{\lrb}[1]{\ensuremath{\left(#1\right)}}
\providecommand{\abs}[1]{\left\lvert #1\right\rvert}

\author{Brett Parker   }
\email{brettdparker@gmail.com}  

\title[Holomorphic curves in exploded manifolds: VFC]{Holomorphic curves in exploded manifolds: virtual fundamental class}

\begin{document}
\maketitle

\begin{abstract}We define Gromov--Witten invariants of exploded manifolds.   The technical heart of this paper is a construction of a  virtual fundamental class $[\mathcal K]$ of any Kuranishi category $\mathcal K$ (which is a simplified, more general version of an embedded Kuranishi structure.) We also show how to integrate differential forms over $[\mathcal K]$ to obtain numerical invariants, and push forward such differential forms over suitable   maps.  We show that such invariants are independent of  any choices, and are compatible with pullbacks, products, and tropical completion of Kuranishi categories.

In the case of a compact symplectic manifold, this gives an alternative construction of Gromov--Witten invariants, including gravitational descendants.
\end{abstract}

\tableofcontents

\newpage

\section{Introduction}

\

In this paper, we construct Gromov--Witten invariants of  exploded manifolds\footnote{ Definitions for exploded manifolds can be found in \cite{iec}. For a short introduction to exploded manifods, see \cite{scgp}. } using embedded Kuranishi structures, defined and constructed in \cite{evc}. As smooth manifolds are a full subcategory of the category of exploded manifolds, this paper gives
an alternative construction --- and proof of the invariance of --- all descendant Gromov--Witten invariants of any compact symplectic manifold.  

Our construction also provides  Gromov--Witten invariants relative normal-crossing divisors.
Given a K\"ahler manifold with normal-crossing divisors, we can explode it, then construct Gromov--Witten invariants of the resulting exploded manifold.   There is a similar construction for symplectic manifolds with normal-crossing symplectic divisors, however we must be more careful here: If we take Definition 2.1 from \cite{zingerNC}, then a simple crossing symplectic divisor is a finite transverse collection of closed codimension $2$ symplectic submanifolds whose  intersection is symplectic, with symplectic orientation agreeing with the intersection orientation. After deforming the symplectic structure in such a configuration, it admits a contractible choice of $\dbar$--log compatible almost complex structure as in Section 14 of \cite{elc}, and we can then apply the explosion functor to get an exploded manifold and define relative Gromov--Witten invariants. The isotopy class of the $\dbar$--log compatible almost complex structure only depends on the isotopy class of the simple crossing symplectic divisor, so again we may define relative Gromov--Witten invariants in this setting as the Gromov--Witten invariants of the associated exploded manifold.

In some cases, the Gromov--Witten invariants defined in this paper coincide with previously-defined Gromov--Witten invariants of symplectic manifolds defined by   Fukaya and Ono in \cite{FO}, McDuff in \cite{McDuff}, Ruan in \cite{Ruanvirtual}, Liu and Tian  in \cite{Liu-Tian}, Siebert in \cite{Siebert}, and Li and Tian in \cite{Tian-Li}. 
I expect that in the algebraic case, the definition of Gromov--Witten invariants given here will agree with the algebraic definition given by Behrend and Fantechi in \cite{BehrendFantechi}, and also log Gromov--Witten invariants defined by Gross and Siebert in \cite{GSlogGW} and Abramovich and Chen in \cite{Chen, acgw}.

\

The technical heart of this paper is a construction of a virtual fundamental class associated to an embedded Kuranishi structure.  As well as the original work of Fukaya and Ono in \cite{FO}, and Li and Tian in \cite{Tian-Li}, there has been much recent work on constructing virtual fundamental cycles from various kinds of Kuranishi structures; see \cite{KFOOO, KMW, KMW1, pardon, MW2, MW3, CLW, joycebook, joyceKS}. In the case of smooth manifolds, this paper provides an alternate construction of a virtual  fundamental cycle from an embedded  Kuranishi structure (constructed in \cite{evc}).

An embedded Kuranishi structure consists of a collection of compatible charts $(\hat f,G, V)$ where $\hat f$ is a family of curves, $G$ is a group of automorphisms of $\hat f$, and $V$ is some $G$--equivariant obstruction bundle over $\hat f$  with a natural $G$--invariant section $\dbar \hat f$. The adjective `embedded' indicates that $\hat f$ is a family of curves in some moduli stack $\dmsw$,  $V$ is defined on a neighourhood of $\hat f$ in $\dmsw$ as a finite-rank sub-bundle of some natural obstruction bundle  with a natural section $\dbar$, and $\hat f/G$  represents the substack $\dbar^{-1}V$. We also require  compatibility of Kuranishi charts --- $(\hat f,G,V)$ is compatible with $(\hat f',G',V')$  if $V$ is a sub-bundle of $V'$ (or visa versa) on their common domain of definition.  So, we have a kind of transition map between charts 
\[\hat f \longleftarrow \hat f\times_{\dmsw}\hat f'\longrightarrow \hat f'\]
where the leftward arrow is a $G$--equivariant $G'$--fold cover of an open subfamily of $\hat f$, and the rightward arrow is a $G'$--equivariant, $G$--fold cover of the  subfamily $(\dbar\hat f')^{-1}V \subset \hat f'$.  

\

We shall often want to construct a global section of some sheaf over our Kuranishi charts. (For example, we might want to perturb $\dbar$ to be transverse to $0$, or construct a smooth function, or use the Chern-Weil construction to obtain the Chern class of a vectorbundle.)    To present a simple and unified construction of such global sections,  we introduce the notion of a $K$--category in section \ref{Kcat section} --- a $K$--category is obtained from an embedded Kuranishi structure by discarding all information apart from the charts $\hat f/G$ and their transition maps. In Proposition \ref{global section}, we prove that, at the expense of shrinking the charts in a $K$--category a little, we can construct a global section of any sheaf  satisfying three basic axioms, called `Patching', `Extension' and `Averaging'.   

Proposition \ref{global section} serves well to construct most of our global sections, but there is one important exception: the sheaf of transverse perturbations of the $\dbar$ equation does not satisfy the Averaging Axiom. 
In section \ref{wb section} we construct  a weighted branched cover $I$ of a Kuranishi category as a  `sheaf' of finite measure-spaces and define a weighted branched section of a sheaf $S$ to be a natural transformation $I\longrightarrow S$. We then prove that the corresponding sheaf of weighted branched sections  $S^{I}$ satisfies the Patching,  Extension, and Averaging axioms if $S$ satisfies the Patching and Extension axioms. We can then use Proposition \ref{global section} to construct global sections of $S^{I}$.

 By the end of section \ref{vmod construction}, we construct the virtual fundamental class $[\dmod]$ of the moduli stack of holomorphic curves --- where $\bullet$ indicates  choices such as specifying the genus, number of marked points and homology class of the curves under study.  This virtual fundamental class $[\dmod]$ is some weighted branched thingy in the moduli stack $\dmsw$, but in section \ref{de rham section} we show how to integrate differential forms  over $[\dmod]$, and also how to push forward differential forms along natural evaluation maps to finite-dimensional exploded manifolds or orbifolds. Such differential forms on $\dmsw$ may be obtained by pulling back differential forms from manifolds or orbifolds under natural evaluation maps, or obtained as Chern classes of any naturally-defined complex vectorbundle over  $\dmsw$, so we can define descendant Gromov--Witten invariants using Chern classes of tautological vectorbundles.

To simplify and emphasize the main points of our construction, we introduce the notion of a Kuranishi category, $\mathcal K$, and construct a virtual fundamental class $[\mathcal K]$ for any such Kuranishi category. Let us describe our results in terms of $[\mathcal K]$.

If $\theta$ is a differential form on $\mathcal K$ and $\mathcal K$ is compact\footnote{A compact Kuranishi category is one in which the subset consisting of holomorphic curves is compact; see Definition \ref{K proper}.} and oriented, then $\int_{[\mathcal K]}\theta$ is defined in section \ref{de rham section}. What type of differential form is $\theta$? Unlike in the case of smooth manifolds, there are several different types of differential forms that are useful on an exploded manifold $\ex B$.  
\[\rof^{*}(\ex B) \hookrightarrow \ro^{*}(\ex B) \hookleftarrow \Omega^{*}(\ex B)\]
All three types coincide with smooth differential forms in the case that $\ex B$ is a smooth manifold. De-Rham cohomology defined using $\Omega^{*}(\ex B)$ is much the same as usual cohomology; see \cite{dre}, Definition 1.2 and Corollary 4.2.  Refined cohomology, $\rh^{*}(\ex B)$, defined using refined differential forms in $\ro^{*}(\ex B)$ is usually infinite-dimensional, but admits pushforwards, and  acts as expected with fiber-products of exploded manifolds; see \cite{dre}, Definition 9.1,  Theorem 9.2, Lemma 9.3, and Lemma 9.5. The cohomology, $\rhf(\ex B)$, defined using refined differential forms generated by functions,  $\rof(\ex B)$, is also compatible with pushforwards and fiber products, but unlike $\rh^{*}(\ex  B)$ and $H^{*}(\ex B)$, is invariant only in families parametrized by connected smooth manifolds, rather that families parametrized by connected exploded manifolds.  See Definition \ref{rof} for $\rof^{*}$. The advantage of differential forms generated by functions is that they are compatible with tropical completion --- this is important for defining the contribution of a tropical curve to Gromov-Witten invariants, and for the tropical gluing formula for Gromow-Witten invariants, equation (1) of \cite{gfgw}. 

The integral $\int_{\mathcal K}\theta$ makes sense for $\theta\in\ro^{*}(\mathcal K)$, and therefore makes sense for any of the above types of differential forms.
 If $\mathcal K$ is complete\footnote{If $\mathcal K$ is an embedded Kuranishi structure for the moduli space of holomorphic curve in some complete exploded manifold, it is complete if and only if the corresponding moduli space of curves is compact. See \cite{cem} for cases in which such compactness holds.} (Definition \ref{K proper}) and $d\theta=0$, then $\int_{[\mathcal K]}\theta$ is independent of all choices in the construction of $\mathcal K$, and depends only on the cohomology class represented by $\theta$. The same holds with the weaker assumption that $\mathcal K$ is compact, and the stronger assumption that $\theta\in\rof^{*}\mathcal K$.  

Given a complete, relatively oriented map $\pi:\mathcal K\longrightarrow \ex A$ to an exploded manifold or orbifold $\ex A$, we can  integrate forms along the fiber of $\pi$ to define a map
\[\pi_{!}:\ro^{*}(\mathcal K)\longrightarrow \ro^{*}(\ex A) \] 
 inducing a map on cohomology independent of all choices involved in the construction of $[\mathcal K]$ or $\pi_{!}$.
\[\pi_{!}:\rh^{*}(\mathcal K)\longrightarrow \rh^{*}(\ex A)\]
\[\text{and }\pi_{!}:\rhf^{*}(\mathcal K)\longrightarrow \rhf^{*}(\ex A)\]
As usual, in the case that $\mathcal K$ is complete, 
\[\int_{[\mathcal K]}\pi^{*}\theta=\int_{\ex A}\theta\wedge\pi_{!}(1)\]
 for any closed differential form $\theta\in\ro^{*}(\ex A)$.
%In Theorem \ref{cobordism theorem}, we prove that pushforwards do not depend on choices involved in the construction of $[\mathcal K]$  --- often $\mathcal K$ is contained in some nice stack $\mathcal X$, and $\pi$  comes from a natural map $\mathcal X\longrightarrow \ex A$; in this case the corresponding maps 
%\[\pi_{!}:\rh^{*}(\mathcal X)\longrightarrow \rh^{*}(\ex A)\text{ and }\pi_{!}:\rhf^{*}(\mathcal X)\longrightarrow \rhf^{*}(\ex A)\]
%only depend on the cobordism class of $\mathcal K$ within $\mathcal X$.

 The following theorem gives that, on the level of cohomology,  $\pi_!$ only depends on the cobordism class of $\mathcal K$.

\begin{thm}\label{thm1} If  $\mathcal K_{0}$ and $\mathcal K_{1}$ are cobordant within a stack $\mathcal X$ with  a  map $\pi:\mathcal X\longrightarrow \ex A$, then, given any construction of $[\mathcal K_{j}]$, the two composite maps
\[\rh^{*}(\mathcal X)\longrightarrow \rh^{*}(\mathcal K_{j})\xrightarrow{\pi_{!}}\rh^{*}\ex A\] 
are equal, and the same holds for the analogous maps
\[\rhf^{*}(\mathcal X)\longrightarrow \rhf^{*}(\mathcal K_{j})\xrightarrow{\pi_{!}}\rhf^{*}\ex A\ .\] 
%Moreover given any complex vectorbundle $W$ over $\mathcal X$, and construction of any characteristic class,  $c(W)$ on $\mathcal K_{j}$ as in Remark \ref{chern},  the maps
%\[\theta\mapsto \pi_{!}(c(W))\wedge \iota_{j}^{*}\theta\]
%\[\begin{tikzcd}[row sep=tiny]\rh^{*}(\mathcal X)\rar &\rh^{*}(\ex A)
%\\\theta \rar[|->] &\pi_{!}(c(W)\wedge \iota_{j}^{*}\theta)
%\end{tikzcd}\]
%are equal on the level of cohomology for $j=0,1$, in both $\rh^*$ and $\rhf^*$. 
%
%If $\mathcal K_{j}$ are only proper over $\ex Z$, and $\ex A$ is a manifold, the corresponding maps
%\[\rhf^{*}(\mathcal X)\longrightarrow H^{*}(\ex A)\]
%are equal. 
\end{thm}

We also prove that Gromov--Witten invariants do not change in families of exploded manifolds, because they are compatible with pullbacks. Given a complete submersion $\mathcal K\longrightarrow \ex Z$ to an exploded manifold, we can pull back $\mathcal K$ over a map $\ex Z'\longrightarrow \ex Z$ to obtain another Kuranishi category $\mathcal K'\longrightarrow \ex Z'$. (For example, $\mathcal K$ might come from holomorphic curves in a family of exploded manifolds parametrized by $\ex Z$. Then $\mathcal K'$ is the Kuranishi category associated to holomorphic curves in the corresponding pulled-back family of exploded manifolds over $\ex Z'$.) When the map $\mathcal K\longrightarrow \ex Z$ factors through a map $\mathcal K\longrightarrow \ex A$, we get the following diagram of maps: 
  \[\begin{tikzcd}\mathcal K'\dar{\pi'}\rar &\mathcal K\dar{\pi} 
\\ \ex A\times_{\ex Z}\ex Z'\rar{y}\dar & \ex A\dar
\\\ex Z'\rar&\ex Z
\end{tikzcd}\]

\begin{thm}\label{thm2} The following diagrams commute:
\[\begin{tikzcd}\rh^{*}(\mathcal K')\dar{\pi'_{!}} &\lar \rh^{*}(\mathcal K)\dar{\pi_{!}} & \rhf^{*}(\mathcal K')\dar{\pi'_{!}} &\lar \rhf^{*}(\mathcal K)\dar{\pi_{!}} 
\\\rh^{*}( \ex A\times_{\ex Z}\ex Z') &\lar{y^{*}} \rh^{*}(\ex A) & \rhf^{*}( \ex A\times_{\ex Z}\ex Z') &\lar{y^{*}} \rhf^{*}(\ex A)
\end{tikzcd}\]
\end{thm}
%
%\begin{thm}\label{thm2}
%Suppose that $\mathcal K$ is complete over an exploded manifold $\ex Z$. Let $\mathcal K'$ be the pullback of $\mathcal K$ over  $\ex Z'\longrightarrow \ex Z$. Given compatible   maps,
%\[\begin{tikzcd}\mathcal K'\dar{\pi'}\rar &\mathcal K\dar{\pi} 
%\\ \ex A\times_{\ex Z}\ex Z'\rar{y}\dar & \ex A\dar
%\\\ex Z'\rar&\ex Z
%\end{tikzcd}\]
%the following diagrams commute:
%\[\begin{tikzcd}\rh^{*}(\mathcal K')\dar{\pi'_{!}} &\lar \rh^{*}(\mathcal K)\dar{\pi_{!}} & \rhf^{*}(\mathcal K')\dar{\pi'_{!}} &\lar \rhf^{*}(\mathcal K)\dar{\pi_{!}} 
%\\\rh^{*}( \ex A\times_{\ex Z}\ex Z') &\lar{y^{*}} \rh^{*}(\ex A) & \rhf^{*}( \ex A\times_{\ex Z}\ex Z') &\lar{y^{*}} \rhf^{*}(\ex A)
%\end{tikzcd}\]
%
%\end{thm}

 To prove gluing theorems, we also need fiber products of Kuranishi categories as well as pullbacks. However,   the (fiber)product of Kuranishi categories  usually has incompatible charts.  In section \ref{weak product section}, this problem is solved by shrinking charts in the product of Kuranishi categories to obtain a `weak' product of Kuranishi categories.
 
 \begin{thm}\label{thm3}Suppose that $\mathcal K$ is a weak product of some finite collection of complete, oriented Kuranishi categories  $\mathcal K_{v}$ with  maps $\pi_{v}:\mathcal K_{v}\longrightarrow \ex A_{v}$, and let $\pi:\mathcal K\longrightarrow \prod_{v}\ex A_{v}$ be the induced  map. Then the following diagrams commute:

\[\begin{tikzcd}\rh^{*}(\mathcal K)\rar{\pi_{!}}&\rh^{*}(\prod_{v}\ex A_{v}) & \rhf^{*}(\mathcal K)\rar{\pi_{!}}&\rhf^{*}(\prod_{v}\ex A_{v})
\\ \prod_{v}\rh^{*}(\mathcal K_{v})\uar\ar{ur}[swap]{\prod_{v}(\pi_{v})_{!}}&&  \prod_{v}\rhf^{*}(\mathcal K_{v})\uar\ar{ur}[swap]{\prod_{v}(\pi_{v})_{!}}\end{tikzcd}\]
\end{thm}
 
 Theorem \ref{thm3} combines with Theorem \ref{thm2} to show that integration over virtual fundamental classes acts as expected under fiber products.  Theorems \ref{thm1}, \ref{thm2} and \ref{thm3} follow immediately from Theorems \ref{cobordism theorem}, \ref{pullback theorem}  and \ref{weak product theorem}.

\

% 
% states that, if $\mathcal K$ is a weak product of  Kuranishi categories, $\mathcal K_{v}$, with complete maps $\mathcal K_{v}\longrightarrow \ex A_{v}$, then the following diagrams commute:
%\[\begin{tikzcd}\rh^{*}(\mathcal K)\rar{\pi_{!}}&\rh^{*}(\prod_{v}\ex A_{v}) & \rhf^{*}(\mathcal K)\rar{\pi_{!}}&\rhf^{*}(\prod_{v}\ex A_{v})
%\\ \prod_{v}\rh^{*}(\mathcal K_{v})\uar\ar{ur}[swap]{\prod_{v}(\pi_{v})_{!}} & & \prod_{v}\rhf^{*}(\mathcal K_{v})\uar\ar{ur}[swap]{\prod_{v}(\pi_{v})_{!}}\end{tikzcd}\]
%
%
%
%\begin{*thm}[\ref{weak product theorem}]Suppose that $\mathcal K$ is a weak product of some finite collection of complete, oriented Kuranishi categories  $\mathcal K_{v}$ with  maps $\pi_{v}:\mathcal K_{v}\longrightarrow \ex A_{v}$, let $\pi:\mathcal K\longrightarrow \prod_{v}\ex A_{v}$ be the induced  map. Then the following diagrams commute:
%
%\[\begin{tikzcd}\rh^{*}(\mathcal K)\rar{\pi_{!}}&\rh^{*}(\prod_{v}\ex A_{v}) & \rhf^{*}(\mathcal K)\rar{\pi_{!}}&\rhf^{*}(\prod_{v}\ex A_{v})
%\\ \prod_{v}\rh^{*}(\mathcal K_{v})\uar\ar{ur}[swap]{\prod_{v}(\pi_{v})_{!}}&&  \prod_{v}\rhf^{*}(\mathcal K_{v})\uar\ar{ur}[swap]{\prod_{v}(\pi_{v})_{!}}\end{tikzcd}\]
%\end{*thm}
%
%

Each exploded manifold $\ex B$ has a tropical part, $\totb{\ex B}$, which describes the (infinitely) large-scale structure of $\totb{\ex B}$;  the fiber, $\ex B\vert _{p}$, of $\ex B\longrightarrow \totb{\ex B}$ over any point $p\in \totb{\ex B}$ is a smooth manifold. The integral of a differential form $\theta$ over an exploded manifold $\ex B$ is defined as $\sum_{p\in\totb{ \ex B}}\int_{\ex B\vert_{p}}\theta$. (Although $\totb{\ex B}$ is generally uncountable, only finitely many terms in this sum are nonzero when the integral of $\theta$ is defined). Moreover, if $\theta\in\rof^{*}(\ex B)$ is closed,  then $\int_{\ex B\rvert_{p}}\theta$ only depends on the homology class of $\theta$ in $\rhf^{*}(\ex B)$. In section  \ref{tropical completion section} we prove analogous results for integration over $[\mathcal K]$. 

One issue is that $\ex B\rvert_{p}$ is usually not compact or complete, even when $\ex B$ is complete. We deal with this issue using the	  tropical completion $\ex B\tc p$ of $\ex B$ at $p$, defined at the start of section \ref{tropical completion section}.  This tropical completion $\ex B\tc p$  always contains $\ex B\rvert_{p}$ as a dense subset, and is complete if $\ex B$ is compact. We can apply tropical completion to maps, differential forms, and Kuranishi categories. If $\theta\in\rof^{*}(\ex B)$, then $\theta\tc p\in\rof^{*}(\ex B\tc p)$ and $\int_{\ex B\rvert_{p}}\theta=\int_{\ex B\tc p}\theta\tc p$.  

Lemma \ref{tropical completion lemma} states that the integral of a closed form $\theta\in\rof^{*}(\mathcal K)$ breaks up into invariantly-defined contributions for each point $p$ in the tropical part of $\mathcal K$, and in particular, 
\[\int_{[\mathcal K]}\theta=\sum_{p\in\totb{\mathcal K}}\int_{[\mathcal K\tc p]}\theta\tc p\ .\]
Similarly, Lemma \ref{tropical completion lemma} states that,  for a complete relatively-oriented map $\pi:\mathcal K\longrightarrow \ex A$,  cohomology class  $\theta\in \rhf^{*}\mathcal K$, and $p'\in\totb{\ex A}$, 
\[\pi_{!}(\theta)\tc {p'}=\sum_{p\in\totb{\pi}^{-1}(p')}(\pi\tc {p})_{!}(\theta\tc {p})\ .\]

This paper concludes with section \ref{GW section}, which summarizes our  construction of Gromov--Witten invariants.

\subsection{Acknowledgements}
This research was supported by the Australian Research Council. Special thanks is due to an anonymous referee, who suggested several improvements to both this paper, and \cite{reg}, on which it depends.

\section{Constructing sections of sheaves on $K$--categories}
\label{Kcat section}

Throughout this paper, we will be using exploded manifolds with the regularity $\C\infty1$, defined in section 7 of \cite{iec}. For all practical purposes, $\C\infty1$  maps are as good as smooth. Let $\dmsw$ indicate some decorated moduli stack of $\C\infty1$ families of not-necessarily-holomorphic stable curves --- see section 11 of \cite{iec} for basic definitions, and section 2 of \cite{evc} for further treatment, including Definition 2.11 of  $\dmsw$. By a sheaf (of sets)  on $\dmsw$, we mean a contravariant functor $S$ from $\dmsw$  (to the category of sets) so that $S$  is a sheaf when restricted to the category of open sub-families of  any family in $\dmsw$. We shall also be interested in sheaves with more restricted domains. 

A stack $\mathcal X$ over the category of $\C\infty1$ exploded manifolds is a category $\mathcal X$  along with a `nice' functor $\ex F$ from $\mathcal X$ to the category of $\C\infty1$ exploded manifolds; see \cite{stacks} for an approachable introduction to stacks, and see section 2 of \cite{evc} for a study of the stack of $\C\infty1$ curves. In this paper, `stack' without further qualification will generally mean a stack over the category of $\C\infty 1$ exploded manifolds.  The `nice' properties of $\ex F$ are loosely paraphrased as follows: families (parametrized by exploded manifolds) glue and  pull back as we expect bundles (parametrized by exploded manifolds) to glue and pull back. Moreover, morphisms between families also restrict, pullback, and glue as expected for fiberwise isomorphisms.  We use $\hat f$ to indicate an object in $\mathcal X$, and call $\hat f$ a family parametrized by $\ex F(\hat f)$ --- this can also be thought of as a map $\ex F(\hat f)\longrightarrow \mathcal X$. The case we are most interested in is when $\mathcal X$ is a stack of curves, so $\hat f$ can be thought of as a family of curves, and any $f\longrightarrow \hat f$ with $\ex F(f)$ a point can be thought of as a curve $f$ in $\hat f$.  As in section 2 of \cite{evc}, by a substack $\mathcal U\subset \mathcal X$, we mean a full subcategory so that $\hat f$ is in $\mathcal U$ if and only if every $f$ in $\hat f$ is in $\mathcal U$.  Then, for any family $\hat f$ in $\mathcal X$, there is a corresponding subset $\mathcal U(\hat f)\subset \ex F(\hat f)$ so that for any $\hat g\longrightarrow \hat f$, $\hat g$ is in $\mathcal U$ if and only if $\ex F(\hat g)$ has image in $\mathcal U(\hat f)$. We call $\mathcal U$ an open substack if $\mathcal U(\hat f)\subset \ex F(\hat f)$ is open for all $\hat f$ in $\mathcal X$.  It is proved in Lemma 2.7 of \cite{evc} that this topology on the moduli stack of $\C\infty1$ curves matches the topology used in \cite{cem} to prove compactness results for the moduli stack of holomorphic curves.

When we take an embedded Kuranishi structure, and throw away all information apart from charts and their embedding into a stack, we obtain a $K$--category, defined below. See also Remark \ref{translation 1} for a definition of a $K$--category using charts not embedded in a stack.  

\begin{defn}\label{K category} A $K$--category   is a full  subcategory $\mathcal K$ of a stack $\mathcal K^{st}$ and a collection of charts $\hat f_{i}/G_{i}$   so that the following holds.
\begin{enumerate}
\item Each  $\hat f_{i}$ is a family in $\mathcal K^{st}$, and $G_{i}$ is a finite group of automorphisms of $\hat f_{i}$. 
\item \label{K category ss}$\hat f_{i}/G_{i}$ represents a substack of $\mathcal K^{st}$.
\item\label{K category t} Each family $\hat f_{i}$ has some fixed dimension, and whenever $\dim \hat f_{i}\leq \dim\hat f_{j}$, 
\[\hat f_{i}\times_{\mathcal K^{st}}\hat f_{j}\longrightarrow \hat f_{i}\]
is a $G_{j}$--fold cover of an open subfamily of $\hat f_{i}$, and
\[\hat f_{i}\times_{\mathcal K^{st}}\hat f_{j}\longrightarrow \hat f_{j}\]
is a $G_{i}$--fold cover of a subfamily of $\hat f_{j}$ --- locally defined by the transverse vanishing of some $\mathbb R$--valued $\C\infty1$ functions on $\ex F(\hat f_{j})$. 
\item The families $\hat f_{i}$ cover $\mathcal K^{st}$. 
\item \label{K category lf}The set of charts is countable and  locally finite ---   for each $i$, there are only finitely many $j$ so that $\hat f_{i}\times_{\mathcal K^{st}}\hat f_{j}$ is nonempty.
\item $\mathcal K$ is the full subcategory of $\mathcal K^{st}$ consisting of families locally isomorphic to some $\hat f_{i}$.
\end{enumerate} 

Given another $K$--category $\mathcal K^{\sharp}$  with charts $\hat f_{i}^{\sharp}/G_{i}$, use the notation $\mathcal K\subset \mathcal K^{\sharp}$ if,
\begin{itemize}
\item $\mathcal K^{st}$ is a substack of $(\mathcal K^{\sharp})^{st}$,
\item $\mathcal K$ is a subcategory of $\mathcal K^{\sharp}$,
\item and $\hat f_{i}$ is a $G_{i}$--equivariant open subfamily of $\hat f_{i}^{\sharp}$.  
\end{itemize}
Moreover, say that $\mathcal K^{\sharp}$ is an extension of $\mathcal K$ and use the notation $\mathcal K\exte\mathcal K^{\sharp}$ if $\mathcal K\subset \mathcal K^{\sharp}$ and  the closure of $\hat f_{i}\subset \hat f_{i}^{\sharp}$ is closed in $(\mathcal K^{\sharp})^{st}$ --- equivalently, if $\hat f_{i}'$ indicates the closure of $\hat f_{i}\subset \hat f_{i}^{\sharp}$,  then, for all $j$,  $\hat f'_{i}\times_{(\mathcal K^{\sharp})^{st}}\hat f_{j}^{\sharp}$ has closed image in $\hat f_{j}^{\sharp}$.

Say that $\mathcal K$ is extendable if there exist extensions $\mathcal K\exte\mathcal K'\exte\mathcal K^{\sharp}$. \footnote{We require two extensions so that any extendable  $K$--category $\mathcal K$ will have an extendable extension --- If $\mathcal K\exte\mathcal K'\exte\mathcal K^{\sharp}$, then we can obtain an extendable extension of $\mathcal K$ by shrinking $\mathcal K'$ appropriately.}

\end{defn}

\begin{defn}A sheaf (of sets)  on a $K$--category, $\mathcal K$, is a contravariant functor, $S$, from $\mathcal K$  (to the category of sets) so that $S$ is a sheaf whenever restricted to the category of open sub-families of a given family in $\mathcal K$.

A sheaf on an extendable $K$--category is a sheaf defined on some  extension $\mathcal K'\rexte\mathcal K$.
 \end{defn}

Extendability is an important property for constructing sections of sheaves on $\mathcal K$. Any extendable  $\mathcal K$ has a partition of unity,   and Proposition \ref{global section} gives a way of constructing global sections of sheaves on $\mathcal K$. For example, if $\mathcal K$ is extendable,  each $G_{i}$ is trivial, and each $\ex F(\hat f_{i})$ is a $n$--dimensional manifold, then the $\hat f_{i}$ provide coordinate charts on a  $n$--dimensional manifold $M$, and $\mathcal K^{st}$ is the moduli stack of maps into $M$. If we drop the condition that $G_{i}$ be trivial, then $\hat f_{i}/G_{i}$ give charts on an orbifold, $M$ (and depending on your position on orbifolds,   $\mathcal K^{st}$ is either that orbifold, or is the category of maps  into $M$).  If, however, $\mathcal K$ is not extendable,  $M$ may not be Hausdorff.

 Global sections of a sheaf on $\mathcal K$ can be  constructed using the Patching, Extension, and Averaging axioms below. To construct Gromov--Witten invariants, we shall use the sheaf $S$ from Definition \ref{Sdef} below. This sheaf obeys the Patching and Extension axioms, but not the Averaging Axiom. Accordingly,  global sections of this sheaf $S$ may not exist, so we shall use weighted branched sections of $S$, which we regard as a natural transformation $I\longrightarrow S$, where $I$ is a  `sheaf'\footnote{This `sheaf' requires scare quotes because it has a restricted domain of definition.  } of finite measure-spaces.

The following axioms  allow the global construction of sections of a sheaf $S$ over an extendable $K$--category $\mathcal K$. 

\begin{itemize}
\item[\textbf{Patching}] Given
 an open cover $\{O_{k}\}$ of $\hat f$, and section $\theta_{k}\in S(O_{k})$ for all $k$, there exists a patched-together section $\theta\in S(\hat f)$ that agrees with $\theta_{k}$ wherever all $\theta_{k}$ agree. 
 
 So,  given any  morphism $\hat g\longrightarrow \hat f$ in $\mathcal K$ and  section $\theta'\in S(\hat g)$, the pullback of $\theta$  to $S(\hat g)$ is $\theta'$ if  the following condition is satisfied for all $k$:  $\theta'$ and $\theta_{k}$ pull back to the same section under the following diagram.
\[\begin{tikzcd}\hat g\times_{\hat f}O_{k}\rar\dar& O_{k}\dar
\\ \hat g \rar   &\hat f\end{tikzcd}\]

\item[\textbf{Extension}] $S(\hat f)$ is nonempty. Moreover, if $\hat g\longrightarrow \hat f$ is a morphism  in $\mathcal K$, then given any $\theta\in S(\hat g)$ and $f$ in $\hat g$,  there exists some $\theta'\in S(\hat f)$ so that the pullback of $\theta'$ to $S(\hat g)$ agrees with $\theta$ in a neighborhood of $f\in \hat g$. 

\item [\textbf{Averaging}] Given any $\theta\in S(\hat f)$ and a finite $G$--action on $\hat f$, there is a $G$--invariant section $\theta'$ agreeing  with $\theta$ wherever $\theta$ was already $G$--invariant --- so, given any $G$--equivariant map $\hat f'\longrightarrow \hat f$  pulling $\theta$ back to a $G$--invariant section, the pullback of $\theta$ agrees with the pullback of $\theta'$ over $\hat f'\longrightarrow \hat f$. 
 
\end{itemize}

For example, if $S(\hat f)$ is the set of $\C\infty1$ real-valued functions on $\ex F(\hat f)$, then $S$  obeys the above axioms.  The Extension Axiom for $S$ follows from Definition \ref{K category} part \ref{K category t}:  given any morphism $\hat g\longrightarrow \hat f$ in $\mathcal K$, the corresponding map  $\ex F(\hat g)\longrightarrow \ex F(\hat f)$ is always locally an isomorphism onto  an exploded submanifold  of $\ex F(\hat f)$. Patching can be achieved using partitions of unity, and averaging achieved by averaging over the action of $G$.

\begin{prop}\label{global section} Let $S$ be a sheaf on an extendable $K$--category, $K_{2}^{\sharp}$, and suppose that $S$ satisfies the Patching, Extension, and Averaging axioms. Given the following inclusions and extensions of $K$--categories, 
\[\begin{tikzcd}\mathcal K_{1}\rar{\exte}\dar[hook]&\mathcal K_{1}^{\sharp}\dar[hook]
\\\mathcal K_{2}\rar{\exte}&\mathcal K_{2}^{\sharp}\end{tikzcd}\]
%\[\mathcal K_{1}\exte \mathcal K_{1}^{\sharp}, \ \ \mathcal K_{2}\exte\mathcal K_{2}^{\sharp}, \ \ \mathcal K_{1}\subset \mathcal K_{2},\ \ \mathcal K_{1}^{\sharp}\subset\mathcal K^{\sharp}_{2}\]
and a global section $\theta$ of $S$ defined on $\mathcal K_{1}^{\sharp}$, there exists a global section of  $S$ defined on $\mathcal K_{2}$ and agreeing with $\theta$ when restricted to $\mathcal K_{1}$. 
\end{prop}

\pf

Use notation $\hat f_{i,\mathcal K_{j}}/G_{i}\subset \hat f_{i,\mathcal K_{j}^{\sharp}}/G_{i}$ for  charts on $\mathcal K_{j}$ and $\mathcal K_{j}^{\sharp}$ respectively, and
 index these charts by the natural numbers. We shall construct our section inductively, but at each step we must shrink the domain of definition a little. Accordingly, for each $(i,j)$ where $j$ is a non-negative integer, choose $G_{i}$--invariant open subfamilies $\hat f_{i,j}\subset\hat f_{i,\mathcal K_{2}^{\sharp}}$ so that

\begin{enumerate}

\item\label{i1} if $j\geq i$, then $\hat f_{i,\mathcal K_{2}}\subset\hat f_{i,j}$;
\item\label{i2} if $j<i$, then $\hat f_{i,\mathcal K_{1}}\subset \hat f_{i,j}\subset \hat f_{i,\mathcal K_{2}^{\sharp}}$, and $\hat f_{i,j}$ contains every curve in the closure of $\hat f_{i,\mathcal K_{1}}$ within $\mathcal K_{2}^{\sharp}$;
\item\label{i3} if $j\neq i$,  then $\hat f_{i,j}\subset \hat f_{i,j-1}$, and $\hat f_{i,j-1}$ contains every curve in the closure of $\hat f_{i,j}$ within $\mathcal K_{2}^{\sharp}$.
\end{enumerate}

Such families can be constructed as follows. Because $K_{2}^{\sharp}$ is extendable, there exists a family $\hat f_{i}^{\sharp}$ containing $\hat f_{i,\mathcal K_{2}^{\sharp}}$  as an open subfamily, so that the closure of $\hat f_{i,\mathcal K_2^{\sharp}}$ within $\hat f_{i}^{\sharp}$ contains the closure of $\hat f_{i,\mathcal K_2^{\sharp}}$ within $\mathcal K_{2}^{\sharp}$ (or any of our other $K$--categories). The same will then hold for the open sub-families $\hat f_{i,j}\subset \hat f_{i,\mathcal K_2^{\sharp}}$.
Choose $\hat f_{i,0}=\hat f_{i,\mathcal K_{1}^{\sharp}}$. Then,   choose a continuous function, $\rho$, on $\hat f_{i}^{\sharp}$ that restricts to be $1$ on $\hat f_{i,K_{1}}\subset \hat f_{i,0}$, and  $0$ outside $\hat f_{i,0}\subset \hat f_{i}^{\sharp}$. With this $\rho$, we then define $\hat f_{i,j}:=\rho^{-1}(j/i,1]$ for $j<i$. Next, define $\hat f_{i,i}:=\hat f_{i,\mathcal K_{2}^{\sharp}}$, and choose another continuous function $\lambda$ on $\hat f_{i}^{\sharp}$, equal to $0$ on $\hat f_{i,K_{2}}$, and $1$ outside $\hat f_{i,i}\subset \hat f_{i}^{\sharp}$. Finally, define $f_{i,i+n}:=\lambda^{-1}[0,2^{-n})$ for positive integers $n$.

Let us construct  our section using the following  inductive step:

\begin{claim}\label{ic} Suppose that a section $\theta_{j-1}$ of $S(\hat f_{i,j-1})$ has been chosen for all $i$, and that these sections are compatible --- so,  given any pair of morphisms $\iota_{1}:\hat f\longrightarrow \hat f_{i,j-1}$ and $\iota_{2}:\hat f\longrightarrow \hat f_{i',j-1}$, we have  $\iota_{1}^{*}\theta_{j-1}=\iota_{2}^{*}\theta_{j-1}$. Suppose moreover that the restriction of $\theta_{j-1}$ to $\hat f_{i,\mathcal K_{1}}$ agrees with our original section $\theta$.

Then there exists a compatible choice of sections $\theta_{j}$ of $S(\hat f_{i,j})$ for all $i$ so that $\theta_{j}\in S(\hat f_{i,j})$ is the restriction of $\theta_{j-1}\in S(\hat f_{i,j-1})$ for all $i\neq j$, and so that the restrictions of  $\theta_{j}$ and  $\theta$ to $S(\hat f_{i,\mathcal K_{1}})$ are equal.
\end{claim}

To prove Claim \ref{ic}, we need only define $\theta_{j}$ on $\hat f_{j,j}$. Let us first construct a local candidate, $\theta'$, for $\theta_{j}$. For a curve $f$ in $\hat f_{j,j}$,  let  $i$ be so that $f\in\hat f_{i, j-1}$, and $\hat f_{i,j-1}$ has the maximal dimension\footnote{Recall, from Definition \ref{K category} part \ref{K category lf}, that $f$ can only be contained in finitely many $\hat f_{i}$.} of all such $\hat f_{k,j-1} $ containing $f$. (If there is no $\hat f_{i,j-1}$ containing $f$, we are free to choose any section of $S$ on a neighborhood if $f$ that does not intersect any $\hat f_{i,j}$ for $i\neq j$.)  By using the Extension Axiom, (or by pulling back $\theta$ from $\hat f_{i,j-1}$), there exists a section $\theta'$ of $S$ on a neighborhood $U_{f}$ of $f$ in $\hat f_{j,j}$ so that, given any pair of morphisms $\iota_{1}:\hat f\longrightarrow U_{f}$ and $\iota_{2}:\hat f\longrightarrow \hat f_{i,j-1}$, $\iota_{1}^{*}\theta'=\iota_{2}^{*}\theta_{j-1}$.  

Let us check that $\theta'$ locally satisfies the required compatibility conditions.
 Because $\hat f_{i}$ has maximal dimension, we may choose $U_{f}$ small enough that Definition \ref{K category} part \ref{K category ss} implies  that, if $\hat f$ is contained in both $U_{f}$ and $\hat f_{i',j}/G_{i'}$ for any other $i'\neq j$, then $\hat f$ is contained in $\hat f_{i,j-1}/G_{i}$.  It follows that $\theta'\in S(U_{f})$ is compatible with $\theta_{j}=\theta_{j-1}\in S(\hat f_{i',j})$ for all $i'\neq j$. We should also ensure that $\theta'$ is compatible with $\theta$ on $\hat f_{j,K_{1}}$. Because $\hat f_{i,j}$ has maximal  dimension,  if $\hat f_{j,j}$ has larger dimension, $f$ is not in $\hat f_{j,\mathcal K_{1}}$. Therefore, we may also choose $U_{f}$ small enough that if $\hat f$ is contained in $U_{f}$ and $\hat f_{j,\mathcal K_{1}}$, then $\hat f$ is contained in $\hat f_{i,j}/G_{i}$, ensuring that $\theta'\in S(U_{f})$ is also compatible with $\theta\in S(\hat f_{j,\mathcal K_{1}})$. 

The Patching Axiom allows us to patch together these sections to a section, $\theta'\in S(\hat f_{j,j})$,  still compatible with $\theta_{j}$  on $\hat f_{i, j}$ for $i\neq j$ and  $\theta$ on $\hat f_{j,\mathcal K_{1}}$. Then the Averaging Axiom gives a $G_{j}$--equivariant section  $\theta_{j}$ still compatible with $\theta_{j}$ on $\hat f_{i,j}$ and  $\theta$ on $\hat f_{j,\mathcal K_{1}}$. Definition \ref{K category} part \ref{K category ss} ensures that  any $G_{j}$--equivariant section of $S(\hat f_{j,j})$ has the required compatibility property that its pullback does not depend on the choice of morphism. 

This completes the proof of Claim \ref{ic}. 

\

Using Claim \ref{ic} inductively, we can construct $\theta_{i}$ for all $i$. Any family $\hat f$ in $\mathcal K_{2}$ is everywhere locally isomorphic to some subfamily of $\hat f_{i,\mathcal K_{2}}\subset\hat f_{i,i}$. We may define $\theta\in S(\hat f)$ to be the pullback   of $\theta_{i}$. This defines the required global section of $S$ over $\mathcal K_{2}$. 
 
 \stop

\begin{remark}Note that Proposition \ref{global section} may also be used to construct global sections of sheaves $S$ on $\mathcal K^{st}$ because any global section of $S$ on $\mathcal K$ automatically pulls back to give a global section of $S$ on $\mathcal K^{st}$.\end{remark}

\begin{lemma}\label{function construction}Consider an extension $\mathcal K\exte\mathcal K^{\sharp}$ of extendable $K$--categories, and a corresponding pair of charts $\hat f_{i}/G_{i}\subset \hat f_{i}^{\sharp}/G_{i}$. Then any  (continuous or $\C\infty1$) function 
\[\rho:\hat f_{i}/G_{i}\longrightarrow \mathbb R\]
 admitting an extension to $\hat f_{i}^{\sharp}/G_{i}$ also extends to  all of $\mathcal K$.

Moreover, if $O\subset \mathcal K^{\sharp}$ is any open subset so that the support of $\rho$ on $\hat f^{\sharp}_{i}$ has closure (within $\mathcal K^{\sharp}$) contained in $O$, then we may choose the corresponding extension $\rho:\mathcal K\longrightarrow \mathbb R$ also to have support with closure contained in $ O$, and if $\rho$ is non-negative, we can choose its extension to also be non-negative. 

\end{lemma}
\pf

We shall construct the required function using Proposition \ref{global section}. Define
\begin{itemize}
\item $\mathcal K_{1}\subset\mathcal K$ to be the category of families locally isomorphic to an open subset of $\hat f_{i}$, 
\item $\mathcal K_{1}^{\sharp}\subset\mathcal K^{\sharp}$ to be the category of families locally isomorphic to an open subset of $\hat f_{i}^{\sharp}$,
\item $\mathcal K_{2}=\mathcal K$,
\item $\mathcal K_{2}^{\sharp}=\mathcal K^{\sharp}$.
\end{itemize}

Let $S$ be the sheaf with $S(\hat f)$  the set of  $\C\infty1$ functions on $\ex F(\hat f)$ whose support has  closure (in $\mathcal K^{\sharp}$) contained  in $O$. The Patching Axiom for $S$ is proved using a partition of unity, and the Averaging Axiom  holds, because it is possible to average a $\C\infty1$ function to obtain a $G$--invariant function. To apply Proposition \ref{global section}, we only need to prove that $S$ satisfies the Extension Axiom.

The Extension Axiom follows from Definition \ref{K category} part \ref{K category t} --- every morphism $\hat f\longrightarrow \hat g$ is locally an isomorphism onto an exploded submanifold (which must be locally closed because $\mathcal K^{\sharp}$ is extendable). In particular, this implies that any (continuous or $\C\infty1$) function defined on $\ex F(\hat f)$ locally extends to a (continuous or $\C\infty1$) function on $\ex F(\hat g)$.  If our original function had support with closure contained in $O$, our local extension may also be chosen to have support with closure contained in $O$. Therefore $S$ obeys the Extension Axiom, and Proposition \ref{global section} tells us that $\rho$ extends to a global section of $S$, which is a map

\[\rho:\mathcal K\longrightarrow \mathbb R \ .\]

The fact that the support of $\rho$ has closure contained in $O$ follows from the analogous fact for each of the individual functions $\hat f_{k}/G_{k}\longrightarrow \mathbb R$, and the local-finiteness condition for $\mathcal K^{\sharp}$ from Definition \ref{K category} part \ref{K category lf}. 

The proof in the case that $\rho$ is non-negative is identical, except we use a sheaf with  $S(\hat f)$ the set of non-negative $\C\infty1$ functions on $\ex F(\hat f)$ whose support has closure contained in $O$.

\stop

\begin{remark}\label{function remark} We can use Lemma \ref{function construction} to construct a non-negative $\C\infty1$ function on $\mathcal K$ with zero set any given closed  substack $C$ of $\mathcal K^{\sharp}\rexte\mathcal K$.

\end{remark}

In particular, Lemma \ref{function construction} implies that any non-negative $\C\infty1$ bump function on $\hat f_{i}^{\sharp}/G_{i}$, whose support  has closure contained in the complement of $C$, extends to a non-negative $\C\infty1$ function $\rho$ on $\mathcal K$, vanishing on $C$. As $\hat f_{i}^{\sharp}/G_{i}\setminus C$ is covered by the support of a countable collection of such bump functions, and $\mathcal K$ has a countable number of charts, there exists a sequence, $\rho_{k}$, of non-negative functions on $\mathcal K$, vanishing on $C$ and with support covering $\mathcal K\setminus C$.  Then there exists a sequence,  $\epsilon_{k}$, of positive numbers so that $\sum_{k}\epsilon_{k}\rho_{k}$ converges to a $\C\infty1$ function $\rho$ on $\mathcal K$. (As with smooth functions on smooth manifolds, convergence to a $\C\infty1$ function is equivalent to convergence in a countable sequence of norms, so we can always ensure convergence of a sum by multiplying each term by a suitably small constant. See Definition 7.5 of \cite{iec}.) Such a $\rho$ is non-negative and has zero set $C$, as required.

\section{Kuranishi categories}

Embedded Kuranishi structures are defined in section 2.9 of \cite{evc}. Each Kuranishi chart $(\mathcal U,V,\hat f/G)$ consists of some open substack $\mathcal U\subset \dmsw$, an obstruction bundle $V$ over $\mathcal U$, and a family of curves $\hat f$ in $\mathcal U$ with automorphism group $G$. The obstruction bundle,  $V$, is a  finite-rank  complex  vectorbundle over $\mathcal U$, and also  a nice\footnote{Technically, $V$ satisfies definitions 2.23, 2.24 and 2.25 of \cite{evc}.} sub-bundle of the sheaf $\mathcal Y$ that is the codomain of the $\dbar$ equation. Moreover,  $\hat f/G$ represents the moduli stack $\dbar^{-1} V\subset \mathcal U$, so $\dbar\hat f$ defines a $G$--equivariant section of the restriction, $V(\hat f)$, of  $V$ to $\hat f$. Kuranishi charts in an embedded Kuranishi structure have to be compatible in the sense that on $\mathcal U_{i}\cap \mathcal U_{j}$, either $V_{i}$ is a sub-bundle of $V_{j}$ or visa-versa. Moreover, our Kuranishi charts  have compatible extensions, so, we can define an extendable Kuranishi category by taking the charts $\hat f_{i}/G_{i}$.

\begin{defn}\label{aKC}Given an embedded Kuranishi structure $\{(\mathcal U_{i},V_{i},\hat f_{i}/G_{i})\}$ on $\dmsw$, define its associated Kuranishi category, $\mathcal K$, to be  the (full) subcategory of $\dmsw$ consisting of families locally isomorphic to an open subfamily of some $\hat f_{i}$, and define $\mathcal K^{st}$ to be the substack of $\dmsw$ consisting of curves isomorphic to curves in some $\hat f_{i}$; so, a family of curves $\hat f$ in $\dmsw$ is in $\mathcal K^{st}$ if each curve $f$ in $\hat f$ is isomorphic to a curve in some $\hat f_i$. This category $\mathcal K$ comes with the following extra structure:
\begin{itemize} 
\item the open substacks $\mathcal U_{i}\cap\mathcal K^{st}$, with the vectorbundles $V_{i}$,
\item the charts  $\hat f_{i}/G_{i}$,
\item the section $\dbar$ of $V_{i}(\hat f_{i})$ over $\ex F(\hat f_{i})$.
\end{itemize}
\end{defn}

In the next definition, we specify the essential properties of this extra structure on $\mathcal K$.

\begin{defn}\label{Kcat} A Kuranishi category is an extendable $K$--category $\mathcal K$ with charts $\hat f_{i}/G_{i}$ (Definition \ref{K category}) along with:
\begin{itemize}
\item open substacks $\mathcal U_{i}\subset\mathcal K^{st}$ containing $\hat f_{i}/G_{i}$ so that each $\mathcal U_{i}$ only intersects finitely many other $\mathcal U_{j}$;
\item constant-rank complex vectorbundles $V_{i}$ on $\mathcal U_{i}$, and, on $\mathcal U_{i}\cap\mathcal U_{j}$, an inclusion of one of $V_{i}$ or $V_{j}$ as a sub-bundle of the other;

\item and  sections $\dbar\hat f_{i}\co \ex F(\hat f_{i})\longrightarrow V_{i}(\hat f_{i})$
\end{itemize}
satisfying the following conditions.
\begin{enumerate}
\item \label{Kcat V}The sections $\dbar\hat f_{i}$ determine a global section, $\dbar$, of the sheaf with sections over $\hat f$ the sections of a vectorbundle, $V(\hat f)$, equal to $V_{i}(\hat f)$ wherever $\hat f$ is locally isomorphic to $\hat f_{i}$, and with pullbacks induced by the inclusions of vectorbundles above. So, $V$ is covariant functor, with $V(\hat f)\longrightarrow \ex F(\hat f)$ a  complex vectorbundle, and $\dbar$ is a natural transformation from $\ex F$ to $V$, with $\dbar\hat f$ a section of $V(\hat f)$. 
%
%\label{Kcat V}a lift of $\ex F$ to a functor $V$ from $\mathcal K$ to the category of complex vectorbundles with morphisms fiberwise-injective  complex-vectorbundle maps so that 
%\begin{itemize}
%\item $V(\hat f)=V_{i}(\hat f)$ if $\hat f$ is locally isomorphic to $\hat f_{i}$,
%\item and if $\hat f\longrightarrow \hat g$ is a morphism in $\mathcal K$, with $\hat f$ and $\hat g$ locally isomorphic to $\hat f_{i}$ and $\hat f_{j}$ repectively, then the map $V(\hat f)\longrightarrow V(\hat g)$ is the inclusion $V_{i}(\hat f)\longrightarrow V_{j}(\hat f)$ followed by the map $V_{j}(\hat f)\longrightarrow V_{j}(\hat g)$;
%\end{itemize}

\item\label{Kcat t}  $\dbar\hat f$ is transverse to $V_{j}(\hat f)\subset V(\hat f)$ when $V_{j}(\hat f)$ is defined and $\dim V_{j}\leq \dim V(\hat f)$; moreover, the intersection of $\dbar\hat f$ with $V_{j}$ is locally isomorphic to $\hat f_{j}$ (and therefore contained in $\mathcal K$).
\item\label{Kcat open} On the other hand, if $ V(\hat f)\subset V_{i}(\hat f)$ --- so $\hat f\in \mathcal U_{i}$ and $\dim V(\hat f)\leq \dim V_{j}(\hat f)$ --- then $\hat f$ is in the substack represented by $\hat f_{i}/G_{i}$ --- so, in $\mathcal K$,  there is a map of a $G_{i}$--fold cover of $\hat f$ to $\hat f_{i}$.

\end{enumerate}

\end{defn}

 \begin{remark} A Kuranishi category, $\mathcal K$, contains more information than a choice of good coordinate system from \cite{KFOOO}. In particular, the vectorbundles $V_{i}$ on open subsets $\mathcal U_{i}$ do not appear there, and the transversality condition \ref{Kcat t} is only required to hold at the intersection of $\dbar$ with $0$ --- without the extensions of our vectorbundles $V_{i}$ from our definition, this condition only makes sense at the intersection of $\dbar$ with $0$. We also explicitly require that $\mathcal K$ is extendable --- I think that this condition may be obtained by shrinking a good coordinate system as in \cite{KFOOshrink}, and  expect the extra data of Definition \ref{Kcat} to be definable from a good coordinate system after shrinking  and making choices using Proposition \ref{global section}.
 \end{remark}
 
\begin{defn}For a Kuranishi category, $\mathcal K$, define $\mathcal K^{hol}\subset\mathcal K^{st}$ to be the substack of $\mathcal K^{st}$ consisting of all holomorphic curves ---  those $f$ in $\hat f$ so that $\dbar\hat f$ vanishes at $f$. Use the induced topology from $\mathcal K^{st}$ on $\mathcal K^{hol}$, so define an open substack of $\mathcal K^{hol}$ to be the intersection of an open substack of of $\mathcal K^{st}$ with $\mathcal K^{hol}$.
\end{defn}

Note that Definition \ref{Kcat} parts \ref{Kcat t} and \ref{Kcat open} imply that the intersection of $\mathcal K^{hol}$ with $\mathcal U_{i}$ is the quotient of $\{\dbar\hat f_{i}=0\}$ by $G_{i}$, so each $\hat f_{i}/G_{i}$ covers an open substack of $\mathcal K^{hol}$.

We need the notion of a Kuranishi category over an exploded manifold or orbifold --- by exploded orbifold, we mean a Deligne-Mumford stack  over the category of exploded manifolds: a stack $\ex Z$ with proper diagonal $\ex Z\longrightarrow \ex Z\times \ex Z$, and  locally represented by $\ex A/G$ where $\ex A$ is an exploded manifold, and $G$ is a finite group acting on $\ex A$; see Remark 2.3 of \cite{evc}.

\begin{defn}\label{K proper}A Kuranishi category over an exploded manifold or orbifold $\ex Z$ is a Kuranishi category $\mathcal K$ along with a submersion $\pi:\mathcal K\longrightarrow \ex Z$.   Say that $\mathcal K$ is proper over $\ex Z$  if $\pi$ restricted to $\mathcal K^{hol}$ is proper --- so given a map $\ex B\longrightarrow \ex Z$ from an exploded manifold $\ex B$ with a compact subset $C$, the stack $\mathcal K^{hol}\times_{\ex Z} C$ is compact.

Say that $\mathcal K$ is complete over $\ex Z$ if it is proper over $\ex Z$ and, for every family $\hat f$ in $\mathcal K$, the map of integral affine spaces $\totb{\ex F(\hat f)}\longrightarrow \totb{\ex Z}$ is complete --- so, the inverse image of any complete subset of $\totb{\ex Z}$ (with its integral affine connection) is a complete subset of $\totb{\ex F(\hat f)}$.

Say that $\mathcal K$ is compact or complete if it is proper or complete respectively over a point. 
\end{defn}

Note that $\mathcal K$ may be complete without $\mathcal K^{hol}$ being complete (in the sense of Definition 3.15 of \cite{iec}). For example consider $\mathcal K$ with a single chart $\hat f$, where $\ex F(\hat f)=\et 1{[0,\infty)}$, $V(\hat f)=\mathbb C\times \et 1{[0,\infty)}$, and $\dbar\hat f(\tilde z)= (\totl{\tilde z},\tilde z)$ (using notation from Example 3.5 of \cite{iec}). Then $\mathcal K^{hol}$ is the stack represented by $\et 1{(0,\infty)}$, which is compact but not complete. On the other hand, $\mathcal K$ is complete by Definition \ref{K proper} because the tropical part of $\et 1{[0,\infty)}$ is the complete polytope $[0,\infty)$. If we perturb $\dbar \hat f$ to be transverse to the zero section, its intersection with the zero section will be complete, so  $\mathcal K^{hol}$ wants to be complete. 

If $\mathcal K$ is the Kuranishi category associated to an embedded Kuranishi structure in some collection of connected components of $\Msw(\hat {\ex B})$, then Lemma 4.2 and Theorem 6.1 of \cite{cem} give conditions under which $\mathcal K^{hol}$ is proper over $\ex B_{0}$. Whenever such a $\mathcal K$ is proper over $\ex B_{0}$, it is also complete over $\ex B_{0}$, because connected families $\hat f$ in $\mathcal K$ either consist of curves with domain $\ex T$, and are parametrized by (a cover of) the quotient of some open subset of $\hat{\ex B}$ by some $\ex T$ action, or  are families of curves in $\hat {\ex B}\longrightarrow \ex B_{0}$ with universal tropical structure. Universal tropical structure is defined in Definition 4.1 and Theorem 3.1 of \cite{uts}; see also Remark 3.3 of \cite{uts}. The $\hat f$ in the Kuranishi charts from \cite{evc} are constructed to have universal tropical structure. In fact,  Theorem 5.3 and Proposition 5.9 of \cite{evc} imply that if  $\hat f/G$ locally represents a closed substack of the moduli stack of curves, (and the curves in $\hat f$ don't have domain $\ex T$), then $\hat f$ must have  universal tropical structure.

\

\subsection{Orienting  Kuranishi charts}

\

We  orient our Kuranishi charts using a canonical homotopy of $D\dbar$ to a complex-linear operator. At holomorphic curves $f$, $T_{f}\dmsw(\ex B)$ has a canonical complex structure, described in section 2.7 of \cite{evc}. In the more general case of curves in a family of targets $\hat {\ex B}\longrightarrow \ex B_{0}$, $T_{f}\dmsw(\hat {\ex B})\ov{\ex B_{0}}$ has a canonical complex structure.\footnote{Given a submersion $\pi:\ex A\longrightarrow \ex B_{0}$, we use the notation $T\ex  A\ov{\ex B_{0}}$ to indicate the vertical tangent bundle of $\ex A$ --- i.e. the kernel of $T\pi$; see  section 2.7.2 of \cite{evc} for a discussion of $T_{f}\dmsw(\hat {\ex B})\ov{\ex B_{0}}$.} In this section, we  discuss a canonical orientation of our Kuranishi charts relative to $\ex B_{0}$. This orientation is only  defined in a neighborhood of the homomorphic curves, so we  construct Gromov--Witten invariants in this neighborhood where the orientation is defined.   

For a Kuranishi chart $(\mathcal U,V,\hat f/G)$,   $V$ is a complex vectorbundle over $\mathcal U$, and  at each holomorphic curve $f$ in $\hat f$, a canonical  linear homotopy of $D\dbar$ to a complex operator is transverse to $V$. As explained in \cite{evc}, this homotopy gives us  a canonical homotopy of $T_{f}\ex F(\hat f)\ov{\ex B_{0}}$ to a complex subspace of $T_{f}\dmsw(\hat {\ex B})\ov{\ex B_{0}}$. In particular, there is a canonical homotopy class of almost complex structures on $T_{f}\ex F(\hat f)\ov{\ex B_{0}}$. Section 8 of \cite{evc}  constructs a  complex structure on $T_{f}\ex F(\hat f^{\sharp}_{i})\ov{\ex B_{0}}$ at all holomorphic curves $f$, so that the following holds:
\begin{itemize}
\item The  complex structure extends to  the vectorbundle $T\ex F(\hat f^{\sharp}_{i})\ov{\ex B_{0}}$ on a neighborhood of the holomorphic curves in $\hat f_{i}^{\sharp}$.
\item Given an open neighborhood $\hat f$ of a holomorphic $f$ in $\hat f^{\sharp}_{i}$ and a map $\hat f\longrightarrow \hat f_{j}^{\sharp}$, the short exact sequence
\[0\longrightarrow T_{f}\ex F(\hat f_{i}^{\sharp})\ov{\ex B_{0}}=T_{f}\ex F(\hat f)\ov{\ex B_{0}}\longrightarrow T_{f}\ex F(\hat f_{j})\ov{\ex B_{0}}\xrightarrow{D\dbar}V_{j}/V_{i}\longrightarrow 0\]
is complex.  
\item The $\mathbb R$--nil vectors in $T_{f}\ex F(\hat f^{\sharp}_{i})\ov{\ex B_{0}}$ (i.e. those that act as the zero derivation on $\mathbb R$--valued functions) are given the canonical complex structure. 
\end{itemize}
Such a complex structure on $T_{f}\ex F(\hat f_{i}^{\sharp})\ov{\ex B_{0}}$ was constructed in Proposition 8.10 of  \cite{evc} by choosing a connection on the inverse image of $V_{i}$ under the homotopy of $D\dbar$ to its complex-linear part. In fact, an appropriate sheaf of such choices  obeys the Patching, Extension and Averaging axioms. 

Does such a complex structure at holomorphic curves $f$ induce an orientation of $T\ex F(\hat f_{i}^{\sharp})\ov{\ex B_{0}}$ elsewhere? To prove this, we construct a global $2$--form, $\alpha$, so that at holomorphic curves $f$,  $\alpha$ is positive on complex planes within $T_{f}\ex F(\hat f_{i}^{\sharp})\ov{\ex B_{0}}$. 

\begin{defn}[Orienting $2$--form]\label{alpha def} Let the sheaf of orienting $2$--forms on  $\hat f$ be the sheaf of  $2$--forms, $\alpha$, on $\ex F(\hat f)$ so that $\alpha$ is positive on holomorphic planes within $T_{f}\ex F(\hat f)\ov{\ex B_{0}}$ for all holomorphic curves $f$ in $\hat f$.

An orienting $2$--form is a global section, $\alpha$, of the sheaf of orienting $2$--forms (over the $K$--category associated to our embedded Kuranishi structure). Say that $\alpha$ is orienting at a curve $f$ if the following holds.
\begin{itemize}
\item Whenever $\hat f_{i}$ contains $f$, some wedge power of $\alpha$ is a volume form on $T_{f}\ex F(\hat f_{i})\ov{\ex B_{0}}$.
\item Whenever $f$ is contained in $\hat f_{i}$ and $\hat f_{j}$, and $V_{i}\subset V_{j}$, the short exact sequence
\[0\longrightarrow T_{f}\ex F(\hat f_{i})\ov{\ex B_{0}}\longrightarrow T_{f}\ex F(\hat f_{j})\ov{\ex B_{0}}\xrightarrow{D\dbar} V_{j}/V_{i}\longrightarrow 0\]
is oriented when the first two terms are given the orientation from $\alpha$, and the last term is oriented by its complex structure. 
\end{itemize}

\end{defn}

To make such a global choice of $\alpha$, use Proposition \ref{global section}. Note that  $2$--forms may be averaged or patched together using a partition of unity, and still satisfy this positivity condition at holomorphic curves. Accordingly, the sheaf of orienting $2$--forms satisfies the Patching and Averaging axioms. Such $2$--forms can also be extended to satisfy the positivity condition, therefore this sheaf  also satisfies the Extension axiom.  By reducing the size of our extensions $\hat f_{i}^{\sharp}$ if necessary, Proposition \ref{global section}  constructs a global section $\alpha$ of the above sheaf of orienting $2$--forms.

Note that $\alpha$ is always orienting at every holomorphic curve. Because the closure of $\hat f_{i}$ is contained in $\hat f_{i}^{\sharp}$,  $\alpha$ is orienting on an open neighborhood of each holomorphic curve in $\hat f_{i}$.

\begin{defn}\label{K orientation}An orientation of a Kuranishi category $\mathcal K$ is,  for all $\hat f$ in $\mathcal K$, an orientation of $\ex F(\hat f)$ so that for all morphisms $\hat f\longrightarrow \hat g$ and  curves $f\in \hat f$, the following short exact sequence is oriented:
\[T_{f}\ex F(\hat f)\longrightarrow T_{f}\ex F(\hat g)\xrightarrow{ D\dbar\hat g}V(\hat g)/V(\hat f)\]

Given a Kuranishi category $\mathcal K$ over $\ex Z$, an orientation of $\mathcal K$ relative to $\ex Z$ is, for all $\hat f$ in $\mathcal K$,  an orientation of $T\ex F(\hat f)\ov{\ex Z}$  so that the following short exact sequence is oriented:
\[T_{f}\ex F(\hat f)\ov{\ex Z}\longrightarrow T_{f}\ex F(\hat g)\ov{\ex Z}\xrightarrow{D \dbar\hat g}V(\hat g)/V(\hat f)\]

\end{defn}

\begin{remark}\label{oriented restriction}Given an embedded Kuranishi structure on $\dmsw(\hat {\ex B})$, we may construct an orienting $2$--form  $\alpha$, and then restrict the associated Kuranishi category to a neighborhood of $\mathcal K^{hol}$  where $\alpha$ is orienting to obtain an extendable Kuranishi categary $\mathcal K$ that is oriented relative to $\ex B_{0}$. 
\end{remark}

\subsection{The sheaf $S$ on a Kuranishi category}

\

We shall need the following information to specify the sheaf, $S$, we use to define Gromov-Witten invariants. 
\begin{defn}[$\mathcal K_{\epsilon}$ and $\mathcal K_{C}$]\label{Kdef} 
Given a Kuranishi category $\mathcal K$, proper and oriented over $\ex Z$, choose continuous functions $\rho_{i}:\mathcal K\longrightarrow [-1,1]$ with the following properties:
\begin{enumerate}
\item \label{Kdef 1} At any holomorphic curve, $\rho_{i}>\frac 12$ for some $i$.
\item\label{Kdef 2} For each $i$, there is some $\mathcal U_{i}$ and associated vectorbundle $V_{i}$ and chart $\hat f_{i}/G_{i}$ from Definition \ref{Kcat} so that the set where $\rho_{i}\geq 0$ is in $\mathcal U_{i}$. 

(We don't require that the map from the indexing set for $\rho_{i}$ to the indexing set for charts on $\mathcal K$ be injective or surjective.)
\item \label{Kdef 3} The subset of $\ex F(\hat f_{i})$ where $\rho_{i}\geq 0$ is compact in the case that $\ex Z$ is a single point, and more generally, the map from this subset to $\ex Z$ is proper.
\item \label{Kdef finite}Over any compact subset of $\ex Z$, there are only finitely many $i$ so that $\rho_{i}$ is somewhere positive.
\end{enumerate}

 \label{K5}Let $\mathcal K_{C}$  be the substack of $\mathcal K^{st}$ comprised of all curves $f$ in some $\hat f_{i}$ where $\rho_{i}(f)\geq \frac 12$ and where $\dbar f\in V_{j}$ wherever any $\rho_{j}>0$.

\

For any $\frac 12>\epsilon>0$, define $\mathcal K_{\epsilon}\subset\mathcal K$ to be the (full) subcategory of $\mathcal K$ consisting of families $\hat f$ so that for some $i$, $\hat f$ is locally isomorphic to an open subfamily of $\hat f_{i}$ where  $\rho_{i}>\epsilon$.
\end{defn}

\begin{remark}\label{KC open}Our virtual moduli space shall be contained in $\mathcal K_{C}$. As well as $\mathcal K_{C}\longrightarrow \ex Z$ being proper, $\mathcal K_{C}$ has the virtue that any family in $\mathcal K_{\epsilon}$ covers an open substack of $\mathcal K_{C}$.
\end{remark}

Note that we can construct functions $\rho_i$ satisfying Definition \ref{Kdef} using Lemma \ref{function construction}: Given any holomorphic curve $f$ in $\hat f_i$, we can choose some $G_i$--equivariant $\rho_i:\ex F(\hat f^\sharp_i)\longrightarrow [-1,1]$ so that the closure of the support of $\rho_i+1$ is in $\mathcal U_i$, and so that $\rho_i(f)>0$ and the subset of $\ex F(\hat f_{i})$ where $\rho_i\geq 0$ is compact. Such a $\rho_i$ can be extended to satisfy condition \ref{Kdef 2} by applying Lemma \ref{function construction} to $\rho_i+1$. Note that condition \ref{Kdef 3} only applies to $\{\rho_i\geq 0\}$ on $\ex F(\hat f_i)$, and not $\ex F(\hat f_j)$, so this extension of $\rho_i$ automatically satisfies condition \ref{Kdef 3}. When $\ex Z$ is compact, a finite collection of such $\rho_i$ will satisfy condition \ref{Kdef 1}, and therefore all conditions of Definition \ref{Kdef}. When $\ex Z$ is not compact, we can choose an exhaustion of $\ex Z$ by compact subsets $\ex Z_k$, and construct each   $\rho_i$ to be  greater than $-1$ only on some $\ex Z_{k+1}\setminus \ex Z_{k-1}$. Then, because the inverse image of  $\ex Z_k$ in $\mathcal K^{hol}$ is compact, there is a collection of such $\rho_i$ satisfying conditions \ref{Kdef 1} and \ref{Kdef finite} so that only finitely many $\rho_i$ are  greater than $-1$ over $\ex Z_{k+1}\setminus \ex Z_{k-1}$. 

Note that although we can construct $\rho_i$ so that $\{\rho_i\geq 0\}\subset \ex F(\hat f_i)$ is always compact, we only require this subset be proper over $\ex Z$ so that Definition \ref{Kdef} is compatible with pullbacks; see  Definition \ref{K pullback}.

We shall define Gromov--Witten invariants using a sheaf of sections of  $V(\hat f)$ over $\hat f$. We need some notion, (condition \ref{1Sdef} of Definition \ref{Sdef} below), of when these sections are `close enough' to the canonical section defined by $\dbar\hat f$. This notion is provided by compatibly choosing metrics on these $V_{i}$ with the property that, where $\abs{\dbar\hat f}\leq 1$, some $\rho_{j}$ from Definition \ref{Kdef} is greater than $\frac 12$. Such a choice ensures that sections that are sufficently close to $\dbar$ have their zero sets contained where these $\rho_{j}>\frac 12$.

\begin{lemma}\label{Kmetric} On $\mathcal K_{\epsilon}$ consider the sheaf $Met$, where  $Met(\hat f)$ is the set of metrics on the vectorbundle $V(\hat f)$ with the property that
\begin{equation}\label{met condition}\text{ on the subset where  }\abs{\dbar\hat f}\leq 1, \text{ some }\rho_{j}>\frac 12\ .\end{equation}

Then $Met$ obeys the Extension, Patching and Averaging axioms, so Proposition \ref{global section} implies that there exists a globally defined metric that is a global section of $Met$.
 \end{lemma}
\pf

As such metrics may be averaged using a partition of unity, $Met$ clearly satisfies the Averaging and Patching axioms. A section of $Met$ locally exists around any curve $f\in \hat f_{i}$, because either $\dbar f>0$ and we can just choose a metric in which $\abs{\dbar f}>1$, or $f$ is holomorphic, and some $\rho_{j}>\frac 12$ around $f$. Now suppose that $V_{i}\subset V_{j}$, and a section of $Met$ on $\hat f_{i}$ has been chosen. Because $\hat f_{i}$ is locally equal to the transversely-cut-out subset of $\hat f_{j}$ where $\dbar \hat f_{j}\in V_{i}$,  a metric on $V_{i}$ over $\hat f_{i}$ can be locally extended to a metric on $V_{j}$ over $\hat f_{j}$. Because our condition (\ref{met condition}) is always satisfied on an open subset, any such extension will locally meet condition (\ref{met condition}), so $Met$ also satisfies the Extension Axiom. 

\stop

\begin{defn}[The sheaf $S$]\label{Sdef} For a choice of $\mathcal K_{\epsilon}$ from Definition \ref{Kdef} and metric from Lemma \ref{Kmetric},  define a sheaf $S$ over $\mathcal K_{\epsilon}$ as follows:
 For $\hat f$ in $\mathcal K_{\epsilon}$, define $S(\hat f)$ to be the set of $\C\infty1$ sections $\nu$ of $ V(\hat f)$ satisfying the following:
\begin{enumerate}
\item\label{0Sdef} On an open neighborhood of where $\rho_{i}\geq 0$ on $\hat f$, $\nu-\dbar\hat f$ is a  section of $V_{i}(\hat f)$. 
\item 
\label{1Sdef} When using the  metric from Lemma \ref{Kmetric}, \[\abs{\dbar-\nu}< 1\] so, in particular,  wherever $\nu=0$, some  $\rho_{i}>\frac 12$.

\item\label{tSdef} $\nu$ is transverse to the $0$-section of  $V(\hat f)$.
\end{enumerate}

\end{defn}

To understand the purpose of item (\ref{0Sdef}) above, note that $V_{i}(\hat f)$ is defined on an open neighborhood of  where $\rho_{i}\geq 0$, but even here, whenever the dimension of $\hat f$ is greater than the dimension of $\hat f_{i}$, $V_{i}(\hat f)\subsetneq V(\hat f)$, and $\dbar\hat f$ is transverse to this sub-bundle.  

 Given any section $\nu$ of $S(\hat f)$, the intersection of $\nu$ with $0$ defines a closed,  $\C\infty1$ exploded submanifold $\nu^{-1}(0)\subset \ex F(\hat f)$.  There is a canonical orientation of $\nu^{-1}(0)$ relative to $\ex Z$ given by the relative orientation of $\hat f$ and the complex orientation of $V(\hat f)$. Intersection with $0$ provides a  natural transformation from $S$ to a sheaf, $E$, defined as follows:
 
 \begin{defn}\label{Edef}Define a sheaf of sets, $E$, on $\mathcal K_{\epsilon}$ as follows:  Let $E(\hat f)$ be the set of  $\C\infty1$ exploded  submanifolds, $X\subset \ex F(\hat f)$, satisfying the following conditions.
  \begin{enumerate}
  \item  $X \subset\ex F(\hat f)$ is closed, and locally defined by the transverse vanishing of some collection of $\C\infty1$ functions.
 \item  $X$ is contained in $\mathcal K_{C}$ (from Definition \ref{Kdef}). 
 \item Wherever $\rho_{i}>0$, $X$ is contained in $\dbar\hat f^{-1}\lrb{V_{i}(\hat f)}$. 
 \item $X$ is oriented relative to $\ex Z$.
\end{enumerate}
Pullbacks in $E$ are naturally defined as inverse images: the pullback of $X$ under $\iota:\hat g\longrightarrow \hat f$ is the inverse image,  $\ex F(\iota)^{-1}(X)$, of $X$ under the induced map $\ex F(\iota):\ex F(\hat g)\longrightarrow \ex F(\hat f)$.
\end{defn}

Is $\ex F(\iota)^{-1}(X)$ really  in $E$? Remark \ref{KC open} and the second condition above ensure that the image of $\ex F(\iota)$ always intersects $X$ in an open subset. As $\ex F(\iota)$ is  locally an isomorphism onto a subset defined by the transverse vanishing of some $\C\infty1$ functions, $\ex F(\iota)^{-1}(X)$ has all the required properties.

We show below that $S$ satisfies the Patching and Extension axioms, but not the Averaging axiom, so global sections of $S$ may not exist. Instead, we will construct  a weighted branched section of $S$ whose intersection with $0$ is a weighted branched section of $E$.

\begin{lemma}\label{S patching} The sheaf $S$ from Definition \ref{Sdef} obeys the Patching Axiom.  

\end{lemma}

\pf

Consider a collection of sections $\nu_{k}$ defined on open subsets $U_{k}$ of $\hat f$. Using a partition of unity $h_{k}$ subordinate to $U_{k}$, these sections $\nu_{k}$ may be averaged to produce a section $\nu:=\sum h_{k}\nu_{k}$ of $V(\hat f)$. Such a section $\nu$ automatically obeys all the conditions to be a section of $S(\hat f)$, except that $\nu$ may not be transverse to $0$.  Our section $\nu$ also obeys the condition required by the Patching Axiom: given any $\hat g\longrightarrow \hat f$, and $\nu'\in S(\hat g)$  agreeing with the pullback of all $\nu_{k}$ (where defined),
 the pullback of $\nu$ is $\nu'$.  Let $X$ be the closure of the image of all $\hat g\longrightarrow \hat f$ with some $\nu'\in S(\hat g)$ satisfying the above condition.

\begin{claim}\label{tpullback}   On $X\subset \hat f$,   $\nu$ is transverse to the $0$--section of $V(\hat f)$.  \end{claim}

  Suppose that $\nu=0$ at $f\in X\subset \hat f$, as otherwise our claim follows trivially. Let $V_{j}$ have the minimal dimension so that $\rho_{j}(f)\geq 0$. On a neighborhood of where $\rho_{j}\geq 0$, $\nu=\dbar$ mod $V_{j}$. Therefore, $\dbar f\in V_{j}$, and mod $V_{j}$,  the derivative of $\nu$ at $f$ is also equal to the derivative of $\dbar$ at $f$.  Definition \ref{Kcat} part \ref{Kcat t} states that  $\dbar\hat f$ at such an $f$ is transverse to $V_{j}$, so we need only check that the derivative of $\nu$ at $f$ surjects onto $V_{j}$.  

Suppose that $U_{k}$ contains $f$. Then, at the image of our $\hat g\longrightarrow \hat f$,  $\nu=\nu_{k}$. It follows that $T\nu=T\nu_{k}$ when restricted to the closure of the image of $T\ex F(\hat g)\longrightarrow T\ex F(\hat f)$. At $f$, this closure must include  $T_{f}\dbar^{-1}(V_{j}(\hat f))$, which equals   $(T_{f}\nu)^{-1}(V_{j}(f))$ and $(T_{f}\nu_{k})^{-1}(V_{j}(\hat f))$. As $\nu_{k}$ is transverse to $0$, $T_{f}\nu_{k}$ restricted to this subspace surjects onto $V_{j}$, therefore $T_{f}\nu$ surjects onto $V_{j}$, so $\nu$ is transverse to $0$ at $f$ as required. This completes the proof of Claim \ref{tpullback}.

\

To complete our proof of Lemma \ref{S patching}, we need the following:
\begin{claim}\label{transverse argument} We can perturb $\nu$ by a small section $v$ of $V(\hat f)$ so that $\nu+v$ is in $S$, and $v$ vanishes on a neighborhood of $X\subset \hat f$.
\end{claim} 

In particular, the section, $\nu +v$, of  $S$  satisfies the conditions required by the Patching axiom. 

To prove Claim \ref{transverse argument}, consider a curve $f\in \hat f\setminus X$. Let $V_{i(f)}$ have the smallest dimension so that  $\rho_{i(f)}(f)\geq 0$. Let $O_f$ be an open neighborhood of $f$ so that
\begin{itemize}\item $V_{i(f)}$ is defined on a neighborhood of the closure of $O_f$,\item  and the closure of $O_f$ does not intersect $X$, and also does not intersect the set where $\rho_j\geq 0$ for any $j$ with $\dim V_j<\dim V_{i(f)}$. 
\end{itemize}
Then $V_{i(f)}$ is defined on $O_f$, and on $O_f$ we can add a small section  of $V_{i(f)}$ to $\nu$ and still satisfy all conditions of Definition \ref{Sdef} apart from possibly (\ref{tSdef}). Condition (\ref{0Sdef}) implies that $\nu$ is equal to $\dbar$ mod $V_{i(f)}$ on $O_f$, so condition (\ref{Kcat t}) of Definition \ref{Kcat} implies that  $\nu$ is transverse to $V_{i(f)}$ on $O_f$, so we can achieve transversality in $O_f$ by adding a small section of $V_i(f)$. We can now complete the proof of Claim \ref{transverse argument} using a standard transversality argument, given below.   

Choose an exhaustion of $\ex F(\hat f)$ by compact subsets, $C_k$. So, $\ex F(\hat f)=\bigcup_kC_k$  where for all $k$,  $C_k$ is compact and contained in the interior of $C_{k+1}$.  Suppose that we have constructed a $v_k$ so that $v_k$ vanishes on a neighborhood of $X$, and $\nu+v_k$ satisfies conditions (\ref{0Sdef}) and (\ref{1Sdef}) of Definition \ref{Sdef}, and is transverse to $0$ on a neighborhood of $C_k$.  Let us construct $v_{k+1}$ satisfying these requirements and equal to $v_k$ when restricted to $C_{k-1}$. Cover  $C_{k+1}$  by the open set where $\nu +v_k$  is transverse to $0$ --- this open set contains $X$ and $C_{k}$ --- and a finite collection of open subsets $O_f$ satisfying the conditions above, and also contained in $C_{k+2}\setminus C_{k}$. Then choose a finite collection of sections $w_1,\dotsc,w_N$ of $V(\hat f)$, each  with support in some $O_f$, and with image in $V_{i(f)}$,   so that the map 
\[\tilde\nu:\mathbb R^N\times \ex F(\hat f)\longrightarrow V(\hat f)\]
\[\text{defined by }\tilde \nu(t_1,\dotsc, t_N,f):=\nu(f)+ v_k(f)+\sum_{i=1}^N t_iw_i(f)\]
is transverse to the zero section of $V$, at least when restricted to an open neighborhood $O$ of  $\mathbb R^N\times C_{k+1}$. To achieve such transversality, it suffices that, restricted to each $O_f$,  the sections $w_1,\dotsc,w_N$ generate $V_{i(f)}$. A finite collection of such sections exist because the closure of $O_f$ is compact and contained in the domain of definition of $V_i$. Then $\tilde \nu^{-1}(0)\cap O$ is a $\C\infty1$ exploded manifold.   The projection $\pi:\tilde\nu^{-1}(0)\cap O\longrightarrow \mathbb R^N$ factors through a smooth map from a manifold on each stratum, so Sard's theorem applies and the critical locus of $\pi$ in $\mathbb R^N$ has measure $0$. Therefore, there exists a regular point of $\pi$,  $(t_1,\dotsc, t_n)\in \mathbb R^N$, arbitrarily close to $0$. For any such regular point, $\nu +v_k+\sum_i t_iw_i$ is transverse to $0$ on some neighborhood of $C_{k+1}$. 

Because the  $w_i$ are compactly supported, Condition (\ref{1Sdef}) of Definition \ref{Sdef} is still satisfied by $\nu +v_k+\sum_i t_iw_i$ so long as $(t_1,\dotsc,t_N)$ is small enough.  Condition (\ref{0Sdef}) holds because each $w_i$ is a section of $V_{i(f)}$ on some $O_f$. Moreover, $v_{k+1}:=v_k+\sum_i t_iw_i$  agrees with $v_k$ on $C_{k-1}$, and vanishes on a neighborhood of $X$ because each $w_i$ is supported in some $O_f$.

The section $v=\lim_{k\to\infty}v_k$ agrees with $v_k$ on $C_{k-1}$ for all $k$. So, $\nu+v$ is transverse to $0$, and is a section of $S$. Moreover, $v$ vanishes on a neighborhood of $X$, so v satisfies the requirements of Claim \ref{transverse argument}. This completes the proof of Claim \ref{transverse argument} and  Lemma \ref{S patching}.

\stop

\begin{lemma}The sheaf $S$ from Definition \ref{Sdef} obeys the Extension Axiom.
\end{lemma}

\pf 

As $S$ obeys the Patching Axiom, we need only verify the local existence of extensions, and the local existence of sections of $S$. 

Sections of $S$ locally exist around any $f\in \hat f$.  Let $V_{i}$ have the smallest dimension so that $\rho_{i}(f)\geq 0$.  If $\dbar f\neq 0$, then in a neighborhood of $f$, $\dbar$ is a section of $S$. If $\dbar f=0$, then $T_{f}\dbar \hat f$ is transverse to $V_{i}$, so $\dbar$ plus some small section of $V_{i}$ will be transverse to the zero section at $f$, and also obey all the other conditions of Definition \ref{Sdef} on a neighborhood of $f$. 

Now suppose that we have a curve $f$ in $\hat g$ and a morphism $\hat g\longrightarrow \hat f$ in $\mathcal K_{\epsilon}$. Again, let $V_{i}$ have the smallest dimension so that $\rho_{i}(f)\geq 0$. Given a section $\nu\in S(\hat g)$, we must construct an extension of $\nu$ to $\nu'$ around $f$ in $\hat f$. As specified by Definitions \ref{K category} and \ref{Kcat}, an open neighborhood of $f$ in $\hat g$ is isomorphic to an open subset of $\hat f_{j} $, and an open neighborhood of $f$ in $\hat f$ is isomorphic to an open subset of $\hat f_{k} $, where $V_{i}\subset V_{j}\subset V_{k}$ on a neighborhood of $f$ in $\hat f_{k}$. Moreover, the morphism $\hat g\longrightarrow \hat f$ is locally an isomorphism onto the transverse intersection of $\dbar\hat f$ with $V_{j}$. Note that $\dbar-\nu$ is a section of $V_{i}\subset V_{j}\subset V_{k}$. We may therefore extend $\dbar-\nu$ to a section $\dbar-\nu'$ of $V_{i}$ on a neighborhood of $f$ within $\hat f$. At $f$, the resulting section $\nu'$ is transverse to $V_{j}$, and transverse to the $0$-section of $V_{j}$ restricted to $(\nu')^{-1}V_{j}$. It is therefore transverse to the zero-section of $V_{k}(\hat f)=V(\hat f)$ at f.  Therefore $\nu'$ is the required local extension of $\nu$ and $S$ obeys the Extension Axiom.

\stop

\begin{remark}\label{extra transversality} We may impose any extra transversality condition satisfied by a generic section of  $S$, and the Patching Axiom and Extension Axiom will still hold, with an identical proof so long the required notion of transversality is independent of extensions. For example, given a submersion $\pi:\mathcal K\longrightarrow \ex X$, we could impose the extra condition that $\pi$ restricted to the intersection of $\nu$ with $0$ must be transverse to a given map of another compact exploded manifold to $\ex X$.
\end{remark}

\

\subsection{Weighted branched sections of sheaves on $K$--categories}
\label{wb section}

\

The sheaf  $S$  from Definition \ref{Sdef} obeys the Patching and Extension axioms, however it does not obey the Averaging Axiom, and $G$--equivariant sections of $S$ may not even exist, so global sections of $S$ can't be constructed.  Instead, we construct weighted branched sections of $S$.  Our approach is essentially  that originally taken by  Fukaya and Ono  in \cite{FO}, however we  use weighted branched sections with a particular branching structure. Our fixed branching structure deals with complications that arise when patching weighted branched sections together. (Even without specifying a particular branching structure, our version of weighted branched sections is subtly different from the definition given by Cieliebak, Mundet i Rivera and Salamon in \cite{wbsubmanifolds} or the  intrinsic definition given by McDuff in \cite{orbifolds}.) 

In the formalism below,  weighted branched sections of a sheaf $S$ over $O$ are labeled by a finite measure space, $(I(O),\mu)$, so each $i\in I(O)$ labels a section,  $\nu(i)$ of $S$ with a weight $\mu(i)$. Given a map $O'\longrightarrow O$, we want to pull back weighted branched sections to $O'$, which requires a measure-preserving map $I(O)\longrightarrow I(O')$. Sometimes this map will send some $i$ and $j$ to the same point, in which case we require that $\nu(i)$ and $\nu(j)$  coincide on a neighborhood of the image of $O'$; our formalism includes an equivalence relation, $\equiv$, on $I(O)$ to keep track of such requirements. 

\begin{defn}\label{wb def}
A weighted branched cover of a $K$--category, $\mathcal K$ (Definition \ref{K category}), is a contravariant functor $I$ with the following domain and target.
\begin{itemize} 
\item  The domain of $I$ is a full subcategory  $\mathcal O_{I}\subset \mathcal K^{st}$ so that the following holds.

\begin{itemize}
\item All $O$ in $\mathcal O_{I}$ are connected.
\item If $O_{1}$ in $\mathcal K^{st}$ is connected and there is a morphism  $\iota:O_{1}\longrightarrow O_{2}$ so that $O_{2}\in \mathcal O_{I}$, then $O_{1}\in\mathcal O_{I}$.
\item For every family $\hat f$ in $\mathcal K$, there is an open cover of $\hat f$  contained in $\mathcal O_{I}$.
\end{itemize}
\item The target of $I$ is the category with objects  triples $(I(O),\mu, \equiv)$
where \begin{itemize}
\item $I(O)$ is a finite set (thought of as a set of sections),
\item $\mu$ is a probability measure on $I(O)$ so that every point in $I(O)$ has positive rational measure,
\item $\equiv$ is an equivalence relation on $I(O)$ (thought of as a warning that our sections are related through needing to be glued somewhere),
\end{itemize}
and morphisms 
\[\iota^{*}:(I(O_{2}),\mu,\equiv)\longrightarrow (I(O_{1}),\mu,\equiv)\]
consisting of measure preserving maps $\iota^{*}:I(O_{2})\longrightarrow I(O_{1})$  so that for any $i,j\in I(O_{2})$, $i\equiv j$ if $\iota^{*}i\equiv \iota^{*}j$.

\end{itemize}

\

We make the additional assumption that, for any curve $f$ and morphism $f\longrightarrow O$ in $\mathcal O_{I}$, there exists a neighborhood $U\subset O$ of $f$ so that $(I(U),\mu,\equiv)\longrightarrow (I(f),\mu,\equiv)$ is an isomorphism. 

Say that $i$ and $j$ in $I(O)$ are separated at $f$ if their images in $I(f)$ are not equivalent.

Say that $I$ has trivial stabilizers if, whenever  $\psi:\hat f\longrightarrow \hat f$ is a nontrivial    automorphism in $\mathcal K\cap \mathcal O_{I}$ that fixes $\totl{f}\in \totl{\ex F(\hat f)}$, then for all $i\in I(\hat f)$, $i$ and $\psi^{*}i$ are separated at $f$. 

\end{defn}

 For example, given a finite group $G$ and a  $G$--fold cover of a manifold, (in other words, a principal $G$--bundle),  a weighted branched cover $I$ could be defined as follows. Let $\mathcal O_{I}$  be the category of maps from connected manifolds $O$ into our space  pulling back our $G$--fold cover to a trivial cover. Then define  $I(O)$  to be the set of sections of this trivial $G$--fold cover over $O$ with the discrete equivalence relation and  the counting measure divided by $\abs G$. The curious reader is invited to think about what goes wrong if we try to define $I(O)$ for disconnected $O$.
 
 We can create another weighted branched cover from the above one by gluing together branches outside some closed set $C$. Define $\mathcal O_{I}$ as above, but now let $I(O)$ have a single point if $O$ does not intersect $C$, and let $I(O)$   be the set of sections over $O$ if $O$ intersects $C$.  Give such an  $I(O)$ the $G$--invariant probability measure; use the discrete equivalence relation if $O$ is contained in the interior of $C$, and use the trivial indiscrete equivalence relation otherwise, so distinct sections of $I(O)$ are separated on the interior of  $C$.
 
 We can build up more complicated examples by taking products of weighted branched covers.

% 
%\begin{defn}A weighted branched section of a sheaf $S$ over $\mathcal K$ is a natural transformation $I\longrightarrow S$, where $I$ is a weighted branched cover of $I$, (and we regard both $I$ and   $S$  as a contravariant functors from $\mathcal O_{I}\subset \mathcal K$ to the category of sets.)
%\end{defn}
%  
% Translated, this says that to specify a weighted branched cover, we need to choose some open cover of $\{\rho_{i}<\epsilon\}\subset \hat f_{i} $ using families $O_{j}$ contained in $\mathcal O_{I}$.  Over $O_{j}$, we have to choose a set of sections of $S$ indexed by $I(O_{j})$, and our choices must be compatible in the sense that given any pair of morphisms $\iota_{1}:O'\longrightarrow O_{j}$ and $\iota_{2}:O'\longrightarrow O_{j'}$, there exists a map $I(O')\longrightarrow S(O)$ so that the following diagram commutes:
% \[\begin{tikzcd}I(O_{j})\dar \rar[swap]{\iota_{1}^{*}} &I(O')\dar&\lar{\iota_{2}^{*}} I(O_{j'})\dar
% \\ S(O_{j}) \rar[swap]{\iota_{1}^{*}} &S(O')&\lar{\iota_{2}^{*}} S(O_{j'})
% \end{tikzcd}\] 

\

If $\hat f$ is in $\mathcal O_{I}\cap\mathcal K$,  define $S^{I}(\hat f)$ to be the set of maps $\nu:I(\hat f)\longrightarrow S(\hat f)$ satisfying the following condition: if two sections, $i$ and $j$ of $I$ are not separated at $f\in \hat f$, then $\nu(i)$ and $\nu(j)$ agree on some  neighborhood of $f$. (Equivalently, for any such $\nu$,  there exists an open cover $\{U\}$ of $\hat f$ so that if $i\equiv j$ in $U$, then $\nu(i)=\nu(j)$ on $U$.) Given any morphism 
$\iota:\hat g\longrightarrow \hat f$ in $\mathcal O_{I}$, the map $\iota^{*}:I(\hat f)\longrightarrow I(\hat g)$ is surjective, and  if $\iota^{*}i=\iota^{*}j$, then $i$ and $j$ are not separated anywhere on the image of $\hat g$, so $\nu(i)$ and $\nu(j)$ agree on the image of $\hat g$. It follows that there is a unique pullback map $\iota^{*}:S^{I}(\hat f)\longrightarrow S^{I}(\hat g)$ fitting into the following commutative diagram:
\[\begin{tikzcd}I(\hat g)\dar{\iota^{*}\nu}&\lar{\iota^{*}}I(\hat f)\dar{\nu}
\\ S(\hat g)&\lar{\iota^{*}} S(\hat f)\end{tikzcd}\] 
Extend the definition of $S^{I}$ to be a sheaf on $\mathcal K$ as follows.

\begin{defn}Given a weighted branched cover $I$ of $\mathcal K$ and  a sheaf $S$ on $\mathcal K$, define the sheaf $S^{I}$ as follows:
\begin{itemize}
\item For a given family $\hat f$ in $\mathcal K$,  let $\mathcal O_{I, \hat f}$ be the category of families in $\mathcal O_{I}\cap\mathcal K$ with a given map to $\hat f$.
\item Define $I_{\hat f}$ as the composition of $I$ with the functor $\mathcal O_{I,\hat f}\longrightarrow \mathcal O_{I}$  forgeting the map to $\hat f$, and define $S_{\hat f}$ as the composition of $S$ with the functor $\mathcal O_{I,\hat f}\longrightarrow \mathcal K$. 
\item Define $S^{I}(\hat f)$ to be the set of natural transformations\footnote{For the purposes of saying what a natural transformation is, consider the codomain of $I_{\hat f}$ and $S_{\hat f}$ to be the category of sets.} $\nu:I_{\hat f}\longrightarrow S_{\hat f}$ so that for any given $O$ in $\mathcal O_{I,\hat f}$, and $i,j\in I(O)$ not  separated at $f$,   $\nu(i)=\nu(j)$ on a neighborhood of $f$.

\item Given $\iota:\hat g\longrightarrow \hat f$, define $\iota^{*}\nu$ to be the pullback of $\nu$ using the obvious functor $\mathcal O_{I,\hat g}\longrightarrow \mathcal O_{I,\hat f}$. 
\end{itemize}
\end{defn}

 If $\hat f$ is in $\mathcal O_{I}\cap \mathcal K$, $S^{I}(\hat f)$ coincides with our easier definition above. This $S^{I}$ is a sheaf because $S$ is a sheaf. Moreover,   we  prove below that $S^{I}$ obeys the Patching, Extension, and Averaging axioms if $S$ obeys the Patching and Extension axioms and $I$  has trivial stabilizers. 

\begin{lemma}\label{patch} If $S$ obeys the Patching Axiom,  then $S^{I}$ does too.
\end{lemma} 

\pf

We must prove the Patching Axiom for $S^{I}$. In particular we must show that given
 an open cover $\{U_{k}\}$ of $\hat f$, and section $\nu_{k}\in S^{I}(U_{k})$ for all $k$, there exists a section $\nu\in S^{I}(\hat f)$ satisfying the following property:
 
 Given any  morphism $\hat g\longrightarrow \hat f$ in $\mathcal K$ and  section $\nu'\in S^{I}(\hat g)$  agreeing with the pullback of $\nu_{k}$ (where defined) for all $k$,
 the pullback of $\nu$ is $\nu'$.

We shall show that the Patching Axiom for $S^{I}$  follows from the following special case:

\begin{claim}\label{two patch} The Patching Axiom holds in the case when the cover $\{U_{k}\}$ has two elements, $\{U_{1},U_{2}\}$,  and $U_{2}$ is in $\mathcal O_{I}$.
\end{claim}

To prove Claim \ref{two patch}, consider a connected component $O$ of $U_{1}\cap U_{2}$. As $O\in \mathcal O_{I}$, $\nu_{i}$ determine maps 
\[\nu_{i}:I(O)\longrightarrow S(O)\ .\] 
 Below, we shall construct $\nu$ on $O$ using the Patching Axiom for $S$ on the individual sections $\nu_{i}(k)$, while ensuring that $\nu$ extends as  $\nu_{1}$ on $U_{1}\setminus U_{2}$, and $\nu_{2}$ on $U_{2}\setminus U_{1}$.  Repeating the construction for every connected component of $U_{1}\cap U_{2}$, we obtain a section $\nu$ of $S^{I}$ on $U_{1}\cup U_{2}$ that agrees with $\nu_{i}$ where only $\nu_{i}$ is defined. To check the Patching Axiom for such a section, it suffices to check the Patching Axiom on $O$ (and the other connected components of $U_{i}\cap U_{j}$) individually.

Identify $I(O)$ with the set $\{1,\dotsc, n\}$. Use the Patching Axiom for $S$ to patch together $\nu_{1}(1)$ and $\nu_{2}(1)$ to create a $\nu(1)\in S(O)$  agreeing with $\nu_{i}(1)$ on a neighborhood of the boundary  of $O$ within $U_{i}$. Such a $\nu(1)$  obeys the requirements of the Patching Axiom for $S^{I}$ on $O$ because it obeys the requirements of the Patching Axiom for $S$: In particular, if $\iota:\hat g\longrightarrow O$  pulls back $\nu_{1}$ and $\nu_{2}$ to $\nu'$ (and $\hat g$ is connected and hence in $\mathcal O_{I}$), then the patching axiom for $S$ implies that $\iota^{*}\nu_i(1)=\nu'(\iota^{*}(1))$.

On some neighborhood $O_{1\equiv 2}$ of the set where the sections $1$ and $2$ of $I(O)$ are not separated, $\nu_{1}(1)=\nu_{1}(2)$ and $\nu_{2}(1)=\nu_{2}(2)$. We can cover $O$ by open subsets $O_{2;1}$, $O_{2;2}$ and $O_{1\equiv 2}$, so that 
\begin{itemize}
\item the closure of $O_{2;i}$  does not intersect the set where $1$ and $2$ are not separated, 
\item and the closure of $O_{2;1}$ within $U_{2}$ does not intersect the boundary of $O$ within $U_{2}$,
\item and the closure of $O_{2;2}$ within $U_{1}$ does not intersect the boundary of $O$ within $U_{1}$. 
\end{itemize}

\includegraphics[scale=0.8]{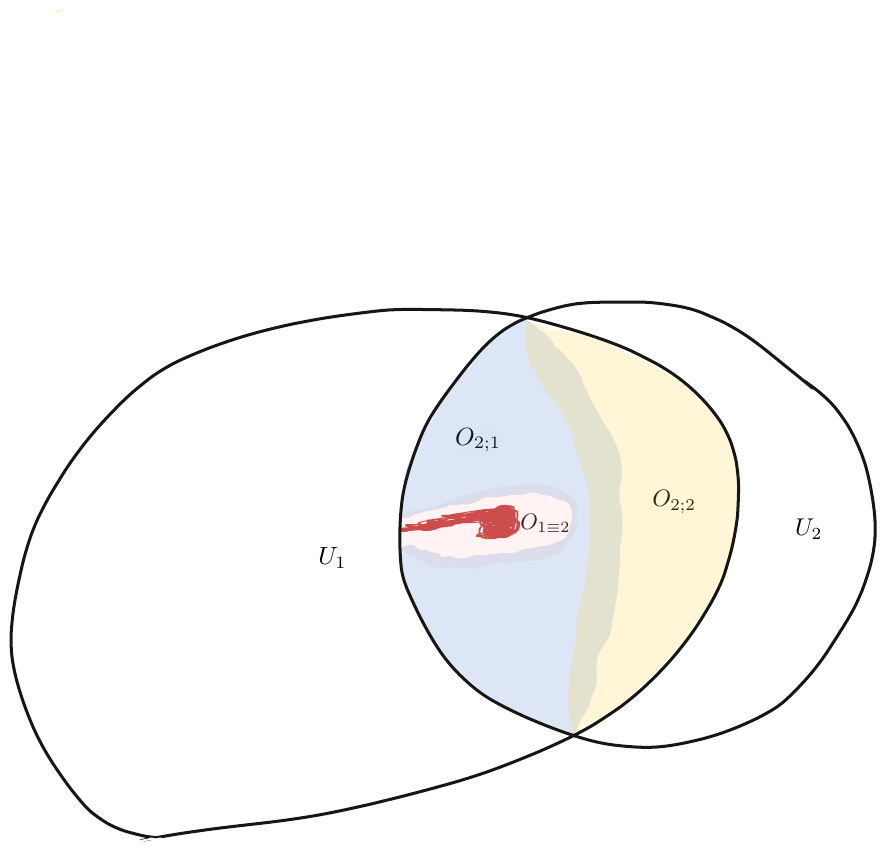}

Such a covering exists  because the set where the sections $1$ and $2$ of $I(O)$ are not separated is closed. 
 Now use the Patching Axiom for $S$ to patch together $\nu_{1}(2)\rvert_{O_{2;1}}$, $\nu_{2}(2)\rvert_{O_{2;2}}$, and $\nu(1)\rvert_{O_{1\equiv 2}}$ to create $\nu(2)$.  Our $\nu(2)$ obeys the requirements of the Patching Axiom on $O$, agrees with $\nu(1)$ on a neighborhood of where $1$ and $2$ are not separated, and agrees with $\nu_{i}(2)$ on a neighborhood of the boundary of $O$ within $U_{i}$.

Inductively continuing this construction gives sections $\nu(k)\in S(O)$, obeying the requirements of the Patching Axiom on $O$, and agreeing with $\nu_{i}(j)$ on a neighborhood of the boundary of $O$ within  $U_{i}$,  so that $\nu(j)=\nu(k)$ on a neighborhood of the set where the sections $j$ and $k$ of $I(O)$ are not separated. In particular, after constructing $\nu(j)$ for all $j<k$, we can cover $O$ by open subsets
\begin{itemize}
\item $O_{j\equiv k}$  for $j<k$, where $\nu_{1}(k)=\nu_{1}(j)$ and $\nu_{2}(k)=\nu_{2}(j)$, covering the set where $j$ and $k$ are not separated;
\item  $O_{k;1}$ with closure in $O$ not intersecting the set where any  $j<k$ is not separated from $k$, and with closure within $U_{2}$ not intersecting the boundary of $O$ within $U_{2}$;
\item and $O_{k;2}$ with closure in $O$ not intersecting the set where any  $j<k$ is not separated from $k$, and with closure within $U_{1}$ not intersecting the boundary of $O$ within $U_{1}$;

\end{itemize}
We can  then construct $\nu(k)$ satisfying the required conditions using the Patching axiom with $\nu(j)\rvert_{O_{j\equiv k}}$ and $\nu_{i}(k)\rvert_{O_{k;i}}$. After completing this construction for all $k\in I(O)$,   we have constructed  a section $\nu\in S^{I}(O)$ agreeing with $\nu_{i}$ on a neighborhood of the boundary of $O$ within $U_{i}$, and  obeying the requirements of the Patching Axiom on $O$: for any map $\iota:\hat g\longrightarrow O$ with connected domain pulling back $\nu_i$ to $\nu'$, $\iota^*\nu(k)=\nu'(\iota^*k)$, so $\iota^*\nu=\nu'$. As $S^{I}$ is a sheaf, there is a section $\nu\in S^{I}(\hat f)$ which agrees with our already constructed $\nu$ in $O$, is similarly defined on every other component of $U_{1}\cap U_{2}$, and which agrees with $\nu_{i}$ elsewhere. As the requirements of the Patching Axiom are local, and such a $\nu$ satisfies the requirements of the Patching Axiom  locally, $\nu$ satisfies the requirements of the Patching Axiom. This completes the proof of Claim \ref{two patch}.
 
 \
 
Now consider an arbitrary open cover of $\hat f$. This arbitrary open cover may be replaced by a countable, locally finite cover $\{U_{i}\}$ consisting of sets in $\mathcal O_{I}$ contained in one of the open sets of the original cover. Then apply Claim \ref{two patch} inductively to patch together  some patched-together section on $\bigcup_{i=1}^{n}U_{i}$ and  $\nu_{n+1}$ on $U_{n+1}$. This procedure constructs a global section $\nu$ of $S^{I}(\hat f)$ obeying the requirements of the Patching Axiom.

\stop  

\begin{lemma}\label{patch extend} If $S$ satisfies the Patching and Extension axioms then $S^{I}$ obeys the Patching and Extension axioms. 
\end{lemma}

\pf

In light of Lemma \ref{patch}, we just need to prove that $S^{I}$ obeys the Extension Axiom. As $S^{I}$ obeys the Patching Axiom, the Extension Axiom follows from the local existence of sections of $S^{I}$, and the local existence of extensions.  

The local existence of sections of $S^{I}$ is easy: For any $O$ in $\mathcal O_{I}$, $S(O)$ is nonempty because $S$ obeys the Extension Axiom.  A constant map $I(O)\longrightarrow S(O)$ suffices to define a section of $S^{I}(O)$. 

Given a morphism $\iota:\hat g\longrightarrow \hat f$ in $\mathcal K$ and $f\in \hat g$, can we  locally extend $\nu\in S^{I}(\hat g)$? This local question can be answered without losing generality by shrinking $\hat g$ and $\hat f$ until  $\hat g$ is in $\mathcal O_{I}$,  $\nu$ agrees on sections  of $I(\hat g)$  not separated at $f$,  and sections of $I(\hat f)$  separated at $f$ are  separated everywhere. Then,  for each section $i$ of $I(\hat f)$, the Extension Axiom for $S$ provides an extension $\nu'(i)$ of $\nu(\iota^{*}i)$. If $i$ and $j$ are not separated at $f$, we can choose $\nu'(i)=\nu'(j)$, because    $\iota^{*}(i)$ and $\iota^{*}(j)$ are not separated at $f$, so $\nu(\iota^{*}i)=\nu(\iota^{*}j)$.   As all other pairs of sections of $I(\hat f)$ are separated everywhere, this $\nu'$ defines a section in $ S^{I}(\hat f)$. 

\stop

\begin{lemma} If $I$ is a weighted branched cover of $\mathcal K$ with trivial stabilizers, and $S$ is a sheaf on $\mathcal K$  satisfying the Patching Axiom, then $S^{I}$ satisfies the Patching and Averaging axioms. 
\end{lemma}

\pf 

Lemma \ref{patch}  implies that $S^{I}$ satisfies the Patching Axiom, so we need only  verify the Averaging Axiom. The idea is to take a given section $\nu$ and patch together pieces of $\nu$  to create a $G$--equivariant section $\nu'$. 

\begin{claim}\label{2average} Suppose that the following holds.
\begin{itemize}
\item $\{U_{1},U_{2}\}$ form a $G$--invariant open cover of $\hat f$.
\item  $\nu_{i}\in S^{I}(U_{i})$.
\item  $\nu_{1}$ is $G$--equivariant.
\item  Each connected component of $U_{2}$ is in $\mathcal O_{I}$.
\item  For each connected component $O$ of $U_{2}$, there exists a subset $I_{0}\subset I(O)$ so that
\begin{itemize}
\item if $i\in I_{0} \subset I(O)$, and $j\equiv i$, then $j\in I_{0}$;
\item the $G$--orbit of $I_{0}$ contains $I(O)$;
\item if the action of some $g\in G$  sends an element of $I_{0}$ to another element of $I_{0}$, then $g$ acts trivially on $O$. 
\end{itemize}

 \end{itemize}
 Then, there exists a $G$--equivariant section $\nu'$ of $S^{I}(\hat f)$, agreeing with $\nu_{1}$ on $U_{1}\setminus U_{2}$ and agreeing with $\nu_{2}$ where $\nu_{2}$ is $G$--equivariant and the same as $\nu_{1}$. More precisely,  given any $G$--equivariant map $\iota:\hat h\longrightarrow \hat f$ with image in $U_{2}$, $\iota^{*}\nu'=\iota^{*}\nu_{2}$ so long as $\iota^{*}\nu_{2}$ is $G$--equivariant and  equals $\iota^{*}\nu_{1}$ on the pullback of $U_{1}$.
\end{claim}

To prove Claim \ref{2average}, we shall first define $\nu'$ on a connected component  $O$ of $U_{2}$.  For  $j\in I_{0}$,  patch together $\nu_{1}(j)$ and $
\nu_{2}(j)$\footnote{Here we have abused notation slightly: $\nu_{2}$ is defined on some open subset of $O$; by $\nu_{2}(j)$ we mean the section which, on connected components, is $\nu_{2}(j')$, where $j'$ is the restriction of $j$.} to obtain $\nu'(j)\in S(O)$ agreeing with $\nu_{1}(j)$ on a neighborhood of the closure of $O\cap U_{1}$ within $U_{1}$.   Similarly to the proof of Claim \ref{two patch},  choose these $\nu'(j)$ so that $\nu'(j)=\nu'(k)$ on a neighborhood of the set where $j$ and $k$ are not separated, and so that the requirements of the Patching Axiom are satisfied by these $\nu'(k)$.

 For any morphism $g$ in $G$, and $j\in I_{0}$ we  can then define $\nu'(g^{*}j)$ to be $g^{*}\nu'(j)$. Whenever $g$ acts nontrivially on $O$, $g^{*}j$ is separated from all $i\in I_{0}$, and the property of being separated is preserved under the action of $G$, so there are no further conditions required of our sections $\nu'(g^{*}j)$, and the resulting $\nu'$ defines a $G$--equivariant section of $S^{I}$ on the orbit of $O$. Because $\nu_{1}$ is equivariant, $\nu'$ agrees with $\nu_{1}$ on an open neighborhood of the closure of the $G$--orbit of $O\cap U_{1}$ within $U_{1}$. Therefore $\nu'$ may be extended to a $G$--equivariant section on the union of the $G$--orbit of $O$ with $U_{1}$ by setting $\nu'=\nu_{1}$ everywhere else.  We may similarly extend the definition of $\nu'$ to all other components of $U_{2}$.
 
 As our $\nu'(j)$ for $j\in I_{0}$ obeyed the requirements of the Patching Axiom, $\nu'$ satisfies the required $G$--equivariant version of the Patching Axiom. This completes the proof of Claim \ref{2average}.
 
 \
 
Construct a $G$--equivariant $\nu'$ inductively using Claim \ref{2average} as follows. Every curve in $\hat f$ has a $G$--invariant open neighborhood $U_{f}$ satisfying the requirements of $U_{2}$ in Claim \ref{2average}: In particular,  $U_{f}$ can be the $G$--orbit of a connected open neighborhood $O$ of $f$ so that the following holds.
 \begin{itemize}
 \item $O$ is invariant under the subgroup $H\leq G$ which is the weak stabilizer of $f$ (consisting of elements  fixing the image of  $f$ in  $\totl{\ex F(\hat f)}$.)
 \item For  $g\in G\setminus H$, $g(O)\cap O=\emptyset$. 
 \item $O$ is in $\mathcal O_{I}$ and is small enough that  $I(O)\longrightarrow I(f)$ is an isomorphism. 
 \end{itemize} 
 For any $i\in I(O)$, $i$ and $g^{*}i$ are separated for any $g\in H$ acting nontrivially on $O$ (because $I$ has trivial stabilizers). As $\equiv$ is a $G$--invariant  equivalence relation on $I(O)$, there is a subset $I_{0}\subset I(O)$ satisfying the requirements of Claim \ref{2average}.
  
  Choose a countable, locally finite open cover $\{O_{i}\}$ of $\hat f$ using sets  satisfying the above conditions on $U_{f}$. Use Claim \ref{2average} with $U_{1}=\emptyset$, $U_{2}=O_{1}$, and $\nu_{2}=\nu$ to construct a $G$--equivariant section $\nu'$ on $O_{1}$. Then inductively apply Claim \ref{2average} with $U_{1}=\bigcup_{i=1}^{n}O_{i}$ and $\nu_{1}$ the already constructed $\nu'$, and $U_{2}=O_{n+1}$ and $\nu_{2}=\nu$. Claim \ref{2average} implies that the resulting patched together $\nu'$ on $\bigcup_{i=1}^{n+1}O_{i}$ is $G$--equivariant and satisfies the condition required by the Averaging Axiom. As the cover $\{ O_{i}\}$ is locally finite, our sequence of $\nu'$s limit to a section $\nu'\in S(\hat f)$ obeying the conditions required by the Averaging Axiom. 

\stop

\section{Construction of $[\mathcal K]$} \label{vmod construction}

In this section we construct a virtual fundamental class $[\mathcal K]$ for a Kuranishi category $\mathcal K$, and show that all choices are cobordant. To state what it means for a choice to be cobordant, we need the notion of a pullback of a Kuranishi category introduced below.

\subsection{Pullbacks of Kuranishi categories}

\

Recall that we are considering Kuranishi categories $\mathcal K$ over $\ex Z$, where $\ex Z$ is an exploded orbifold ---  a Deligne-Mumford stack over the category of exploded manifolds. The following is a standard definition of a representable map of stacks:

\begin{defn} \label{representable} A map \[\ex Z'\longrightarrow \ex Z\] of Deligne-Mumford stacks over the category of exploded manifolds is  representable if, given any submersion $\ex A\longrightarrow \ex Z$ from an exploded manifold $\ex A$, the fiber product $\ex A\times_{\ex Z}\ex Z'$ is  represented by an exploded manifold. 

For more general stacks, say that a submersion $\mathcal X'\longrightarrow\mathcal X$ of stacks is representable if given any map $\ex A\longrightarrow \mathcal X$ from an exploded manifold $\ex A$,  the fiber product $\mathcal X'\times_{\mathcal X}\ex A$ is also represented by an exploded manifold.

\end{defn}

In the language of orbifolds, a representable map is one that is injective on stabilizers. For example, if $\ex Z'$ is an exploded manifold, $\ex Z'\longrightarrow \ex Z$ is always representable. In the definition of pullback of $\mathcal K$ over $\ex Z'\longrightarrow \ex Z$ below, we require that $\ex Z'\longrightarrow \ex Z$ is representable so that $\hat f_{i}\times_{\ex Z}\ex Z'$ is a family of curves parameterized by an exploded manifold, rather than a more general Deligne-Mumford stack. A notion of pullback for non-representable maps $\ex Z'\longrightarrow \ex Z$ exists, but there would either be some extra choices, or we would need arbitrary Deligne-Mumford stacks to take the place of $\hat f/G$.

\begin{defn} \label{K pullback}
Given any Kuranishi category $\mathcal K$ over $\ex Z$ and representable map $\ex Z'\longrightarrow \ex Z$,  the pullback of $\mathcal K$ is a Kuranishi category $\mathcal K'$ over $\ex Z'$.
\[\begin{tikzcd}\mathcal K'\dar\rar&\mathcal K\dar
\\ \ex Z'\rar&\ex Z\end{tikzcd}\]
where\begin{itemize}\item $(\mathcal K')^{st}$ is defined to be $\mathcal K^{st}\times_{\ex Z}\ex Z'$,\item $\hat f'_{i}$ is $\hat f_{i}\times_{\ex Z}\ex Z'$ and $G_{i}'=G_{i}$. (For the purposes of defining $\mathcal K'$, we discard the indices $i$ for which $\hat f_{i}\times_{\ex Z}\ex Z'$ is empty.)
\item $\mathcal U'_{i}$ and $V'_{i}$ are the pullback of $\mathcal U_{i}$ and $V_{i}$ under the map $(\mathcal K')^{st}\longrightarrow \mathcal K^{st}$
\item $V'$ is defined as in Definition \ref{Kcat} part \ref{Kcat V}, so $V'(\hat f)=V_{i}'(\hat f)$ if $\hat f$ is locally isomorphic to $\hat f_{i}'$. In particular, $V'(\hat f\times_{\ex Z}\ex Z')$ is the pullback of $V(\hat f)$.
\item $\dbar'(\hat f\times_{\ex Z}\ex Z')$ is the pullback of $\dbar\hat f$ under the natural map $\hat f\times_{\ex Z}\ex Z'\longrightarrow \hat f$.
\end{itemize}
Similarly, given a map $\mathcal K^{st}\longrightarrow \mathcal X$ and a representable submersion $\mathcal X'\longrightarrow \mathcal X$, define the pullback $\mathcal K'$ of $\mathcal K$ as above with $\mathcal X$ taking the place of $\ex Z$.
\end{defn}

For example, $\mathcal X$ might be the moduli stack of $\C\infty1$ curves in $\ex B$, and $\mathcal X'$ might be the moduli stack of $\C\infty1$ curves in some refinement of $\ex B$, or the moduli stack of $\C\infty1$ curves in $\ex B$ with an extra choice of marked point.

\begin{remark}
Pullback of Kuranishi categories is compatible with many notions: \begin{enumerate}\item
Any extension (Definition \ref{K category}) of $\mathcal K$  also pulls back to an extension of $\mathcal K'$. 
\item An orientation of $\mathcal K$ over $\ex Z$ (Definition \ref{K orientation}) pulls back to an orientation of $\mathcal K'$ over $\ex Z'$. 

Similarly, if the fibers of $\mathcal X'\longrightarrow\mathcal X$ are oriented, an orientation of $\mathcal K$ pulls back to an orientation of $\mathcal K'$.
\item $\mathcal K'$ is proper or complete over $\ex Z'$ (Definition \ref{K proper}) if $\mathcal K$ is proper or complete over $\ex Z$. 

Similarly, if $\mathcal X'\longrightarrow \mathcal X$ is proper or complete, $\mathcal K'$ is proper or complete if $\mathcal K$ is.
\item All choices involved in defining $\mathcal K_{\epsilon}$ (Definition \ref{Kdef}) pull back to define $\mathcal K_{\epsilon}'$, which happens to coincide with the pullback of $\mathcal K_{\epsilon}$ in the above sense. (When we pull back $\mathcal K_{\epsilon}$, we drop the indices $j$ for which the pullback of $\rho_{j}$ is negative.) Also $\mathcal K'_{C}$ is the pulback of $\mathcal K_{C}$.
\item Any sheaf $S$ defined on $\mathcal K^{st}$ pulls back to a sheaf defined on $(\mathcal K')^{st}$, and any global section of such an $S$ over $\mathcal K$ pulls back to a global section over $\mathcal K'$.
\item Any weighted branched cover $I$ of $\mathcal K_{\epsilon}$ (Definition \ref{wb def}) pulls back to a weighted branched cover $I'$ of $\mathcal K'_{\epsilon}$.
\item For $S$ the sheaf from Definition \ref{Sdef}, if   $\ex Z'\longrightarrow \ex Z$ is a submersion, any global section of $S^{I}$ over $\mathcal K_{\epsilon}$ pulls back to a global section of $S'^{I'}$ over $\mathcal K'_{\epsilon}$.

Similarly, given any representable submersion $\mathcal X'\longrightarrow \mathcal X$, any global section of $S^{I}$ over $\mathcal K_{\epsilon}$ pulls back to a global section of $S'^{I'}$ over $\mathcal K'_{\epsilon}$
\end{enumerate}
\end{remark}

\begin{defn}\label{cobordant K}Say that two Kuranishi categories $\mathcal K_{0}$ and $\mathcal K_{1}$  that are proper (or complete) and oriented over $\ex Z$ are cobordant over $\ex Z$ if there exists a Kuranishi category $\mathcal K$ proper (or complete) and oriented over $\ex Z\times \mathbb R$ so that $\mathcal K_{i}$ (with its orientation relative to $\ex Z$) is the pullback of $\mathcal K$ under the inclusion of $\ex Z$ over  $i\subset \mathbb R$. Call $\mathcal K$ a cobordism between $\mathcal K_{i}$.

If $\mathcal K_{i}^{st}$ are substacks of a stack $\mathcal X$ with a map to $\ex Z$, say that $\mathcal K_{i}$ are cobordant within $\mathcal X$ if $\mathcal K^{st}$ is a substack of $\mathcal X\times \mathbb R$ and the map $\mathcal K\longrightarrow \ex Z\times \mathbb R$ is the restriction of the map $\mathcal X\times\mathbb R\longrightarrow \ex Z\times \mathbb R$. 
\end{defn}

\subsection{Construction of $[\mathcal K]$}

\begin{lemma}\label{wb cover}Suppose that a Kuranishi category $\mathcal K$ is proper and oriented over $\ex Z$.  Then there exists a $\mathcal K_{\epsilon}$  and a weighted branched cover $I$ of $\mathcal K_{\epsilon}$ with trivial stabilizers,   satisfying the conditions of definitions \ref{Kdef} and   \ref{wb def}.

Moreover, any two choices of $\mathcal K_{\epsilon}$ and $I$   are cobordant in  the following sense: Suppose that $\mathcal K_{0}$ and $\mathcal K_{1}$ are cobordant (within $\mathcal X$) over $\ex Z$. Then given any two choices of $\mathcal K_{i,\epsilon}$ and weighted branched covers $I_{i}$, there is some cobordism $\mathcal K$ (within $\mathcal X$) between $\mathcal K_{i}$, along with a construction of $\mathcal K_{\epsilon}$ and a weighted branched cover $I$ so that $\mathcal K_{i}$,   $\mathcal K_{i,\epsilon}$ and $I_{i}$ are the pullback of $\mathcal K$,  $\mathcal K_{\epsilon}$, and $I$ respectively under the inclusions of $\ex Z$ over $i\in\mathbb R$. 

\end{lemma}

\pf 

To construct $\mathcal K_{\epsilon}$, we must choose functions $\rho_{i}$ as in  Definition \ref{Kdef}. As discussed after Definition \ref{Kdef}, such functions exist because of  Lemma \ref{function construction}, the fact that Kuranishi categories are extendable,  our assumption that $\mathcal K^{hol}\longrightarrow \ex Z$ is proper, and the fact that each chart $\hat f_{i}/G_{i}$ covers an open substack of $\mathcal K^{hol}$. 

  Definition \ref{K category} part \ref{K category ss} states that that $\hat f_{i}/G_{i}$ represents a substack of $\mathcal K^{st}$. Suppose that $\hat f_{i}$ has the largest dimension of any family in $\mathcal K$ containing $f$. Then some $G$--invariant neighborhood  of $f$ in $\hat f_{i}$ covers an open substack $U$ of $\mathcal K_{\epsilon}^{st}\subset\mathcal K^{st}$. Moreover, because Kuranishi categories are, by definition, extendable, there is a neighborhood $U'$ of $f$ with closure contained in $U$. 

Define a weighted branched cover $I_{f}$ of $\mathcal K_{\epsilon}$ as follows: Define $\mathcal O_{I_{f}}$ to be the full subcategory of $\mathcal K_{\epsilon}^{st}$ with objects all connected families $O$  so that either
\begin{itemize}
\item $O$ does not intersect  the closure of $U'$ --- in this case $I_{f}(O)$ is the probability space with a unique element,
\item or $O$ intersects  the closure of $U'$, is contained entirely inside $ U$, and the corresponding $G$--fold cover $O\times_{\hat f_{i}/G_{i}}\hat f_{i}$  of $O$ defined by the following fiber product diagram is trivial.
\[\begin{tikzcd}O\times_{\hat f_{i}/G_{i}}\hat f_{i}\rar\dar&\hat f_{i}\dar
\\ O\rar &\hat f_{i}/G_{i}\end{tikzcd}\]
 In this case, define $I_{f}(O)$ be the set of sections of  the above $G_{i}$--fold cover  of $O$   with the $G_{i}$--invariant probability measure.
If $O$ is contained in $U'$, then use the discrete equivalence relation on $I_{f}(O)$, otherwise, use the trivial indiscrete equivalence relation on $I_{f}(O)$. 
 
\end{itemize}

Given any morphism $\iota: O'\longrightarrow O$ in $\mathcal O_{I_{f}}$, $\iota^{*}$ is either uniquely determined because $I_{f}(O')$ has a unique point, or $\iota^{*}$ is induced by the natural map between $G_{i}$--fold covers. To check that $I_{f}$ defines a weighted branched cover, we must still check that, given any morphism $g\longrightarrow O$ from a curve in $\mathcal O_{I_{f}}$, there exists a neighourhood $O'$ of the image of $g$ so that the induced map $I_{f}(O')\longrightarrow I_{f}(g)$ is an isomorphism.
If $g$ is in the boundary of $U'$, then the corresponding map $I_{f}(O)\longrightarrow I_{f}(g)$ is an isomorphism; if $g$ is not in the boundary of $U'$, then $I_{f}(O)\longrightarrow I_{f}(g)$ is an isomorphism so long as $O$ is small enough not to intersect the boundary of $U'$. Therefore, $I_{f}$ obeys the final condition we required of weighted branched covers. 

Because $\hat f_{i}/G_{i}$ represents a substack of $\mathcal K^{st}$,     $I_{f}$ is separating at any curve in $U'$.

Now, make  corresponding choices of $I_{f_{i}}$, $ U_{i}$ and $U'_{i}$ so that $\{U'_{i}\}$ is a  cover of $\mathcal K_{\epsilon}$ and so that each $U_{i}$ intersects $U_{j}$ for only finitely many $j$. Then define $I:=\prod_{i} I_{f_{i}}$ as follows:
\begin{itemize}
\item $\mathcal O_{I}:=\bigcap_{i}\mathcal O_{I_{f_{i}}}$. In particular, $O$ is in $\mathcal O_{I}$ if it is connected and is contained in $ U_{i}$ whenever it intersects $U'_{i}$. As each curve has a neighborhood that intersects only finitely many $U'_{i}$,  each family in $\mathcal K_{\epsilon}$ still has an open cover contained in $\mathcal O_{I}$. 
\item For $O\in\mathcal O_{I}$, 
\[I(O):=\prod_{i}I_{f_{i}}(O)\ .\]
As $O$ is in only finitely many $U_{i}$, $I_{f_{i}}$ is a probability space with a unique element for all but finitely many $i$, so $\prod_{i}I_{f_{i}}(O)$ is still a finite probability space.
 \end{itemize}
Sections in $I(O)$ are separated if and only if their image is separated in some $I_{f_{i}}(O)$. Because $I_{f_{i}}$ has trivial stabilizers on $U'_{i}$,  $I$  has trivial stabilizers on $U'_{i}$ for each $i$, so $I$ has trivial stabilizers on all of $\mathcal K_{\epsilon}$.

\

It remains to prove that any two choices of $\mathcal K_{i,\epsilon}$ and $I_{i}$ are cobordant. We may reparameterize the original cobordism between $\mathcal K_{i}$ to obtain a Kuranishi category $\mathcal K$ over $\ex Z\times \mathbb R$ with the extra property that, for $i\in\{0,1\}$, there exist neighborhoods $N_{i}$ of $i\in\mathbb R$ so that the pullback of $\mathcal K$ to $\ex Z\times N_{i}$ is the pullback of $\mathcal K_{i}$ under the projection $\ex Z\times N_{i}\longrightarrow \ex Z$. Using $\ex Z\times N_{i}\longrightarrow \mathcal K_{i}$, pull back the functions used to define $\mathcal K_{i,\epsilon}$. Then, extend these functions by taking the minimum of each function and  a smooth function on $\mathbb R$ equal to 1 at $i$ and $-1$ outside of $N_{i}$. Following this,  use Lemma \ref{function construction} to construct the other functions required to define $\mathcal K_{\epsilon}$, (as discussed following Definition \ref{Kdef}), and  choose these other functions to be  $-1$ on the inverse image of $\{0,1\}$. This construction of $\mathcal K_{\epsilon}$ pulls back to the construction  of $\mathcal K_{i,\epsilon}$ as required.

Now  define $I$. Use the pullbacks of $I_{i}$ over $N_{i}$ to define weighted branched covers $I_{i}'$ as follows:
Choose open neighborhoods  $N_{i}'\subset N_{i}$ of $0$ or $1$ respectively so that  $\bar {N_{i}'}\subset N_{i}$. Let $O$ in $\mathcal K^{st}$ be in $\mathcal O_{I'_{i}}$ in the following two cases:
\begin{itemize}
 \item If $O$ is connected and does not intersect the pullback of  $\bar{N_{i}'}$, define $I'_{i}(O)$ to be the probability space with a unique element.
 \item If $O$ is contained in the inverse image of  $ N_{i}$, and projects to $O'$ within $\mathcal K^{st}_{i}$ so that $O'\in\mathcal O_{I_{i}}$,
 \begin{itemize}
 \item  define $I'_{i}(O)$ to be the probability space with a unique element if $O$ does not intersect the pullback of $\bar{ N_{i}'}$,
 \item define $I'_{i}(O)$ to be  $I_{i}(O')$ with the same measure, but the trivial indiscrete equivalence relation if $O$ intersects the pullback of the boundary of $ N_{i}'$,
 \item define $I'_{i}(O)$ to be  $I_{i}(O')$ with the same measure and equivalence relation if $O$ is contained in the pullback of $ N_{i}'$.
\end{itemize}
\end{itemize}

Now $I'_{0}\times I'_{1}$ is a weighted branched cover of $\mathcal K_{\epsilon}$, and the pull back of  $I'_{0}\times I'_{1}$ under the inclusion of $\ex Z$  over $i$ is $I_{i}$.

$I'_{0}\times I'_{1}$ has trivial stabilizers restricted to the pullback of $N_{i}'$, but may not have trivial stabilizers elsewhere. We can make our weighted branched  cover $I$  have trivial stabilizers   by multiplying $I'_{0}\times I'_{1}$ by other weighted branched covers $I_{f}$ as above. We can achieve this while choosing $I_{f}$ to be trivial when restricted the inverse image of $i\in\mathbb R$, so that $I$ pulls back to give $I_{i}$ as required.

\stop

\begin{remark} Similarly, there exists a weighted branched cover of any extendable $K$--category, and as above, any two such weighted branched covers are cobordant.
\end{remark}

\begin{defn}\label{vdm} Given a Kuranishi category $\mathcal K$, proper and oriented over $\ex Z$, construct the virtual class $[\mathcal K]$ as follows: 
 
\begin{enumerate}
\item Choose $\mathcal K_{\epsilon}$ as in Definition \ref{Kdef}.
\item Choose a separating weighted branched cover $I$ of $\mathcal K_{\epsilon}$ as in Lemma \ref{wb cover}.
\item Consider the sheaf $S^{I}$ of weighted branched sections of $S$ from Definition \ref{Sdef}. We have proved that $S^{I}$ satisfies the Patching, Extension, and Averaging axioms. Choose a global section of $S^{I}$ over $\mathcal K_{\epsilon}$, as allowed by Proposition \ref{global section}. (This may involve increasing $\epsilon$ slightly.)
\item Recall that intersection with $0$ defines a natural transformation $S\longrightarrow E$, where $E$ is the sheaf of oriented sub-families from Definition \ref{Edef}. Intersection with $0$ therefore defines a natural transformation $S^{I}\longrightarrow E^{I}$, so the intersection with $0$ of our weighted branched section of $S$ defines a weighted branched section of $E$. Use the notation $[\mathcal K]$ for such a section. 
\end{enumerate}
\end{defn}

 Lemma \ref{wb cover} along with Proposition \ref{global section} and the fact that $S^{I}$ obeys the Patching, Extension and Averaging axioms imply that any two weighted branched sections of $E$ defined using the above procedure are cobordant.

\section{Representing Gromov--Witten invariants using DeRham cohomology}\label{numerical GW}
\label{de rham section}
\subsection{Differential forms and DeRham cohomologies on stacks and Kuranishi categories}

\

Differential forms on a stack (over the category of smooth manifolds or exploded manifolds) form a sheaf --- to each family $\hat f$, we associate the differential forms on $\ex F(\hat f)$. 
Any notion for differential forms commuting with pullbacks (such as exterior differentiation, wedge products, sums) also makes sense for such differential forms on stacks.  In particular, it is possible to do DeRham cohomology with differential forms on a stack. In the case of a stack (over the category of manifolds) with infinite stabilizers, the resulting cohomology will be smaller than the `correct' cohomology, explained in \cite{behrendCOS}, however I do not know how to imitate the constructions of \cite{behrendCOS} in the infinite-dimensional setting of $\dmsw$.
For exploded manifolds,  we use restricted types of differential forms for defining cohomology theories; see  \cite{dre}.   We need several species of differential forms, all of which are just smooth differential forms  on smooth manifolds.

\begin{defn}[$\Omega^{*}(\ex B)$]\label{Omega def} Let  $\Omega^{k}(\ex B)$ be the space of $\C\infty1$ differential $k$--forms $\theta$ on an exploded manifold $\ex B$ so that for all integral vectors $v$, the differential form $\theta$ vanishes on $v$, and for all maps $f:\et 1{(0,\infty)}\longrightarrow \ex B$, the differential form $\theta$ vanishes on all vectors in the image of $df$. 

Similarly, a differential form on a stack $\mathcal X$ is in $\Omega^{*}(\mathcal X)$ if  it  is in $\Omega^{*}(\ex F(\hat f))$ for all $\hat f$ in $\mathcal X$.

Denote by $\Omega^{k}_{c}(\ex B)\subset \Omega^{k}(\ex B)$ the subspace of forms with complete support.\footnote{A form has complete support if  the set where it is non zero is contained inside a complete subset of $\ex B$ --- in other words, a compact subset with tropical part consisting only of complete polytopes. 

Because completely supported forms do not always pull back to be completely supported, there is no analogous notion of forms with complete support on $\dmsw$. See \cite{behrendCOS} for a version of compactly supported forms on finite-dimensional differential stacks.
} 
 Denote the homology of $(\Omega^{*}(\ex B),d)$ or $(\Omega^{*}(\mathcal X),d)$,  by $H^{*}(\ex B)$ or  $H^{*}(\mathcal X)$ respectively.
\end{defn}

\begin{remark}\label{chern} Given a complex vectorbundle $W$ over an extension $\mathcal K^{\sharp}$ of  $\mathcal K$, we may represent the Chern classes of $W$ in $H^{*}(\mathcal K)$ as closed differential forms.  As unitary metrics and connections may be constructed  as sections of sheaves satisfying the Patching, Extension, and Averaging axioms, they may be constructed on $\mathcal K$ using Proposition \ref{global section}, then we may construct the Chern classes of $V$ over $\mathcal K$ using the Chern-Weil construction. 
\end{remark}

We use $\Omega^{*}$ instead of all $\C\infty1$ differential forms in order to use a version of Stokes' theorem, Theorem 3.4  in \cite{dre}. We shall also wish to use integration along the fiber. For this we shall need the following more general types of differential forms. (For integration along the fiber, see Theorem 9.2 of \cite{dre}.)

\begin{defn}[Refined forms] \label{refined def}
  A refined form $\theta\in\ro^{*}(\ex B)$ is a choice $\theta_{p}\in\bigwedge T^{*}_{p}(\ex B)$ for all $p\in\ex B$ satisfying the following condition: given any point $p\in \ex B$, there exists an open neighborhood $U$ of $p$, a complete, surjective, equidimensional submersion
 \[r:U'\longrightarrow U\]
  and a form $\theta'\in\Omega^{*}(U')$ which is the pullback of $\theta$. In other words,    if $v$ is any vector on $U'$ so that $Tr(v)$ is a vector based at $p$, then 
  \[\theta'(v)=\theta_{p}(Tr(v))\ .\] 
 As refined forms pull back to refined forms, there is an analogous notion of refined forms on any Kuranishi category or stack over the category of $\C\infty1$ exploded manifolds.

 A refined form  $\theta\in\ro^{*}(\ex B)$ is completely supported if there exists some complete subset $V$ of an exploded manifold $\ex C$ with a map $\ex C\longrightarrow \ex B$ so that $\theta_{p}=0$ for all $p$ outside the image of $V$. Use the notation $\ro^{*}_{c}$ for completely supported refined forms. (There is no analogous notion on $\dmsw$.)
 
 Denote the homology of $(\ro^{*}(\ex B),d)$ by $\rh^{*}(\ex B)$ and  $(\ro^{*}_{c}(\ex B),d)$ by $\rh^{*}_{c}(\ex B)$.
 
 \end{defn}

The Poincar\'e dual to a map $\ex C\longrightarrow \ex B$  is correctly viewed as a refined differential form in $\rh^{*}(\ex B)$; see Lemma 9.5 of \cite{dre}. As with all types of forms considered in this paper, refined forms admit the usual operations of wedge product, pullbacks, exterior derivatives and contraction with any $\C\infty1$ vectorfield --- for a discussion of wedge products, see section 9 of \cite{dre}; contraction with any $\C\infty1$ vectorfield is defined because any equidimensional  submersion lifts any such vectorfield uniquely to  $\C\infty1$ vectorfield.

\

In the coming sections,  we will define integration and pushforwards of differential forms using $[\mathcal K]$. If $\mathcal K$ is complete (Definition \ref{K proper}) and contained in the stack $\mathcal X$, then integration over $[\mathcal K]$ defines a map
    
\[H^{*}(\mathcal X)\xrightarrow{\ \ \ \int_{[\mathcal K]} \ \ \  } \mathbb R\] or more generally a map \[\rh^{*}(\mathcal X)\xrightarrow{\ \ \ \int_{[\mathcal K]} \ \ \  } \mathbb R\]
The first map factors through $H^{*}\longrightarrow \rh^{*}$, so the second map contains more information than the first. Given a complex vectorbundle $V$ over $\mathcal X$, we may define more maps $\rh^{*}(\mathcal X)\longrightarrow \mathbb R$ by taking the product with Chern classes of $V$ before integrating. 

 If $\mathcal K$ is complete over $ \ex Z$ and contained in a stack $\mathcal X$,   integration along the fiber of the map $[\mathcal K]\longrightarrow \ex Z$ defines a map 
    \[\rh^{*}(\mathcal X)\longrightarrow \rh^{*}(\ex Z)\ .\]
 
 These maps are compatible with base changes
\[\begin{tikzcd}\mathcal X'\rar \dar&\mathcal X\dar
\\ \ex Z'\rar&\ex Z\end{tikzcd}\]
in the sense that the following diagram commutes: 
\begin{equation}\label{GW base change}\begin{tikzcd}\rh^{*}(\mathcal X') \dar&\lar \rh^{*}(\mathcal X)\dar
\\ \rh^{*}(\ex Z')&\lar\rh^{*}(\ex Z)\end{tikzcd}\end{equation}

In the case that $\mathcal K$ is not complete, we need to use the following, more restrictive, types of differential forms to define invariants:

\begin{defn}[Differential forms generated by functions]\label{rof} A differential form is generated by functions if it is locally equal to a form constructed from $\C\infty1$ functions using the operations of exterior differentiation and wedge products. Use the notation $\rof^{*}\subset\ro^{*}$ for the set of refined forms with differential forms generated by functions playing the role of $\Omega^{*}$ in Definition \ref{refined def}.

Use $\rhf^{*}$ to denote the homology of $(\rof^{*},d)$.
%\[\ro^{*}\xrightarrow{\ \ \ d\ \ \ } \ro^{*-1}\]
%by the image of 
%\[\rof^{*+1}\xrightarrow{\ \ \ d\ \ \ }\ro^{*}\]

\end{defn}

 Examples of differential forms generated by functions are the Poincar\'e dual to a point, the Chern class defined using the Chern Weil construction, and any smooth differential form on a smooth manifold. 
 
 \begin{remark} \label{gbf} Differential forms generated by functions could equivalently be defined as $\C\infty1$ differential forms that vanish on all $\mathbb R$--nil vectors.\footnote{A $\mathbb R$--nil vector is a vector $v$ for which $df(v)=0$ for all differentiable $\mathbb R$--valued functions.} This follows from the proof of Lemma 4.1 in \cite{dre}.
 \end{remark}

In the case that $\mathcal K\subset \mathcal X$ is compact but not complete, integration over $[\mathcal K]$ defines a map  \[\rhf^{*}(\mathcal X)\xrightarrow{\ \ \ \int_{[\mathcal K]} \ \ \  } \mathbb R\ .\] If $\mathcal K$ is complete, then the above map factorises through the map $\rhf^{*}\longrightarrow \rh^{*}$. In the case that $\mathcal K$ is not complete,  these invariants only behave well in families parametrized by smooth (not exploded) manifolds. In particular, if $\hat{\ex B}\longrightarrow \ex B_{0}$ is a family of targets parametrized by a smooth manifold $\ex B_{0}$, then Gromov--Witten invariants  give a map $\rhf^{*}(\dmsw(\hat{\ex B}))\longrightarrow H^{*}(\ex B_{0})$. These Gromov--Witten invariants are invariant under base-changes in the sense that a diagram analogous to diagram (\ref{GW base change}) commutes. 

\

\subsection{Integrating over $[\mathcal K]$  for compact $\mathcal K$}

\

\

Recall that $[\mathcal K]$ is a section of $E^{I}$, where $E$ is defined in Definition \ref{Edef}, and $I$ is a weighted branched cover of $\mathcal K_{\epsilon}$ (Definition \ref{wb def}). So $[\mathcal K]$ is a natural transformation $I\longrightarrow E$. In particular, given any family $O$ in $\mathcal K_{\epsilon}\cap\mathcal O_{I}$, we have a finite probability space $(I(O),\mu)$ and a map from $I(O)$ to the set  of complete subfamilies of $O$ contained in $\mathcal K_{C}$. Each such family $[\mathcal K](i)$ has a canonical orientation relative to $\ex Z$, so in the case that $\ex Z$ a point, or oriented,  $[\mathcal K](i)$ is oriented.

\

We need two different notions of a partition of unity. By a partition of unity subordinate to a given covering   of $\mathcal K^{st}$ by open substacks we  mean a collection of $\C\infty1$ functions $\rho_{i}$ on $\mathcal K^{st}$ with support contained in these open substacks so that $\sum\rho_{i}=1$. Remark \ref{function remark} implies that we may construct such partitions of unity as usual (under the assumption that $\mathcal K$ is extendable). For integrating differential forms over $[\mathcal K]$, we need a different notion of partition of unity on $\mathcal K_{C}$.   When $\mathcal K$ is compact in the sense of Definition \ref{K proper}, $\mathcal K_{C}$, from Definition \ref{Kdef},  is compact and has a finite open cover consisting of families in $\mathcal K_{\epsilon}\cap\mathcal O_{I}$. The partition of unity on $\mathcal K_{C}$ defined below has functions living on such an open cover. This definition is similar to the notion of a partition of unity on an \'etale proper groupoid given in Definition 22 of \cite{behrendCOS}.

\begin{defn}\label{new partition def} A partition of unity on $\mathcal K_{C}$ is a family $O$ in $\mathcal K$, and a  $\C\infty1$ function 
\[r:O\longrightarrow \mathbb R\]
so that the support of $r$ on any connected component $O_{k}$ of $O$ is compact, and so that  for any family $\ex F$ in $\mathcal K_{C}$, the pullback of  $r$ to $O\times_{\mathcal K^{st}}\ex F$ has proper support over $\ex F$, and has pushforward to $\ex F$ equal to $1$.
\[\begin{tikzcd}O\times_{\mathcal K^{st}}\ex F\rar{\pi_{2}}\dar{\pi_{1}}&\ex F\dar
\\ O\rar & \mathcal K^{st}\end{tikzcd}\]
\[(\pi_{2})_{!}\pi_{1}^{*}r=1\]

Say that this partition of unity is compatible with a weighted branched cover $I$ of $\mathcal K_{\epsilon}$ if each connected component $O_{k}$ of $O$ is in $\mathcal O_{I}\cap \mathcal K_{\epsilon}$. If  $[\mathcal K]$ is a natural transformation $[\mathcal K]:I\longrightarrow E$, say that a partition of unity is compatible with $[\mathcal K]$ if it is compatible with $I$. 

\end{defn}

\begin{lemma}\label{partition}Given any weighted branched cover $I$ of $\mathcal K_{\epsilon}$, there exists a partition of unity on $\mathcal K_{C}$ compatible with $I$. Moreover, given any curve $f\in\mathcal K_{C}$, there exists an open neighborhood $N$ of $f$ within $\mathcal K_{C}$ and a partition of unity $r:\coprod_{k}O_{k}\longrightarrow \mathbb R$ satisfying the following: 
\begin{itemize}
\item There is a group $G_{1}$ of automorphisms of $O_{1}$ so that $O_{1}/G_{1}$ represents a substack of $\mathcal K^{st}$.
\item On  the intersection of $O_{1}$ with $N$, $r= 1/\abs{G_{1}}$.
\item On the intersection of $O_{k}$ with $N$ for $k\neq 1$, $r=0$.
\end{itemize}
\end{lemma}

\pf
Choose a locally finite cover of $\mathcal K_{C}$ by families $O_{k}\in\mathcal K_{\epsilon}\cap \mathcal O_{I}$ with automorphism groups $G_{k}$ so that $O_{k}/G_{k}$ represents a substack of $\mathcal K^{st}$. We can do this so that $O_{1}/G_{1}$ compactly contains a neighborhood $N$ of the given curve $f\in\mathcal K_{C}$.

Using Lemma \ref{function construction}, we may choose $\C\infty1$ functions $r_{k}:\mathcal K_{\epsilon}\longrightarrow [0,1]$ satisfying the following conditions:
\begin{itemize}
\item The support of $r_{k}$ intersected with $\mathcal K_{C}$ is compactly contained within $O_{k}/G_{k}$.
\item The support of $r_{k}$ intersects the support of $r_{j}$ for only finitely many $j$.
\item  $\sum_{k}r_{k}>0$ on  $\mathcal K_{C}$.
\item The support of $r_{1}$ contains $N$, and all other $r_{j}$ vanish on $N$.
\end{itemize}

As $\mathcal K_{C}$ is a closed substack, Remark \ref{function remark} implies that there exists a $\C\infty1$ function $R:\mathcal K_{\epsilon}\longrightarrow [0,1]$ with zero set $\mathcal K_{C}$.

Then set 
\[r=\frac {r_{k}}{\abs {G_{k}}(R+\sum_{k}r_{k})}\text{ on }O_{k}.\]
Given a curve $f$ in $\mathcal K_{C}$, the fiber product of $f$ with $O_{k}$ consists of $\abs{G_{k}}$ points if $f$ is also contained in $O_{k}$, and is otherwise empty, so
\[\int_{f\times_{\mathcal K^{st}} O_{k}}r=\frac{r_{k}(f)}{\sum_{k}r_{k}(f)}\]
and
\[\int_{f\times_{\mathcal K^{st}}O}r=\frac {\sum_{k}r_{k}(f)}{\sum_{k}r_{k}(f)}=1\]
as required.

\stop

\

\begin{defn}\label{intdef} Given  $\theta\in \ro^{*} (\mathcal K)$, for $\mathcal K$ oriented and compact (Definition \ref{K proper}), define the integral of $\theta$ over $[\mathcal K]$ as follows:
Choose a partition of unity $r:\coprod_{k}O_{k}\longrightarrow \mathbb R$ compatible with $[\mathcal K]$. Then  define
\[\int_{[\mathcal K]}\theta:=\sum_{k}\sum_{i\in I(O_{k})}\mu(i)\int_{[\mathcal K](i)}r\theta\]
where $[\mathcal K](i)$ is the closed oriented exploded submanifold  of $O_{k}$ that is the image of $i$ under the map $[\mathcal K]:I(O_{k})\longrightarrow E(O_{k})$, and $\mu$ is the probability measure on the indexing set $I(O_{k})$. Note that $r\theta$ has compact support on $[\mathcal K](i)$, so the integral of $r\theta$ is defined as in \cite{dre}. Note also that the above sum is finite. 
\end{defn}

As $\ro^{*}$ contains $\Omega^{*}$ and $\rof^{*}$, the above definition also works for $\theta$ in $\Omega^{*}$ or $\rof$.

 \begin{lemma}\label{iwd} $\int_{[\mathcal K]}\theta$ does not depend on the choice of partition of unity.
 \end{lemma}
 
 \pf Let $r:\coprod_{k}O_{k}\longrightarrow \mathbb R$ and $r':O'\longrightarrow \mathbb R$ be two partitions of unity compatible with $[\mathcal K]$. We have
 \[\int_{[\mathcal K]}\theta:=\sum_{k}\sum_{i\in I(O_{k})}\mu(i)\int_{[\mathcal K](i)}r\theta\]
and because $r'$ is a partition of unity, 
\[\int_{O'\times_{\mathcal K^{st}}[\mathcal K](i)}r'r\theta=\int_{[\mathcal K](i)}r\theta\ .\]
We need to examine $O'\times_{\mathcal K^{st}}[\mathcal K](i)$. Let $O$ be a connected component of $O'\times_{\mathcal K^{st}}O_{k}$, so $O\in\mathcal O_{I}\cap \mathcal K_{\epsilon}$, and comes with a morphism $\iota:O\longrightarrow O_{k}$. For any $j\in \iota^{*}i$,  $[\mathcal K](j)=\iota^{-1}([\mathcal K](i))$, and therefore $[\mathcal K](j)$ is some collection of connected components of $ O'\times_{\mathcal K^{st}}[\mathcal K](i)$. If $\iota^{*}i$ always had a unique element $j$, the union of all such $[\mathcal K](j)$ would be $ O'\times_{\mathcal K^{st}}[\mathcal K](i)$. In the general case, $\iota^{*}i$ is some finite set with measure $\mu(i)$, so we get
\[\mu(i)\int_{[\mathcal K](i)}r\theta=\sum_{j\in \iota^{*}i}\mu(j)\int_{[\mathcal K](j)}r'r\theta\]
 where the sum is over all connected components $O$ of $O'\times_{[\mathcal K^{st}]}O_{k}$ and $j$ in the inverse image of $i$ within $I(O)$. Taking the sum of this expression over all $O_{k}$ and $i\in I(O_{k})$ gives
 \[\int_{[\mathcal K]}\theta=\sum_{O}\sum_{j\in I(O)} \mu(j)\int_{[\mathcal K](j)}r'r\theta\]
where the above sum now is over all connected components $O$ of $O'\times_{\mathcal K^{st}}\coprod_{k}O_{k}$. The above expression is symmetric in the two partitions of unity, therefore the integral is independent of the choice of partition of unity.

\stop

\begin{lemma}\label{complete stokes} If $\theta\in \ro^{*}(\mathcal K)$, and $\mathcal K$ is complete,  then
\[\int_{[\mathcal K]}d\theta=0\ .\]
\end{lemma}

\pf

The fact that $\mathcal K$ is complete implies that the support of each $r$ on each connected component of $ O$ is complete.
Because exterior differentiation and integration are linear we may use a partition of unity on $\mathcal K^{st}$ to reduce to the case that $\theta$ has small support. In particular, we can assume that the support of $\theta$ is small enough that Lemma \ref{partition} gives a partition of unity $r:\coprod_{k}O_{k}\longrightarrow \mathbb R$  so that $r\theta$ has support contained in $O_{1}$ on a subset where $r$ is equal to $1/\abs G$. Then
\[\int_{[\mathcal K]}d\theta=\sum_{i\in I(O_{1})}\frac{\mu(i)}{\abs G}\int_{[\mathcal K](i)}d\theta=0\]
because $\theta\in \ro_{c}^{*}([\mathcal K](i))$, and  such forms satisfy Stokes' theorem; see \cite{dre}.

\stop

\begin{lemma}\label{compact stokes} If $\theta\in \rof^{*}(\mathcal K)$ and $\mathcal K$ is compact but not necessarily complete,    then
\[\int_{[\mathcal K]}d\theta=0\ .\]
\end{lemma}

\pf 

As in the proof of Lemma \ref{complete stokes}, this  lemma reduces to the case of proving that $\int_{[\mathcal K](i)}d\theta=0$ for any compactly supported $\theta\in\rof^{*}[\mathcal K](i)$. By using a partition of unity on $[\mathcal K](i)$, we may reduce to the case that $\theta$ is supported on an open subset $U$ of $[\mathcal K](i)$, and pulls back under a refinement map $U'\longrightarrow U$ to a compactly supported form generated by functions on $U'$. We may use a partition of unity on $U'$ to reduce to the case that $\theta$ is compactly supported within a single coordinate chart $V$ on $U'$. The fact that $\theta$ is generated by functions implies that $\theta$ is pulled back from a differential form on $\mathbb R^{n}$ under an embedding $\totl{V}\longrightarrow \mathbb R^{n}$. Then, the usual Stokes' theorem implies that  $\int_{V}d\theta=0$. 

\stop

\begin{cor} If $\mathcal K$ is complete, then integration over $[\mathcal K]$ defines a map \[\rh^{*}(\mathcal K)\longrightarrow\mathbb R\ .\] If $\mathcal K$ is compact,  integration over $[\mathcal K]$ defines a map 
\[\rhf^{*}(\mathcal K)\longrightarrow \mathbb R\ .\]
\end{cor}

If $\mathcal K$ is contained in a stack $\mathcal X$, we shall show that the resulting maps $\rh^{*}(\mathcal X)\longrightarrow \mathbb R$ and $\rhf(\mathcal X)\longrightarrow \mathbb R$ only depend on the cobordism class of $\mathcal K$ within $\mathcal X$, and in particular, are independent  of the choices involved in the construction of $[\mathcal K]$.

\subsection{Pushing forward cohomology classes }

\

To define pushforwards of differential forms along maps $[\mathcal K]\longrightarrow \ex X$,  we need integration along the fiber, constructed in Theorem 9.2 of  \cite{dre}. Given any oriented submersion of exploded manifolds, $\psi:\ex X\longrightarrow \ex Y$, there exists a linear chain map 
\[\psi_{!}:\ro_{c}^{*}(\ex X)\longrightarrow \ro_{c}^{*'}(\ex Y) \]
uniquely determined by the usual property of integration along the fiber, namely, 
\[\int_{\ex X}\psi^{*}\alpha\wedge \beta=\int_{\ex Y}\alpha\wedge \psi_{!}\beta\ .\]
Above,  the notation $*'$ emphasizes  that $\psi_{!}$ does not preserve degree --- it shifts it by $\dim Y-\dim X$.
Later, in Lemma \ref{tc pushforward}, we show that $\psi_!$ sends $\rof^*(\ex X)$ into $\rof^*(\ex Y)$.
We only need that $\ex X$ is oriented relative to $\ex Y$; when $\ex Y$ is not oriented, the above expression should either be interpreted locally in $\ex Y$ with a choice of orientation, or $\alpha$ needs to be a form twisted by the orientation line bundle of $\ex Y$. For integration along the fiber to work, we require our forms to have  complete support (or at least complete support relative to the target $\ex Y$.)

\begin{remark} 
This pushforward is also defined in the case that $\ex X$ and $\ex Y$ are exploded orbifolds (i.e. Deligne-Mumford stacks in the category of exploded manifolds.)  In particular, we can define $\psi_{!}$ so that given any pullback diagram, 
\[\begin{tikzcd}\ex X'\dar\rar{\psi'} &\ex Y'\dar
\\ \ex X\rar{\psi}&\ex Y\end{tikzcd}\]
where $\ex X'$ and $\ex Y'$ are exploded manifolds,  the following diagram commutes:
\[\begin{tikzcd}\ro^{*}(\ex X')\rar{\psi'_{!}} &\ro^{*}(\ex Y')
\\ \ro^{*}(\ex X)\uar\rar{\psi_{!}}&\ro^{*}(\ex Y)\uar\end{tikzcd}\]

\end{remark}

If all that can be guaranteed is compact (but not complete) support, integration along the fiber will not give a well behaved form on the target $\ex Y$. For example, consider the proper but incomplete map given by the inclusion of $\et 1{(0,1)}$ into $\ex T$. Integrating the constant function $1$ along the fiber of this map gives a discontinuous function on $\ex T$. To get around this problem, we may use differential forms generated by functions, and restrict to the case that the target of our map is a manifold.

\begin{lemma} Given any oriented submersion from an exploded manifold to a smooth manifold,   \[\psi:\ex X\longrightarrow M\]
and a compactly supported form $\beta\in \rof^{*}\ex X$, integration along the fiber of $\psi$ gives a form $\psi_{!}\beta\in\Omega^{*}_{c} M$ uniquely determined by the property that, for any differential form $\alpha$ on $M$ (twisted by the orientation line bundle if necessary), 
\[\int_{\ex X}\psi^{*}\alpha\wedge \beta=\int_{M}\alpha\wedge\psi_{!}\beta\ .\]
As usual, integration along the fiber is a chain map, so $d\psi_{!}\alpha=\psi_{!}d\alpha$.
\end{lemma}

\pf

Forms on smooth manifolds are uniquely determined by their integrals against all other forms, so if there exists a form $\psi_{!}\beta$ satisfying the required property, it is unique. As the required property is linear, we can use a partition of unity to reduce to the case that $\beta$ is compactly contained in a single coordinate chart.  As $\beta$ is a refined form, we may need to refine this coordinate chart and use a further partition of unity to reduce to the case that $\beta$ is an (unrefined) form compactly supported on some standard coordinate chart isomorphic to $\mathbb R^{k}\times \mathbb R^{n}\times \et mP$, where the submersion to $M$ is modeled on the projection to $\mathbb R^{k}$. Integrating along the fiber, and any integral on this coordinate chart, involves a contribution from each $0$--dimensional stratum of $P$, and no other strata contribute. 

For each $0$--dimensional stratum $i$ of $P$, let $P_{i}$ be the union of all strata in $P$ with closure intersecting $i\in P$.  We may construct a tropical completion\footnote{ See section \ref{tropical completion section} for a more thorough discussion of tropical completion.} $\check P_{i}$ of $P_{i}$ as the union of all rays that start at $i\in P$ and intersect $P$ in more than one point. There is a natural inclusion of $\et m{P_{i}}$ into $\et m{\check P_{i}}$, and any compactly supported form $\beta$ on $\mathbb R^{k+n}\times \et m{P_{i}}$  extends uniquely to a completely supported form $\beta_{i}$ in $\mathbb R^{k+n}\times \et m{\check P_{i}}$. As the submersion to $M$ is constant on the $\et mP$ fibers, it also extends uniquely to a submersion $\psi^{i}:\mathbb R^{n+k}\times \et m{\check P_{i}}\longrightarrow M$. As integration along the fiber is defined and satisfies the required properties on completely supported forms, we may define 
\[\psi_{!}\beta:=\sum_{i}\psi^{i}_{!}\beta_{i}\ .\]

Then 
\[\int_{M}\alpha\wedge\psi_{!}\beta=\sum_{i}\int (\psi^{i})^{*}\alpha\wedge\beta_{i}\]
but the integral of $\psi^{*}\alpha\wedge \beta$ has contributions from each $0$ dimensional stratum $i$ of $P$, each equal to the integral of $(\psi^{i})^{*}\alpha\wedge \beta_{i}$, so 
\[\int_{M}\alpha\wedge\psi_{!}\beta=\int \psi^{*}\alpha\wedge\beta \]
as required.  The fact that $\psi_{!}$ commutes with exterior differentiation follows from the analogous fact for $\psi^{i}_{!}$.

\stop

\

Let $\ex A\longrightarrow \ex Z$ be some family of  exploded manifolds or orbifolds over $\ex Z$, oriented relative to $\ex Z$. For example, in the case of a Kuranishi category describing families of curves in $\hat{\ex B}\longrightarrow \ex B_{0}$,  $\ex Z=\ex B_{0}$, and $\ex A$ might be  $\ex B_{0}$, or the $n$--fold fiber product of $\hat{\ex B}$ over $\ex B_{0}$, or the product of $\ex B_{0}$ with  a component of (the explosion of)  Deligne-Mumford space. In this section, we shall consider pushing forward cohomology classes from $\rh^{*}(\mathcal K)$ along a map,\footnote{By a map from $K$-category or Kuranishi category to $\ex A$, we mean a map $\mathcal K^{st}\longrightarrow \ex A$ --- this entails a compatible choice of $\C\infty1$ map $\ex F(\hat f)\longrightarrow \ex A$ for every family $\hat f$ in $\mathcal K^{st}$.} $\pi:\mathcal K\longrightarrow \ex A$, compatible with a complete submersion $\mathcal K\longrightarrow \ex Z$. 
\[\begin{tikzcd}\mathcal K\rar{\pi}\ar{dr} &\ex A\dar
\\ & \ex Z \end{tikzcd}\]

\begin{defn}\label{pushforward}Given  a  map $\pi$ compatible with a submersion as follows, 
\[\begin{tikzcd}\mathcal K\rar{\pi}\ar{dr} &\ex A\dar
\\ & \ex Z \end{tikzcd}\]
in the case that $\mathcal K$ is complete and oriented over $\ex Z$,  define a pushforward map 
\[\pi_{!}:\ro^{*}(\mathcal K)\longrightarrow \ro^{*'}(\ex A)\]
and in the  case where $\mathcal K$ is proper over $\ex Z$  but $\ex A$ is a manifold, define 
\[\pi_{!}:\rof^{*}(\mathcal K)\longrightarrow \Omega^{*'}(\ex A)\ .\]
In either case, $\pi_{!}$ is defined as follows:

\begin{enumerate}
\item Choose an oriented vectorbundle $W\longrightarrow \ex A$ along with a map 
\[x: W\longrightarrow \ex A\]
so that $x$ restricted to the zero section is the identity, and satisfies an extra transversality condition described below.

\item   Choose a partition of unity $r:\coprod_{k}O_{k}\longrightarrow \mathbb R$ compatible with $[\mathcal K]$. 

\item \label{pushforward t}For any $i\in I(O_{k})$, the natural transformation $[\mathcal K]:I\longrightarrow E$ gives a complete subfamily $[\mathcal K] (i)\subset O_{k}$. Use the notation $\pi(i)$   for $\pi$  restricted to $[\mathcal K](i)$.

Define a map $\hat\pi(i):\pi(i)^{*}W\longrightarrow \ex A$ by the composition below.
\[\begin{tikzcd}\pi(i)^{*}W\dar \rar\ar[bend left]{rr}{\hat \pi(i)}& W\dar \rar{x}& \ex A
\\{} [\mathcal K](i) \rar{\pi(i)}& \ex A\end{tikzcd}\]
The extra transversality assumption we require of $x$ is that these maps $\hat\pi(i)$ are submersions. This condition is always satisfied in the case that $x$ is a submersion restricted to each fiber.  

\item Given any $\theta$ in  $\ro^{*}(\pi^{*}W)$  or $\rof^{*}(\pi^{*}W)$ respectively, with compact support on fibers of $\pi^{*}W\longrightarrow \mathcal K$, define 
\[\hat\pi_{!}(\theta):=\sum_{O_{k}}\sum_{i\in I(O_{k})}\mu(i)\hat{\pi}(i)_{!}(r\theta)\] 
 Note that integration along the fiber of $\hat\pi(i)$ requires an orientation relative to $\ex A$. Both $\ex A$ and the family $[\mathcal K](i)$ are canonically oriented relative to $\ex Z$, and  $W$ is an oriented vectorbundle,  so $\hat \pi(i)$ has a canonical relative orientation.

\item Choose a closed form $e\in\Omega^{*}W$ with fiberwise compact support so that $e$ represents the Thom class of $W\longrightarrow \ex A$.
\item For any $\theta$ in $\ro^{*}(\mathcal K)$ or $\rof^{*}(\mathcal K)$ respectively, consider $\theta$ and $e$ as forms on $\pi^{*}W$, then define
\[\pi_{!}(\theta):=\hat \pi_{!}(\theta\wedge e)\ .\]

\end{enumerate}
\end{defn}

Because the Thom form, $e$, is generated by functions,  Lemma \ref{tc pushforward} implies that $\pi_!$ sends forms in $\rof^*(\mathcal K)$ to $\rof^*(\ex A)$.

\begin{lemma}\label{pushforward poui}The map $\hat\pi_{!}$ from Definition \ref{pushforward} is independent of choice of partition of unity. 
\end{lemma}

\pf 

The proof is identical to the proof of Lemma \ref{iwd}, except integration over $[\mathcal K](i)$ is replaced by pushing forward via $\hat \pi(i)_{!}$.

\stop

\begin{lemma}\label{d commute} The map $\hat\pi_{!}$ from Definition \ref{pushforward} commutes with exterior differentiation.

\end{lemma}

\pf

As $d$ and $\hat \pi_{!}$ are linear,  we may assume that the support of $\theta$ is small enough that Lemma \ref{partition} gives a compatible partition of unity $r:\coprod_{k}O_{k}\longrightarrow \mathbb R$ so that $r\theta$  has support where $r=1/\abs G$ in $O_{1}$.
\[\begin{split}\hat\pi_{!}(d\theta)&=\frac 1{\abs {G}}\sum_{i\in I(O_{1})}\mu(i)\hat{\pi}(i)_{!}d\theta
\\ &=d\lrb{\frac 1{\abs {G}}\sum_{i\in I(O_{1})}\mu(i)\hat{\pi}(i)_{!}\theta}
\\ &=d\hat \pi_{!}(\theta)
\end{split}\]

\stop

As an immediate corollary, $\pi_{!}$ induces maps 
\[\pi_{!}:\rh^{*}(\mathcal K)\longrightarrow \rh^{*}(\ex A)\]\[\text{ and }\pi_!:\rhf^*(\mathcal K)\longrightarrow \rhf^*(\ex A)\]
in the case that $\mathcal K$ is complete over $\ex Z$, and a map
\[\pi_{!}:\rhf^{*}(\mathcal K)\longrightarrow H^{*}(\ex A)\]
in the case that $\ex A$ is a smooth manifold, and $\mathcal K$ is proper  over $\ex Z$. Of course, if $\mathcal K$ is contained in a stack $\mathcal X$, we are usually interested in the precomposition of the above maps with  $H^{*}(\mathcal X)\longrightarrow H^{*}(\mathcal K)$.

\begin{lemma}\label{pi} On the level of cohomology, the map 
$\pi_{!}$
does not depend on the choice of $W$, $x$, or $e$ in Definition \ref{pushforward}.
\end{lemma}

\pf

For a given $W$, any two different choices of $e$ differ by $d\alpha$ where $\alpha$ is some fiberwise compactly supported differential form on $W$. The difference between $\pi_{!}(\theta)$ defined using these two choices of $e$ is therefore $\hat \pi_{!}(\theta\wedge d\pi^{*}\alpha)$, which vanishes in cohomology  whenever $\theta$ is closed because of Lemma \ref{d commute}.

Suppose that we have two homotopic choices of $x$ for a  given $W$. We then have a homotopy $\hat \pi_{t}$ between $\hat \pi$ defined using each choice of $x$. If $\theta$ is closed,  $(\hat\pi_{t})_{!}(\theta\wedge \pi^{*}e)$  then gives a closed form on $\ex A\times [0,1]$  restricting to $\ex A\times\{0,1\}$ to give  $\pi_{!}(\theta)$ defined using the two different choices  of $x$. It follows that the cohomology class of $\pi_{!}(\theta)$ is independent of choice of homotopic $x$.

Given two choices $(W_{1},x_{1},e_{1})$ and $(W_{2},x_{2},e_{2})$, consider $W:=W_{1}\oplus W_{2}$ with $x$ defined by the projection to $W_{1}$ followed by $x_{1}$, and $e$ given by the wedge product of the pullback of $e_{1}$ and $e_{2}$ to $W_{1}\oplus W_{2}$. We may then factorize $\hat \pi(i)$ as 
\[\begin{tikzcd}\pi(i)^{*}W_{1}\oplus W_{2}\rar{p(i)} \ar[bend left]{rr}{\hat \pi(i)} & \pi(i)^{*}W_{1} \rar[swap]{\hat \pi_{1}(i)} & \ex A\end{tikzcd}\]
so $\hat \pi(i)_{!}$ factorizes as $\hat\pi_{1}(i)_{!}\circ p(i)_{!}$. The fact that $e_{2}$ is a Thom class implies that
\[p(i)_{!}(p(i)^{*}\theta \wedge e_{2})=\theta\]
so if we abuse notation as in Definition \ref{pushforward} and regard the form $\theta$ on $\mathcal K$ as also living on $\pi^{*}W_{1}$ and $\pi^{*}W_{1}\oplus W_{2}$,  then 
\[\hat\pi(i)_{!}(\theta\wedge e_{1}\wedge e_{2})=\hat\pi_{1}(i)_{!}(\theta\wedge e_{1})\]
It follows that $\pi_{!}(\theta)$ defined using $\hat \pi$ and $\hat \pi_{1}$ is identical. The projection to $W_{1}$ followed by $x_{1}$ is homotopic to the projection to $W_{2}$ followed by $x_{2}$, therefore the same argument identifies the cohomology class of $\pi_{!}(\theta)$ defined using $W_{1}\oplus W_{2}$ with that defined using $W_{2}$.  We have now shown that on the level of cohomology, 
$\pi_{!}$
does not depend on $W$, $x$, or $e$.

\stop

\begin{lemma} If $\mathcal K$ is complete, the following equation holds when $\theta$ is any closed differential form in $\ro^{*}(\ex A)$.
\[\int_{[\mathcal K]}\pi^{*}\theta=\int_{\ex A}\theta\wedge \pi_{!}(1)\]
 In the case that  $\mathcal K$ is compact and  $\ex A$ is a smooth manifold, the above equation holds for   $\theta$  any closed  differential form on $\ex A$.
\end{lemma}

\pf 

Because $e$ represents the Thom class for $W$, 
\[\int_{[\mathcal K](i)}r\pi^{*}\theta=\int_{\pi(i)^{*}W}r\pi^{*}\theta\wedge e\ , \]
where we have abused notation a little to indicate the lift of $\pi^{*}\theta$ to $\pi^{*}W$ simply as $\pi^{*}\theta$. This is not to be confused with the pullback of $\theta$ to $\pi^{*}W$ using $\hat \pi$. The definition of integration along the fiber gives that
\[\int_{\pi(i)^{*}W}\hat \pi^{*}\theta\wedge re=\int _{A}\theta \wedge\hat\pi(i)_{!}(re)\ .\]
Therefore, 
\[\int_{[\mathcal K]}\pi^{*}\theta-\int_{A}\theta\wedge \pi_{!}(1)=\sum_{O_{k}}\sum_{i\in I(O_{k})}\mu(i)\int_{\pi(i)^{*}W}r(\pi^{*}\theta-\hat\pi^{*}\theta)\wedge e\ .
\]
As the vectorbundle map $W\longrightarrow \ex A$ and the map $x:W\longrightarrow \ex A$ are homotopic, there is a differential form $\alpha$ in $\Omega^{*}W$ so that $d\alpha$ is the difference between $x^{*}\theta$ and the lift of $\theta$ using the vectorbundle map $W\longrightarrow \ex A$.

In particular, on $\pi^{*}W$
\[d\pi^{*}\alpha=\hat \pi^{*}\theta-\pi^{*}\theta\ .\]

The same argument as the proof of Lemma \ref{complete stokes} gives that the integral of $d(\pi^{*}\alpha\wedge e)$ over the pullback of $W$ to $[\mathcal K]$ vanishes, so 
\[\sum_{O_{k}}\sum_{i\in I(O_{k})}\mu(i)\int_{\pi(i)^{*}W}r(\pi^{*}\theta-\hat\pi^{*}\theta)\wedge\pi^{*} e=0\ .\]
Therefore we get the required equality:
\[\int_{[\mathcal K]}\pi^{*}\theta=\int_{A}\theta\wedge \pi_{!}(1)\]

\stop

\begin{thm}\label{cobordism theorem}If $\mathcal K_{0}$ and $\mathcal K_{1}$, complete and oriented over $\ex Z$, are cobordant within a stack $\mathcal X$ with  a  map $\pi:\mathcal X\longrightarrow \ex A$, then, given any construction of $[\mathcal K_{j}]$, the two composite maps
\[\rh^{*}(\mathcal X)\xrightarrow{\iota_{j}^{*}} \rh^{*}(\mathcal K_{j})\xrightarrow{\pi_{!}}\rh^{*}\ex A\] 
are equal, and the same holds for the analogous maps
\[\rhf^{*}(\mathcal X)\xrightarrow{\iota_{j}^{*}} \rhf^{*}(\mathcal K_{j})\xrightarrow{\pi_{!}}\rhf^{*}\ex A\ .\] 
Moreover given any complex vectorbundle $W$ over $\mathcal X$, and construction of any characteristic class,  $c(W)$ on $\mathcal K_{j}$ as in Remark \ref{chern},  the maps
\[\theta\mapsto \pi_{!}(c(W))\wedge \iota_{j}^{*}\theta\]
%\[\begin{tikzcd}[row sep=tiny]\rh^{*}(\mathcal X)\rar &\rh^{*}(\ex A)
%\\\theta \rar[|->] &\pi_{!}(c(W)\wedge \iota_{j}^{*}\theta)
%\end{tikzcd}\]
are equal on the level of cohomology for $j=0,1$, in both $\rh^*$ and $\rhf^*$. 

If $\mathcal K_{j}$ are only proper over $\ex Z$, and $\ex A$ is a manifold, the corresponding maps
\[\rhf^{*}(\mathcal X)\longrightarrow H^{*}(\ex A)\]
\end{thm}
are equal.

\pf

This theorem follows from the fact that any two choices of $[\mathcal K_{j}]$ and $c(W)$ are cobordant.

More precisely,  
Lemma \ref{wb cover} gives that  we can choose a cobordism $\mathcal K\subset \mathcal X\times \mathbb R$ between $\mathcal K_{j}$, and construct $\mathcal K_{\epsilon}$ and $I$ on $\mathcal K$  restricting to the given $\mathcal K_{j,\epsilon}, I_{j}$. By reparametrising $\mathbb R$ if necessary, we may also assume that there are connected open neighborhoods $U_{0}$ of $0$ and $U_{1}$ of $1$ so that,
on each of these neighborhoods, $\mathcal K_{j,\epsilon}$ and $I_{j}$ are pullbacks of the respective choices under the maps $\mathcal X\times U_{j}\longrightarrow \mathcal X$.  
 
 Proposition \ref{global section} then implies that by possibly shrinking $U_{i}$, and by increasing $\epsilon$ to $\epsilon'<\frac 12$, we may choose a global section $\nu$ of $S^{I}$ over $\mathcal K_{\epsilon'}$ so that  $\nu$ restricted to $U_{j}$ is the pullback of the corresponding section $\nu_{j}$ under the maps $\mathcal K_{\epsilon'}\rvert_{ U_{j}}\longrightarrow \mathcal K_{j}$.  (Note that increasing $\epsilon'<\frac 12$ does not affect $[\mathcal K]$ at all, because $[\mathcal K]$ has image inside $\mathcal K_{C}\subset \mathcal K_{\epsilon'}$.) Similarly, as the sheaves of unitary metrics and connections on $W$ satisfy the Patching, Extension, and Averaging axioms, we may choose a unitary metric and connection on $W$ over $\mathcal K_{\epsilon'}$ that restricted to $U_{j}$ is the pullback of the corresponding choices used to define $c(W)$. It follows that defining $c(W)\in \rof^{*}\mathcal K_{\epsilon'}$ using the Chern-Weil construction on $\mathcal K_{\epsilon'}$, we obtain a form that, restricted to $U_{j}$, is the pullback of the corresponding forms on $\mathcal K_{j,\epsilon'}$.   
 
  As allowed by Lemma \ref{pi}, we may assume that the auxiliary choices of $(W,x,e)$ from Definition \ref{pushforward} used to define our two choices of  $\pi_{!}$ are the same, and we may also assume that $W=W'\times \mathbb R$ and that $x$ factors through $W\times\mathbb R\longrightarrow W'$ and is a submersion restricted to fibers.  If we pull back these   choices to define
 $\pi'_{!}:\ro^{*}(\mathcal K_{\epsilon'})\longrightarrow \ro^{*}(\ex A\times\mathbb  R)$ as in Definition \ref{pushforward}, the pulled back $x$ may not satisfy the transversality from condition \ref{pushforward t} of Definition \ref{pushforward} outside of $\ex A\times U_{j}$. After possibly shrinking $U_{j}$, we can modify this pulled back $x$ outside of $\ex A\times  U_{j}$ to satisfy this condition (by changing $x$ in the $\mathbb R$ direction so that it becomes a submersion restricted to fibers where necessary).

Restricting attention to $U_{j}\subset \mathbb R$, we may consider $[\mathcal K]$ on families of the form $O\times U_{j}$. Here $I(O\times U_{j})$ may be identified with $I(O)$ from the corresponding construction of $[\mathcal K_{j}]$.  Moreover, for all $i\in I(O\times U_{j})=I(O)$, $[\mathcal K](i)=[\mathcal K_{j}](i)\times U_{j}$.  

It follows that given any $\theta$ in  $\ro^{*}(\mathcal X)$ or $\rof^{*}(\mathcal X)$ respectively and  $1$--form $\theta_{0}$  supported inside $U_{j}\subset \mathbb R$, 
\[\pi'_{!}(\theta\wedge \theta_{0})=\pi_{!}(\theta)\wedge \theta_{0}\]
where  we abuse notation to think of $\theta$ also  on $\mathcal X\times \mathbb R$, and $\theta_{0}$ as also in both $\rof^{1}(\mathcal X\times \mathbb R)$ and $\rof^{1}(\ex A\times \mathbb R)$. The above equality holds for the two different choices of $\pi_{!}$, depending on the choice of $U_{j}$. The fact that $\pi'_{!}$ permutes with $d$ then implies that on the level of cohomology, $\pi_{!}$ does not depend on these choices. With an even more egregious abuse of notation, it also follows that 
\[\pi'_{!}(\theta\wedge c(W)\wedge \theta_{0})=\pi_{!}(\theta\wedge c(W))\wedge \theta_{0}\]
so on the level of cohomology,  the two different maps $\theta\mapsto \pi_{!}(c(W)\wedge \theta)$ are equal for $j=0,1$.

\stop

\begin{prop}\label{sbc} Suppose that $\mathcal K$ is complete over $\ex Z$. Given a representable  submersion $\ex Z'\longrightarrow \ex Z$ (Definition \ref{representable}), let $\mathcal K'$ be the pullback of $\mathcal K$. Let $\pi$ and $\pi'$ be compatible   maps:
\[\begin{tikzcd}\mathcal K'\dar{\pi'}\rar &\mathcal K\dar{\pi} 
\\ \ex A\times_{\ex Z}\ex Z'\rar{y}\dar & \ex A\dar
\\\ex Z'\rar&\ex Z
\end{tikzcd}\]
Then the following diagrams commute.
\[\begin{tikzcd}\rh^{*}(\mathcal K')\dar{\pi'_{!}} &\lar \rh^{*}(\mathcal K)\dar{\pi_{!}} & \rhf^{*}(\mathcal K')\dar{\pi'_{!}} &\lar \rhf^{*}(\mathcal K)\dar{\pi_{!}} 
\\\rh^{*}( \ex A\times_{\ex Z}\ex Z') &\lar{y^{*}} \rh^{*}(\ex A) & \rhf^{*}( \ex A\times_{\ex Z}\ex Z') &\lar{y^{*}} \rhf^{*}(\ex A)
\end{tikzcd}\]

\end{prop}

\pf

 Every step of the construction of $[\mathcal K]$ pulls back under representable submersions. (In general,  sections of $S$ from Definition \ref{Sdef} may not pull back to transverse sections, but this is not a problem in the case that $\ex Z'\longrightarrow \ex Z$ is a submersion.)  Therefore, we may take $[\mathcal K']$ to be the pullback of $[\mathcal K]$.

Explicitly, if $I'$ is the pullback of the weighted branched cover $I$ used to define $[\mathcal K]$,   for any  $O\in \mathcal O_{I}$, there is a map $\iota^{*}:I(O)\longrightarrow I'(O\times_{\ex Z}\ex Z')$, (so long as $O\times_{\ex Z}\ex Z'$ is connected) and given any $i\in I(O)$, 
\[[\mathcal K'](\iota^{*}i)=[\mathcal K](i)\times_{\ex Z}\ex Z'\ .\]
When $O\times_{\ex Z}\ex Z'$ is not connected, the left hand side of the above expression must be replaced with $\coprod_{\iota}[\mathcal K'](\iota^{*}i)$ where the $\iota$  are the maps from each connected component of $O\times_{\ex Z}\ex Z'$ to $O$. For notational conveniance, we shall continue with the case that $O\times_{\ex Z}\ex Z'$ is connected, and simply use $i$ to indicate $\iota^{*}i$.

The choices of $W$, $x$ and $e$ from Definition \ref{pushforward} are more problematic to pull back. Choose an oriented vectorbundle $W$ over $\ex A$ and map $x:W\longrightarrow \ex A$ satisfying the conditions of Definition \ref{pushforward}. Then consider the pullback $y^{*}W$ of $W$ to $\ex A\times_{\ex Z}\ex Z'$. If $\ex Z'\longrightarrow \ex Z$ is not proper, then the map $x$ may not correspond compatibly to a map $x':y^*W\longrightarrow \ex A\times_{\ex Z}\ex Z'$. On the other hand, the fact that $\ex Z'\longrightarrow \ex Z$ is a submersion implies that we may choose $x':y^{*}W\longrightarrow \ex A\times_{\ex Z}\ex Z'$ so that restricted to some neighborhood of the zero section in $y^{*}W$,  the outer circuit of the following diagram commutes. 
\[\begin{tikzcd}y^*W\dar\dar[bend right, swap]{x'} \rar&W\dar\dar[bend left]{x} 
\\  \ex A\times_{\ex Z}\ex Z'\rar{y}& \ex A\end{tikzcd}\]
If we also choose  $x'$ to be a fiber preserving, submersive reparametrization of the pullback of $x$, then $x'$ will still satisfy all conditions of Definition \ref{pushforward}, including the transversality condition from part \ref{pushforward t}.
  Because the outer circuit of the above diagram commutes on a neighborhood of the zero section, there is a natural map of this neighborhood   into a different fiber product, $W\fp xy (A\times_{\ex Z}\ex Z')$ that uses the map $x:W\longrightarrow \ex A$ instead of the vectorbundle projection. Because $x$ and $x'$ coincide with the vectorbundle projections  when restricted to the zero section, restricted to a sufficiently small neighborhood $N$ of the zero section, this natural map is  an isomorphism onto an open subset of $W\fp xy (A\times_{\ex Z}\ex Z')$.

Given any compactly contained open subset $O$ of $A\times_{\ex Z}\ex Z'$, by modifying $x'$ out of the neighborhood $N$ if necessary, we may assume that there exists a neighborhood $U$ of the zero section in $W$ so that $y^{*}U\cap (x')^{-1}O\subset N$, so we can consider $U\fp xy O$ to be an open subset of $N\subset y^{*}W$.

Recall from Definition \ref{pushforward} that $\pi_{!}$ is defined using pushforward along the maps $\hat \pi(i)$ defined using the following diagram.
\[\begin{tikzcd}\pi(i)^{*}W\dar\rar\ar[bend left]{rr}{\hat \pi(i)} & W\dar\rar{x}&\ex A
\\{} [\mathcal K](i)\rar{\pi(i)}&\ex A\end{tikzcd}\]
Similarly define $\hat\pi'(i)$ using the following diagram.
\[\begin{tikzcd}(y\circ\pi'(i))^{*}W\dar\rar\ar[bend left]{rr}{\hat \pi'(i)} & y^{*}W\dar\rar{x'}&\ex A\times_{\ex Z}\ex Z'
\\{}[\mathcal K'](i)\rar{ \pi'(i)}&\ex A\times_{\ex Z}\ex Z'\end{tikzcd}\]
Now, consider the following diagram.
\[\begin{tikzcd}(y\circ \pi'(i))^{*}W\rar(\hat y)\dar\ar[bend right,swap]{dd}{\hat\pi'(i)}& \pi(i)^{*}W\dar\ar[bend left]{dd}{\hat\pi(i)}
\\ y^{*}W\dar{x'}\rar& W\dar{x}
\\ \ex A\times_{\ex Z}\ex Z'\rar{y}&\ex A\end{tikzcd}\]
The top square above is a fiber product diagram, and the bottom square is isomorphic to a fiber product diagram when $W$ is restricted to $U$, $A\times_{\ex Z}\ex Z'$ is restricted to $O$, and $y^{*}W$ is restricted accordingly.  Therefore, the outer circuit of the above diagram is a fiber product diagram when $A\times_{\ex Z}\ex Z'$ is restricted to $O$, $\pi(i)^{*}W$ is restricted to $\pi(i)^{*}U$, and $(y\circ \pi'(i))^{*}W$ is restricted accordingly. 

Lemma 9.3 from \cite{dre}  implies that the following diagram commutes:
\[\begin{tikzcd}\ro(\pi(i)^{*}U\fp {\hat \pi(i)}y O)\dar{\hat{\pi}'(i)_{!}}&\lar{\hat y^{*}} \ro (\pi(i)^{*}U)\dar{\hat\pi(i)_{!}}
\\ \ro(O)&\lar{y^{*}}\ro(\ex A)
\end{tikzcd}\]
Therefore the following diagram obeys the  restricted commutativity condition described below:
\[\begin{tikzcd} \ro((y\circ\pi'(i))^{*}W)\dar{\hat\pi'(i)_{!}}&\lar{\hat y^{*}} \ro(\pi(i)^{*}W)\dar{\hat{\pi}(i)_{!}}
\\ \ro(\ex A\times_{\ex Z}\ex Z')&\lar{y^{*}}\ro(\ex A)
\end{tikzcd}\]
If $\theta\in\ro(\pi(i)^{*}W)$ has support compactly contained in $\pi(i)^{*}U$, then $\hat \pi'(i)_{!}\hat y^{*}\theta$ restricted to $O$ is equal to $y^{*}\hat \pi(i)_{!}\theta$. 

Define $\pi_{!}$ as in Definition \ref{pushforward} using $W$, $x$, and a Thom form $e$ supported in $U$. If we define $\pi'_{!}$ analogously using the pullback of $W$, $x'$, and the pullback of $e$, then  the following diagram commutes. 
\[\begin{tikzcd} & \ro^{*}(\mathcal K_{\epsilon}')\ar[swap]{dl}{\pi'_{!}} &\lar \ro^{*}(\mathcal K_{\epsilon})\dar{\pi_{!}} 
\\\ro^{*}( \ex A\times_{\ex Z}\ex Z')\rar &\ro^{*}(O) &\lar{y^{*}} \ro^{*}(\ex A)
\end{tikzcd}\]
As integration along the fiber and pullbacks preserve the subspace of refined differential forms generated by functions, the same diagram commutes with $\rof^*$ replacing $\ro^*$. 

If we can choose $O\subset \ex A\times_{\ex Z}\ex Z'$ to be a smooth retraction of $\ex A\times_{\ex Z}\ex Z'$ so that the restriction maps $\rh^{*}( \ex A\times_{\ex Z}\ex Z')\longrightarrow \rh^{*}(O)$ and $\rhf^{*}( \ex A\times_{\ex Z}\ex Z')\longrightarrow \rh^{*}(O)$ are isomorphisms, it follows  the required diagrams commute: 
\[\begin{tikzcd}\rh^{*}(\mathcal K'_{\epsilon})\dar{\pi'_{!}} &\lar \rh^{*}(\mathcal K_{\epsilon})\dar{\pi_{!}} & \rhf^{*}(\mathcal K'_{\epsilon})\dar{\pi'_{!}} &\lar \rhf^{*}(\mathcal K_{\epsilon})\dar{\pi_{!}} 
\\\rh^{*}( \ex A\times_{\ex Z}\ex Z') &\lar{y^{*}} \rh^{*}(\ex A) & \rhf^{*}( \ex A\times_{\ex Z}\ex Z') &\lar{y^{*}} \rhf^{*}(\ex A)
\end{tikzcd}\]

More generally, $\rh^{*}(\ex A\times_{\ex A}\ex Z')$ is the inverse limit of the system   $\rh^{*}(O)$ for all open $O\subset \ex A\times_{\ex A}\ex Z'$  with compact closure, and the same holds for $\rhf^*$. (It is obvious that $\ro^{*}(\ex A\times_{\ex A}\ex Z')$ is the inverse limit of $\ro^{*}(O)$, and the system $\ro^{*}(O)$ satisfies the Mittag-Leffler  condition, so taking homology is compatible with inverse limits; the same holds for refined forms generated by functions.) It follows that the above diagrams always commute.

\stop

\begin{thm}\label{pullback theorem}
Suppose that $\mathcal K$ is complete over $\ex Z$. Let $\mathcal K'$ be the pullback of $\mathcal K$ over a representable map $\ex Z'\longrightarrow \ex Z$ (Definition \ref{representable}). Given compatible   maps,
\[\begin{tikzcd}\mathcal K'\dar{\pi'}\rar &\mathcal K\dar{\pi} 
\\ \ex A\times_{\ex Z}\ex Z'\rar{y}\dar & \ex A\dar
\\\ex Z'\rar&\ex Z
\end{tikzcd}\]
the following diagrams commute:
\[\begin{tikzcd}\rh^{*}(\mathcal K')\dar{\pi'_{!}} &\lar \rh^{*}(\mathcal K)\dar{\pi_{!}} & \rhf^{*}(\mathcal K')\dar{\pi'_{!}} &\lar \rhf^{*}(\mathcal K)\dar{\pi_{!}} 
\\\rh^{*}( \ex A\times_{\ex Z}\ex Z') &\lar{y^{*}} \rh^{*}(\ex A) & \rhf^{*}( \ex A\times_{\ex Z}\ex Z') &\lar{y^{*}} \rhf^{*}(\ex A)
\end{tikzcd}\]

\end{thm}

\pf 

Extend $\ex Z'\longrightarrow \ex Z$ to a  map  $\ex X\longrightarrow \ex Z$ where \begin{itemize}
\item $\ex X$ is a vectorbundle over $\ex Z'$, 
\item restricted to some tubular neighborhood $O_{1}$ of the zero section, the map $\ex X\longrightarrow \ex Z$ factors through $\ex X\longrightarrow \ex Z'\longrightarrow \ex Z$, 
\item restricted to some  other open set $O_{2}\subset \ex X$, $\ex X\longrightarrow \ex Z$ is a submersion with nonempty convex fibers.
\end{itemize}
Because $\ex X$ is a vectorbundle over $\ex Z'$, and $\ex Z'\longrightarrow \ex Z$ is representable,  our extension $\ex X\longrightarrow \ex Z$ is also representable.  Use $\mathcal K_{\ex X}$ and  $\mathcal K_{ O_{i}}$ to denote the respective pullbacks of $\mathcal K$. Proposition \ref{sbc} implies that all inner loops in the following diagram commute, 
\begin{equation}\label{cd1}\begin{tikzcd}\rh^{*}(\mathcal K_{O_{1}})\dar &\lar \rh^{*}(\mathcal K_{\ex X})\rar \dar&\rh^{*}(\mathcal K_{ O_{2}})\dar &\lar \rh^{*}(\mathcal K)\dar\ar[bend right]{ll}\ar[bend right]{lll} 
\\\rh^{*}(O_{1}\times_{\ex Z}\ex A)&\lar \rh^{*}(\ex X\times_{\ex Z}\ex A)\rar &\rh^{*}(O_{2}\times_{\ex Z}\ex A)&\lar \rh^{*}\ex A\ar[bend left]{ll}\ar[bend left]{lll}
\end{tikzcd}\end{equation}\
and the same holds with $\rhf^*$ replacing $\rh^*$.
Because the restriction maps $\rh^{*}(\ex X\times_{\ex Z}\ex A)\longrightarrow \rh^{*}(O_{i}\times_{\ex Z}\ex A)$ are isomorphisms,  every loop in the above diagram commutes, and the same holds with $\rhf^*$ replacing $\rh^*$.  Proposition \ref{sbc} also implies  the commutativity of the inner square of following diagram
\begin{equation}\label{cd}\begin{tikzcd}\rh^{*}(\mathcal K_{ O_{1}})\dar\ar[bend left]{r}&\lar \rh^{*}(\mathcal K')\dar
\\ \rh^{*}(O_{1}\times_{\ex Z}\ex A)\rar[bend right]&\lar \rh^{*}(\ex Z'\times_{\ex Z}\ex A)\end{tikzcd}\end{equation}
and the analogous diagram with $\rhf^*$ replacing $\rh^*$.
The horizontal maps in the above diagram are defined using the projection $O_{1}\longrightarrow \ex Z'$ and inclusion $\ex Z'\longrightarrow O_{1}$. We shall now show that each pair of horizontal maps in the diagram above are inverse isomorphisms, proving that the above diagram commutes. In each case,  the leftward map followed by the rightward is the identity on forms, so it remains to prove that the other  composition is the identity on the level of cohomology.  

Let $\Psi_{t}:O_{1}\longrightarrow O_{1}$ for $t\in[0,1]$ be any fiber-preserving smooth homotopy that is the identity at $t=0$ and the projection onto $\ex Z'\subset O_{1}$ at $t=1$. Then, given a form $\theta\in \ro^{*}O_{1}$, define 
\[K \theta=\int_{0}^{1}\Psi_{t}^{*}i_{\frac{\partial \Psi_{t}}{\partial t}}\theta dt\ .\]  
Then 
\[(dK+Kd)\theta=\Psi_{1}^{*}\theta-\Psi_{0}^{*}\theta\ .\]
It follows that if $\theta$ is closed, it represents the same cohomology class as $\Psi_{1}^{*}\theta$. As the homotopy preserves the fibers of the map to $\ex Z$, it induces a smooth homotopy from the identity on $O_{1}\times_{\ex Z}\ex A$ to the projection to  $\ex Z'\times_{\ex Z}\ex A$. Making the above argument with this new homotopy gives that   the composition of the lower pair of horizontal maps in  the above diagram is the identity on the level of cohomology. The same holds with $\rhf^*$ replacing $\rh^*$ because $K$ preserves the subspace of forms generated by functions. It remains to show the same for the upper pair of maps.

Let us consider the homotopy induced by $\psi_{t}$ on $\mathcal K_{O_{1}}$, and define an operator analogous to the above $K$. Start with a family of curves $\hat f$ in $\mathcal K_{ O_{1}}$. The map $\mathcal K_{O_{1}}^{st}\longrightarrow (\mathcal K')^{st}$ applied to  $\hat f$  gives a family of curves $\hat f_{0}$ with $\ex F(\hat f_{0})=\ex F(\hat f)$. We may then pull this family $\hat f_{0}$ back to get a bigger family of curves $\hat f_{0}\times_{\ex Z}O_{1}$ in $\mathcal K_{ O_{1}}^{st}$ with domain $\ex F(\hat f)\times _{\ex Z'}O_{1}$. Alternately, using the inclusion $\hat{\ex B}'\subset \hat O_{1}$, we can consider $\hat f_{0}$ as a family of curves in $\mathcal K_{O_{1}}^{st}$. As the homotopy $\Psi_{t}:O_{1}\longrightarrow O_{1}$ preserves fibers of the map $O_{1}\longrightarrow \ex Z'$, it induces a homotopy of $\ex F(\hat f_{0}\times_{\ex Z'}O_{1})$ that is the identity at $t=0$, and the projection onto $\ex F(\hat f_{0})$ at $t=1$. Again call this homotopy $\Psi_{t}$, and define $K:\ro^{*}\ex F(\hat f_{0}\times_{ \ex Z'}O_{1})\longrightarrow \ro^{*-1}\ex F(\hat f_{0}\times_{ \ex Z'}O_{1})$ as in the above equation. 
 Then, given  any $\theta\in\rh^{*}(\mathcal K_{\hat O_{1}})$, we can define $K\theta$ on $\ex F(\hat f)$ by defining $K\theta$ on $\ex F(\hat f_{0}\times_{\ex Z'}O_{1})$ as above, then pulling back this $K\theta$ to $\ex F(\hat f)$ via the natural inclusion $\hat f\longrightarrow \hat f_{0}\times_{\ex Z'}O_{1}$. Such a $K\theta$ gives a well defined form in $\ro^{*}(\mathcal K_{O_{1}})$ because given any map  of curves $\hat g\longrightarrow \hat f$, the corresponding diagram 
 \[\begin{tikzcd} \ex F(\hat g_{0}\times_{ \ex Z'}O_{1})\rar \dar{\Psi_{t}}& \ex F(\hat f_{0}\times_{ \ex Z'} O_{1})\dar{\Psi_{t}}
 \\ \ex F(\hat g_{0}\times_{ \ex Z'}O_{1})\rar&\ex F(\hat f_{0}\times_{ \ex Z'}O_{1})\end{tikzcd}\]
commutes.  $dK\theta$ is the difference between $\theta$ and the composition of the topmost pair of  maps in  diagram (\ref{cd}) applied to $\theta$, so this topmost pair of maps are inverse isomorphisms on the level of cohomology. As all our maps preserve $\rof^*\subset \ro^*$, the same holds with $\rhf^*$ replacing $\rh^*$.
  
  As the horizontal maps in diagram (\ref{cd}) are inverse isomorphisms, diagram (\ref{cd}) commutes, and combines with the outermost loop of diagram (\ref{cd1}) to give our required commutative diagram:
   \[\begin{tikzcd}\rh^{*}(\mathcal K')\dar{\pi'_{!}} &\lar \rh^{*}(\mathcal K)\dar{\pi_{!}} 
\\\rh^{*}( \ex A\times_{\ex Z}\ex Z') &\lar{y^{*}} \rh^{*}(\ex A)
\end{tikzcd}\]
  The same argument applies to the analogous diagram replacing $\rh^*$ with $\rhf^*$.
  
  \stop

  \section{Weak  products of Kuranishi categories}
 \label{weak product section}
 
 \

In order to discuss gluing theorems for Gromov--Witten invariants, it is useful to have a notion of fiber products of Kuranishi structures or Kuranishi categories. The tropical gluing formula for Gromov--Witten invariants, equation (1) of \cite{gfgw}, is proved by expressing a Kuranishi category associated to some tropical curve $\gamma$ as a (weak) fiber product of Kuranishi categories associated to the vertices, $v$, of $\gamma$. We can decompose such a (weak) fiber product as a pullback, dealt with in Theorem \ref{pullback theorem}, and a (weak) product, defined below  and dealt with in Theorem \ref{weak product theorem}.

Because of our choice that all Kuranishi charts in an embedded Kuranishi structure should be compatible, and the related part \ref{K category t} of Definition \ref{K category}, and part \ref{Kcat open} of Definition \ref{Kcat},  the product of Kuranishi categories is no longer a Kuranishi category, because the obvious candidates for charts on the product will no longer be compatible.\footnote{Most other authors work with a version of Kuranishi structures that avoid this problem. Our embedded Kuranishi structure could be regarded as akin to a choice of good coordinate charts in Fukaya, Oh, Ohta and Ono's approach, \cite{KFOOO}.} 

Nevertheless, we can construct a weak product of Kuranishi categories by taking the product of charts, and then shrinking these charts appropriately to make them compatible. 

\begin{defn}\label{weak product}Say that $\mathcal K$ is a weak product of some finite collection $\mathcal K_{v}$ of Kuranishi categories if the following holds:
\begin{enumerate}
\item $\mathcal K^{st}\subset \prod_{v}\mathcal K_{v}^{st}$
\item $\mathcal K^{hol}=\prod_{v}\mathcal K^{hol}$
\item For each chart $(\mathcal U_{i},V_{i},\hat f_{i}/G_{i})$ on $\mathcal K$, there are corresponding charts $(\mathcal U_{i_{v}},V_{i_{v}},\hat f_{i_{v}}/G_{i_{v}})$ on $\mathcal K_{v}$ so that
\begin{itemize}
\item $\mathcal U_{i}\subset\prod_{v}\mathcal U_{i_{v}}$;
\item on $\mathcal U_{i}$, $V_{i}=\oplus_{v}V_{i_{v}}$;
\item $G_{i}=\prod_{v}G_{i_{v}}$;
\item $\hat f_{i}$ is a $G_{i}$--equivariant  open subfamily of $\prod_{v}\hat f_{i_{v}}$;
\item $\dbar\hat f_{i}$ is equal to $\prod_{v}\dbar\hat f_{i_{v}}$ restricted to $\ex F(\hat f_{i})\subset \prod_{v}\ex F(\hat f_{i_{v}})$.
\end{itemize}

\end{enumerate}
%Suppose that $\prod_{v}\mathcal K_{v}\longrightarrow \ex Z$ is a Kuranishi category over an orbifold $\ex Z$, and  is a weak product of Kuranishi categories $\mathcal K_{v}$. Given a representable map $\ex Z'\longrightarrow \ex Z$ (Definition \ref{representable}), we may construct the pullback $\mathcal K'$ of $\prod_{v}\mathcal K_{v}$ as in Definition \ref{K pullback}.
%\[\begin{tikzcd}\mathcal K'\dar\rar&\ex Z'\dar
%\\ \prod_{v}\mathcal K_{v}\rar &\ex Z\end{tikzcd}\]   
%Say that such a $\mathcal K$ is a weak fiber product of $\mathcal K_{v}$ over $\ex Z$. \marginpar{meh... maybe remove def of weak fiber product?}
\end{defn}

Note that if $\mathcal K_{v}$ are proper or complete and  oriented over $\ex Z_{v}$, then $\mathcal K$ is proper or complete and oriented over $\prod_{v}\ex Z_{v}$. Given   maps $\mathcal K_{v}\longrightarrow \ex A_{v}$, there is an obvious induced  map $\mathcal K\longrightarrow \prod_{v}\ex A_{v}$.

\begin{thm}\label{weak product theorem}Suppose that $\mathcal K$ is a weak product of some finite collection of complete, oriented Kuranishi categories  $\mathcal K_{v}$ with  maps $\pi_{v}:\mathcal K_{v}\longrightarrow \ex A_{v}$, let $\pi:\mathcal K\longrightarrow \prod_{v}\ex A_{v}$ be the induced  map. Then the following diagrams commute:

\[\begin{tikzcd}\rh^{*}(\mathcal K)\rar{\pi_{!}}&\rh^{*}(\prod_{v}\ex A_{v}) & \rhf^{*}(\mathcal K)\rar{\pi_{!}}&\rhf^{*}(\prod_{v}\ex A_{v})
\\ \prod_{v}\rh^{*}(\mathcal K_{v})\uar\ar{ur}[swap]{\prod_{v}(\pi_{v})_{!}}&&  \prod_{v}\rhf^{*}(\mathcal K_{v})\uar\ar{ur}[swap]{\prod_{v}(\pi_{v})_{!}}\end{tikzcd}\]
\end{thm}

\pf

Using Lemma \ref{function construction}, and the fact that $\mathcal K^{hol}=\prod_{v}\mathcal K^{hol}_{v}$, we may construct functions $\rho_{i,v}:\mathcal K_{v}\longrightarrow \mathbb R$ satisfying the requirements of Definition \ref{Kdef} so that the corresponding functions $\rho_{i}:=\min_{v}\rho_{i,v}$ on $\mathcal K$ also satisfy the requirements of Definition \ref{Kdef}, where the chart $\hat f_{i}/G_{i}$ associated with $\rho_{i}$  is locally the product of the charts $\hat f_{i_{v}}/G_{i_{v}}$ associated with $\rho_{i,v}$. Indeed, each holomorphic curve $f=\prod_v f_v$ in $\hat f_i\subset \prod_v\hat f_{i_v}$ has   some ($G_i$--invariant) product neighborhood $\prod U_v\subset \hat f_i$; we can choose $G_{i_v}$--invariant functions $\rho_{f,v}:\prod\hat f_{i_v}\longrightarrow [-1,1]$ with $\rho_{f,v}(f_v)>\frac 12$ so that $\rho_{f,v}+1$ is compactly supported within $U_v$, then use Lemma \ref{function construction} on $\rho_{f,v}+1$ to extend these functions to functions satisfying Definition \ref{Kdef}  part \ref{Kdef 2}. The corresponding continuous function $\rho_f:=\min_v\rho_{f,v}$ is also greater than $\frac 12$ on $f$, and satisfies Definition \ref{Kdef} parts \ref{Kdef 2} and \ref{Kdef 3}.  Our functions $\rho_i$  and $\rho_{i,v}$ will be chosen from such functions  $\rho_f$ and $\rho_{f,v}$ for holomorphic $f$. As we have assumed that  $\mathcal K_{v}$ are complete, $\mathcal K^{hol}$ is covered by the sets where $\rho_i>0$ for a finite number of such $\rho_i$ so we can achieve  Definition \ref{Kdef} parts \ref{Kdef 1} and \ref{Kdef finite} with this finite collection of functions,  $\rho_i$ and $\rho_{i,v}$.

 The corresponding $\mathcal K_{\epsilon}^{st}$ is a substack of $\prod_{v}\mathcal K_{v,\epsilon}^{st}$, because the subset of $\hat f_{i}$ where $\rho_{i}>\epsilon$ is equal to the product of the subsets of $\hat f_{i_{v}}$ where $\rho_{i,v}>\epsilon$, (but $\prod_{v}\mathcal K_{v,\epsilon}$ contains more families than $\mathcal K_{\epsilon}$ because it also includes the products of the relevant subsets of $\hat f_{i_{v}}$ where $i$ may depend on $v$).
 
  Choose separating  weighted branched covers $I_{v}$ of $\mathcal K_{v,\epsilon}$, and let $I$ be the weighted branched cover of $\mathcal K_{\epsilon}$ that is the product of the pull back of these $I_{v}$ via the maps $\mathcal K_{\epsilon}^{st}\longrightarrow \mathcal K_{v,\epsilon}^{st}$. The fact that $I_{v}$ are separating implies that $I$ is also separating. 

We must also choose metrics on $V$ as in Lemma \ref{Kmetric}. We need to do this so that `small enough' perturbations of $\dbar$ on $\mathcal K_{v,\epsilon}$ pull back to `small enough' perturbations of $\dbar$ on $\mathcal K_{\epsilon}$. The required condition on $\mathcal K_{\epsilon}$ is that wherever $\abs{\dbar\hat f_{i}}\leq 1$, some $\rho_{j}>\frac 12$, and the corresponding condition on $\mathcal K_{v,\epsilon}$ is that wherever $\abs{\dbar \hat f_{i_{v}}}\leq 1$, some $\rho_{j_{v}}>\frac 12$. 
If there are $n$ indices $v$, we need to strengthen these latter conditions so that wherever $\abs{\dbar\hat f_{i_{v}}}\leq n$ for all $v$, there exists some $j$ so that $\rho_{j_{v}}>\frac 12$ for all $v$.
Essentially, we need that  $\{\abs{\dbar\hat f_{i_{v}}}\leq n\}$ is contained in some open neighborhood $U_{v}$ of $\mathcal K^{hol}_{v}$ and that $\prod_{v}U_{v}$ is covered by the sets where $\rho_{j_{v}}>\frac 12$ for all $v$. This is possible because $\prod_{v}\mathcal K^{hol}_{v}$ is compact and  and the sets where all $\rho_{j_{v}}>\frac 12$ form an open cover of $\prod_{v}\mathcal K^{hol}_{v}\subset \prod_{v}\mathcal K_{v}^{st}$.
%  More explicitly, we can construct $U_{v}$ as follows: Identify the indices $v$ with $\{1,\dots,n\}$, let the set of indices $\{j\}$ of $\rho_{j}$ be $J_{1}$ and consider the set of all possible nested subsets $J_{1}\supset J_{2}\supset J_{3}\dotsb \supset J_{n}\supset J_{n+1}$  so that for each $J_{k}$ with $k\leq n$, the sets where $\rho_{j_{v}}>\frac 12$ for all $v\geq k$ form an open cover of $\prod_{v\geq k}\mathcal K_{v}^{hol}$ indexed by $j\in J_{k}$. Then on $\mathcal K_{v}$, the set \[U_{v}:=\bigcap_{J_{v}}\lrb{\bigcup_{J_{v+1}\subset J_{v}}\lrb{\bigcap_{j\in J_{v+1}}\{\rho_{j_{v}}>1/2\}}}\]
% is an open neighborhood of $\mathcal K_{v}^{hol}$. Moreover, $\prod_{v}U_{v}$ is an open neighborhood of $\prod_{v}\mathcal K_{v}^{hol}$ that is contained in the union of $\prod_{v}\{\rho_{j_{v}}>1/2\}$. 
Once such $U_{v}$ have been chosen, the proof of Lemma \ref{Kmetric} applies to show that there exists a global choice of metric on $V_{v}$ over $\mathcal K_{v,\epsilon}$ so that  $\abs{\dbar \hat f_{i_{v}}}\leq n$ is contained in $U_{v}$. 
 
 If we use $1/n$ times the product of these metrics on $V$ over $\mathcal K_{\epsilon}$, then this metric satisfies the conditions  of Lemma \ref{Kmetric}, and the pullback of any product of sections of $V_{v}$ that are smaller than $1$ is a section of $V$ that is smaller than $1$.
 
We have now made all the choices to define the sheaves $S$ and $S_{v}$ so that any choice of sections of $S_{v}$ over $\hat f_{i_{v}}$ for all $v$ pull back to a section of $S$ over $\hat f_{i}$.  Recall from Definition \ref{Sdef} that a section of $S$ over $\hat f_{i}$ is a section $\nu$  of $V_{i}$ so that the following holds:
\begin{itemize}
\item $\dbar\hat f_{i}-\nu$ is contained in $V_{j}\subset V_{i}$  on a neighborhood of wherever $\rho_{j}\geq 0$; a stronger condition is satisfied by sections pulled back from $\hat f_{i_{v}}$. These pulled back sections will be contained in $\oplus_{v} V_{j_{v}}$ on some neighborhood of where $\rho_{j_{v}}\geq 0$ for all $v$. This is stronger because  $j_{v}$ is allowed to depend on $v$.
 \item $\nu$ is close to $\dbar\hat f_{i}$ in the sense that $\abs{\dbar\hat f_{i}-\nu}<1$; we have chosen our metrics so that pullbacks of sections of $\prod_{v}S_{v}$ satisfy this condition.   
 \item $\nu$ is transverse to 0; pullbacks of sections of $\prod_{v}S_{v}$ also clearly satisfy this condition. 
 \end{itemize}
 
 Similarly, if a family $\hat f$ in $\mathcal K_{\epsilon}$ is a product of families $\hat f_{v}$ in $\mathcal K_{v,\epsilon}$ any choice of section of $S_{v}$ over $\hat f_{v}$ for all $v$ pulls back to a section of $S$ over $\hat f$.
 
 As our weighted branched cover $I$ of $\mathcal K_{\epsilon}$ is the pull back of the corresponding weighted branched covers $I_{v}$ of $\mathcal K_{\epsilon}$, we may also pull back a choice of global section of $S_{v}^{I_{v}}$ over $\mathcal K_{v,\epsilon}$ for all $v$ to give a global section of $S^{I}$ over $\mathcal K_{\epsilon}$. In particular, suppose that $\hat f$ in $\mathcal K_{\epsilon}$ is a product of families $\hat f_{v}$ in $\mathcal K_{\epsilon,v}\cap \mathcal O_{I_{v}}$, then $\hat f$ is in $\mathcal O_{I}$, and $I(\hat f)=\prod_{v} I_{v}(\hat f_{v})$. On $\hat f_{v}$ a section of $S_{v}^{I_{v}}$ is a choice of section $\nu_{v}(i)$ of $S_{v}(\hat f_{v})$ for all $i\in I_{v}$ so that if $i$ and $j$ are not separated at $f\in\hat f_{v}$, then $\nu_{v}(i)=\nu_{v}(j)$ on a neighborhood of $f$. Taking the product over $v$ of such sections  gives  sections $\nu(i)$ of $S_{v}(\hat f)$ for all $i$ in $I(\hat f)$, where if $i$ is the product of $i_{v}$, then $\nu(i)$ is the product of $\nu_{v}(i_{v})$. Again, $\nu(i)$ will be equal to $\nu(j)$ on a neighborhood of any curve $f\in \hat f$ where $i$ is not separated from $j$, so $\nu$ defines a section of $S^{I}(\hat f)$.   
 
 Any map $\hat f\longrightarrow \hat g$  in $\mathcal K_{\epsilon}$ between families of curves that are products of families in  $\mathcal K_{v,\epsilon}$ corresponds to maps $\hat f_{v}\longrightarrow \hat g_{v}$. Global sections of $S_{v}^{I_{v}}$ are compatible with pullbacks, so the corresponding sections of  $S^{I}(\hat f)$ and $S^{I}(\hat g)$ are compatible with pullbacks. As any family in $\mathcal K_{\epsilon}$ is locally a product of families in $\mathcal K_{v,\epsilon}$, it follows that pulling back global sections of $S_{v}^{I_{v}}$ give global sections of $S^{I}$.
 
 Intersecting a global section of $S^{I}$ with the zero section  gives the virtual  class $[\mathcal K]$; $\nu(i)$ intersected with $0$ gives $[\mathcal K](i)$, an oriented family in $\mathcal K_{C}\subset\mathcal K_{\epsilon}^{st}$. Unfortunately, the map $\mathcal K^{st}\longrightarrow\prod_{v}\mathcal K_{v}^{st}$ does not pull back $\prod_{v}\mathcal K_{v,C}$ to be contained in $\mathcal K_{C}$, but actually, we have that these families $[\mathcal K_{v}](i)$ are contained in the substack of $\mathcal K'_{v,C}\subset \mathcal K_{v,C}$ where $\dbar<1$, and we have constructed our metrics so that $\prod_{v}\mathcal K'_{v,C}$ pulls back to be contained in $\mathcal K_{C}$.

Make the auxiliary choices $(W,x,e)$ required by Definition \ref{pushforward} for defining $\pi_{!}$ and $(\pi_{v})_{!}$ so that the diagram 
\begin{equation}\label{product diagram}\begin{tikzcd}\ro^{*}(\mathcal K)\rar{\pi_{!}}&\ro^{*}(\prod_{v}\ex A_{v})
\\ \prod_{v}\ro^{*}(\mathcal K_{v})\uar \ar{ur}[swap]{\prod_{v}(\pi_{v})_{!}}\end{tikzcd} \end{equation}
 commutes. In particular, choose $(W_{v},x_{v},e_{v})$ satisfying the requirements of Definition \ref{pushforward} so that $x_{v}:W_{v}\longrightarrow \prod_{v}\ex A_{v}$ is a submersion restricted to fibers. Then we can take $W=\prod_{v}W_{v}$, $x=\prod_{v}x_{v}:W\longrightarrow \prod_{v}\ex A_{v}$, and $e=\prod_{v}e_{v}$. As $x$ is a submersion restricted to fibers, $(W,x,e)$ also satisfies the requirements of Definition \ref{pushforward}. To avoid needing to think about orientation complications, we may also choose each $W_{v}$ to have even rank so that the Thom forms $e_{v}$ are even dimensional.
 
 We need to prove that Diagram (\ref{product diagram}) commutes. As the maps in Diagram (\ref{product diagram}) above are linear, we may use partitions of unity (on $\mathcal K_{v}^{st}$) to reduce to the case of differential forms with small support on $\mathcal K_{C}$ and $\mathcal K'_{v,C}$. In particular, without losing generality, assume  the following. 
 \begin{itemize}
 \item $\theta\in\ro^{*}\mathcal K$ has support on $\mathcal K_{C}$ compactly  contained in $\hat f/G$.
 \item $\hat f$ is in $\mathcal K_{\epsilon}\cap\mathcal O_{I}$ and is a product of families $\hat f_{v}$ in $\mathcal K_{v,\epsilon}\cap\mathcal O_{I_{v}}$.
 \item $\hat f/G$ and $\hat f_{v}/G_{v}$ represent substacks.
 \item $G=\prod_{v}G_{v}$.
 \item $\theta$ is the product of the pull back of $\theta_{v}\in \ro^{*}(\mathcal K_{v})$.
 \item $\theta_{v}$ has support on $\mathcal K'_{v,C}$ compactly contained in $\hat f_{v}/G_{v}$.
 \item There are partitions of unity  (Definition \ref{new partition def}) compatible with $[\mathcal K]$  and $[\mathcal K_{v}]$ so one  of the connected families $O_{1}$ involved is $\hat f$ (respectively $\hat f_{v}$), and so that the corresponding function $r$ is $1/\abs G$ (respectively $1/\abs {G_{v}}$) when restricted to the support of $\theta$ (respectively $\theta_{v}$) on $\hat f$ (respectively $\hat f_{v}$). 
 \end{itemize}
  Then $I(\hat f)=\prod_{v}I_{v}(\hat f_{v})$, and for $i=\prod_{v}i_{v}\in I(\hat f)$, we have $[\mathcal K](i)=\prod_{v}[\mathcal K_{v}](i_{v})$. We have also chosen $(W,x)$ and $(W_{v},x_{v})$ so that the map $\hat\pi(i):\pi(i)^{*}W\longrightarrow \prod_{v}\ex A_{v}$ is equal to the product of the maps $\hat \pi_{v}(i_{v}):\pi_{v}(i)^{*}W_{v}\longrightarrow \prod_{v}\ex A_{v}$ from Definition \ref{pushforward} part \ref{pushforward t}. It follows that $\hat\pi(i)_{!}$ is the product of $\hat \pi_{v}(i_{v})_{!}$. Then, for $\theta$ and $\theta_v$ satisfying the conditions listed above, 
 \begin{equation}\label{6}\begin{split}\pi_{!}(\theta)&=\frac 1{\abs G}\sum _{i\in I(\hat f)}\mu(i)\hat\pi(i)_{!}(\theta\wedge e)
 \\&=\prod_{v}\frac 1{\abs {G_{v}}}\sum _{i_{v}\in I_{v}(\hat f_{v})}\mu(i_{v})\hat\pi_{v}(i_{v})_{!}(\theta_{v}\wedge e_{v})
 =\prod_{v}(\pi_{v})_{!}(\theta_{v})\ .
 \end{split}\end{equation}
 
 Note that $\pi_{!}$ does not depend  on the choice of partition of unity,  so Diagram (\ref{product diagram}) commutes on forms that are a sum of forms obeying the condition listed above and hence satisfying equation \ref{6}. Any choice of form in $\prod_{v}\ro^{*}(\mathcal K_{v})$ may be expressed as a sum of forms obeying the conditions listed above and forms that for some $v$ vanish on $\mathcal K'_{v,C}$. Both $\pi_{!}$ and $\prod_{v}(\pi_{v})_{!}$ vanish on such vanishing forms, therefore Diagram (\ref{product diagram}) commutes as required.
 
On the level of cohomology,  $\pi_{!}$ and $\prod_{v}(\pi_{v})_{!}$  are independent of any choices, therefore in general, the following diagrams commute
 \[\begin{tikzcd}\rh^{*}(\mathcal K)\rar{\pi_{!}}&\rh^{*}(\prod_{v}\ex A_{v}) & \rhf^{*}(\mathcal K)\rar{\pi_{!}}&\rhf^{*}(\prod_{v}\ex A_{v})
\\ \prod_{v}\rh^{*}(\mathcal K_{v})\uar\ar{ur}[swap]{\prod_{v}(\pi_{v})_{!}}&&  \prod_{v}\rhf^{*}(\mathcal K_{v})\uar\ar{ur}[swap]{\prod_{v}(\pi_{v})_{!}}\end{tikzcd}\]
as required.

\stop

\section{Tropical completion}
\label{tropical completion section}

 Given a point $p$ in the tropical part $\totb{\ex B}$ of a compact exploded manifold $\ex B$, the set of points $\ex B\rvert_{p}\subset \ex B$ with tropical part $p$ is a possibly non-compact manifold. The tropical completion $\ex B\tc{p}$ of $\ex B$ at $p$ is a canonical  way of completing $\ex B\rvert_{p}$ to a complete exploded manifold. It should be thought of as a minimal way of adding `structure at infinity' to $\ex B\rvert_{p}$ in a way determined by $\ex B$.  In this section, we show that tropical completion is compatible with our construction of virtual class $[\mathcal K]$, as well as the associated integration and pushforward of forms from $\mathcal K$. This is used in the proof of the tropical gluing formula, equation (1) of \cite{gfgw}, where the contribution of a tropical curve $\gamma$ is defined using tropical completion. 

\begin{defn}[Tropical completion in a coordinate chart]\label{tc chart} 
The tropical completion of  a polytope $P$ in $\mathbb R^{m}$ at a point $p\in P$ is a polytope $\check P_{p}\subset \mathbb R^{m}$ which is the union of all rays in $\mathbb R^{m}$  beginning at $p$ and intersecting $P$ in more than one point.

Given an open subset  $U\subset \mathbb R^{n}\times\et mP$ and a point $p\in \totb U$.   let $ \bar U_{p}$ indicate the closure within $U$ of all points with tropical part $p$.   Note that $\bar U_{p}$ is naturally contained in both $\mathbb R^{n}\times \et mP$ and $ \mathbb R^{n}\times \et m{\check P_{p}}$. Define the tropical completion of the coordinate chart $U$ at $p\in\totb U$
\[ U\tc{p}\subset \mathbb R^{n}\times \et m{\check P_{p}}\]
to be the smallest open subset of $\mathbb R^{n}\times \et m{\check P_{p}}$ that contains $\bar U_{p}$.

Similarly, given a countable collection of points $X\subset \totb{ U}$, define 
\[U\tc X:=\coprod_{p\in X}U\tc p\ .\]

\end{defn}

Tropical completion in coordinate charts is functorial: given a map $\totb f:P\longrightarrow Q$ and point $p\in P$,  there is a unique map $\totb{ f\tc{p}}: \check P_{p}\longrightarrow \check Q_{f(p)} $  restricting to  $P\subset \check P_{p}$ to be equal to $\totb f$, and
similarly,  given a smooth or $\C\infty1$ map of coordinate charts
 \[f:U\longrightarrow U'\]  and a point $p$ in $\totb U$, there is a unique map 
\[ f\tc{p}: U\tc{p}\longrightarrow  U'\tc{f(p)}\]
so that $ f\tc{p}$  is equal to $f$ when restricted to the inverse image of $p$ in $U$.   Of course, the tropical part of $ f\tc{p}$ is equal to the map $\totb { f\tc{p}}$ above. This map $ f\tc{p}$ is smooth or $\C\infty1$ if $f$ is, is an isomorphism onto an open subset of $ U'\tc{f(p)}$ if $f$ is an isomorphism onto an open subset of $U'$, and is complete if $f$ is proper. There is therefore a functorial construction of the tropical completion $\ex B\tc{p}$ of any exploded manifold  $\ex B$ at a point $p\in \totb{\ex B}$ defined by applying tropical completion to coordinate charts. 

In particular, let $ O\longrightarrow \ex B$ be a cover of $\ex B$ by a collection of coordinate charts. We can  recover $\ex B$  by gluing together these coordinate charts using the transition maps, encapsulated in $O\times_{\ex B}O\rightrightarrows O$ --- in other words, $\ex B$ is a pushout as in the following diagram.
\[\begin{tikzcd}\ex B&\lar O
\\ \uar O& O\times_{\ex B}O\uar{t}\lar{s}\end{tikzcd}\]
For such a pushout to exist as an exploded manifold, we require that $O$ and $O\times_{\ex B}O$ are exploded manifolds,  $s$ and $t$ are \'etale (or local isomorphisms), and that $(s,t):O\times_{\ex B}O\longrightarrow O\times O$ is injective and proper. Such properties are preserved by tropical completion.

\begin{defn}[tropical completion of an exploded manifold]\label{tropical completion} Given $p\in\totb{\ex B}$, the tropical completion of $\ex B$ at $p$ is defined as follows. Choose an open cover of $\ex B$ by coordinate charts $\pi\co O\longrightarrow \ex B$ with transition maps encoded by $s,t\co O\times_{\ex B}O\rightrightarrows O$. Then $\ex B\tc p$ is the pushout created using the tropical completion of $s$ and $t$ at the inverse image of $p$ --- so $\ex B\tc p$ has charts and transition maps the tropical completion of those from $\ex B$.
\[\begin{tikzcd}\ex B\tc p&\lar O\tc{\totb\pi^{-1}p}
\\ O\tc{\totb\pi^{-1}p}\uar& \lrb{O\times_{\ex B}O}\tc{(\totb\pi\circ\totb s)^{-1}p}\uar \lar\end{tikzcd}\]
If $X\subset \totb {\ex B}$ is a countable set of points in the tropical part of $\ex B$, the tropical completion of $\ex B$ at $X$ may be defined as above, with $X$ replacing $p$. The resulting exploded manifold ${\ex B}\tc{X}$ is the disjoint union of $\ex B\tc p$ for all $p\in X$.
\end{defn}

Clearly, $\ex B\tc p$ is well defined independent of the particular collection of coordinate charts chosen, because applying tropical completion to extra coordinate charts will just give compatible coordinate charts on $\ex B\tc p$.

If $\ex B\rvert_{p}$ indicates the subset of $\ex B$ with tropical part $p$, note that   $\ex B\tc{p}$ always contains $\ex B\rvert_{p}$ as a dense subset. If $\ex B$ is basic, then ${\ex B}\tc{p}$ also contains a copy of the closure of  $\ex B\rvert_{p}$ in $\ex B$. (In the case that $\ex B$ is not basic, a single point in the closure of $\ex B\rvert_{p}$ may correspond to multiple points in $\ex B\tc p$.)

The construction of tropical completions is functorial in the sense that given a map
\[f:\ex A\longrightarrow \ex B\]
and a point $p\in \totb{\ex A}$, there exists a unique map
\[ f\tc{p}:{\ex A}\tc{p}\longrightarrow {\ex B}\tc{f(p)}\]
 restricting to be $f$ on $\ex A\rvert_{p}\subset {\ex A}\tc{p}$.

Note that $\ex B\tc p$ is complete if $\ex B$ is compact, and $f\tc p$ is complete if $f$ is proper.

Because  $\mathbb R$ is also an exploded manifold (with tropical part a single point, and tropical completion at this point still $\mathbb R$), the tropical completion of any $\mathbb R$--valued function on an exploded manifold is still a $\mathbb R$--valued function. Moreover,  given any tensor $\theta$ on $\ex B$ such as an almost complex structure, metric, or differential form, there is a unique tensor $ \theta\tc{p}$ on ${\ex B}\tc{p}$ restricting to be $\theta$ on $\ex B\rvert_{p}$. If $\theta$ is a refined differential form on $\ex B$, the restriction of $\theta\tc{p}$ to $\ex B\rvert_{p}$ no longer uniquely specifies $\theta\tc p$; in this case $ \theta\tc{p}$ is the unique refined differential form on ${\ex B}\tc{p}$ so that if $\theta$ pulls back to an honest form $\theta'$ under $r:U'\longrightarrow \ex B$, then $ \theta\tc{p}$ pulls back to $\theta'\tc{p'}$ under $ r\tc{p'}: U'\tc{p'}\longrightarrow  {\ex B}\tc{p}$, where $\totb r(p')=p$. Be warned that $\theta\tc {p}$ may not be in $\Omega^{*}\ex B\tc p$, even if $\theta\in\Omega^{*} \ex B$, because $\theta\tc{p}$ may not vanish on the image of all $\et 1{(0,\infty)}$. If, however,  $\theta\in \rof^{*}\ex B$, then $\theta\tc p$ is in $\rof^{*}\ex B\tc p$, so differential forms generated by functions are more compatible with tropical completion. In fact, Remark \ref{gbf} implies that $\theta\in \ro ^*(\ex B)$ is in $\rof^*(\ex B)$ if and only if $\theta\tc p\in \ro^*(\ex B)$ for all $p\in \totb{\ex B}$. 

\begin{lemma}\label{tc pushforward}Given any submersion, $f:\ex A\longrightarrow \ex B$, and form $\theta\in \rof_c^*(\ex A)$, the pushforward of $\theta$ is in $\rof^*_c{\ex B}$, and for any $p\in \ex B$, 
\[(f_!\theta)\tc p=\sum_{p'\in\totb f^{-1}p}(f\tc{p'})_!\theta\tc {p'} \ .\]
\end{lemma}

\pf As specified by Theorem 9.2 of \cite{dre}, $f_!\theta$ is uniquely determined by the property that 
\[\int_{\ex A} f^*\alpha \wedge \theta=\int_{\ex B} \alpha\wedge f_!\theta\]
for all $\alpha\in \ro^*(\ex B)$.  For any such $\alpha$ supported only in $\ex B\vert_p$, the pullback of $\alpha$ is supported only over the inverse image of $p$, so
\[\int_{\ex B\tc p}\alpha \tc p\wedge (f_!\theta)\tc p= \sum_{p'\in \totb f^{-1}p}\int_{\ex A\tc {p'}}(f\tc{p'})^*\alpha\tc p\wedge \theta\tc {p'}=\int_{\ex B\tc p} \alpha\tc p\wedge\sum_{p'\in \totb f^{-1}p}(f\tc{p'})_!\theta\tc {p'}\ .\]
 As any form on the smooth manifold $\ex B\vert_p$ is determined by its integral against compactly supported forms, and the restriction of $f_!\theta$ to $\ex B\tc p$ uniquely determines $(f_!\theta)\tc p$, if follows that
 \[(f_!\theta)\tc p=\sum_{p'\in\totb f^{-1}p}(f\tc{p'})_!\theta\tc {p'} \]
 as required. This equation then implies that $(f_!\theta)\tc p\in \ro^*(\ex B\tc p)$ for all $p\in \totb{\ex B}$, so Remark 5.5 implies that $f_!\theta\in \rof^*(\ex B)$. 
 
 \stop

We can define the tropical completion of an orbifold  $\mathcal X$ similarly, by applying tropical completion to an atlas: In particular, any surjective \'etale map from a (not necessarily connected) exploded manifold $U$ to $\mathcal X$ defines an \'etale proper\footnote{ An \'etale proper groupoid in the category of exploded manifolds is a groupoid with objects and morphisms parametrized by  exploded manifolds $\mathcal X_{0}$, $\mathcal X_{1}$, all structure maps morphisms in the category of exploded manifolds, source and target maps $s,t:\mathcal X_{1}\longrightarrow \mathcal X_{0}$ \'etale (ie locally isomorphisms), and the map $(s,t):X_{1}\longrightarrow X_{0}^{2}$ proper. The correct generalization  of `proper' for exploded manifolds is usually `complete', however if $s$ and $t$ are \'etale, then $(s,t)$ being proper is equivalent to it being complete. } groupoid $\mathcal X'$ with objects parametrized by $U$, and morphisms parametrized by $U\times_{\mathcal X}U$. We can recover $\mathcal X$ from this \'etale proper groupoid  $\mathcal X'$ by taking the stack of principal $\mathcal X'$--bundles. Taking tropical completion preserves the structure equations of a groupoid and the property of maps being \'etale or proper. So,  we can apply tropical completion to  $U\times_{\mathcal X}U\rightrightarrows U$ to obtain a new \'etale proper groupoid representing the tropical completion of our stack.

To define tropical completion of an orbifold, we need a notion of the tropical part of an orbifold. The following defines the tropical part of a stack as a set. The tropical part of a stack obviously has a bit more structure, but we shall not need it.

\begin{defn}As a set, the tropical part $\totb{\mathcal X}$ of a stack $\mathcal  X$ over the category of exploded manifolds is the set of path connected components of $\mathcal X$. Say that two  given points in $\mathcal X$ are connected by a path if there is a map of $\mathbb R$ into $\mathcal X$ that restricts to $\{0,1\}$ to be (isomorphic to) the two given points. 

Define the tropical part $\totb{\mathcal K}$ of a $K$--category $\mathcal K$ to be equal to the tropical part of $\mathcal K^{st}$.
\end{defn}

The following definition formalises the idea that tropical completion of an orbifold $\mathcal X$ is achieved by applying  tropical completion to coordinate charts on $\mathcal X$.

\begin{defn}[tropical completion of an orbifold] Given an exploded orbifold $\mathcal X$, and a countable subset $A\subset \totb{\ex X}$, 
the tropical completion $\mathcal X\tc A$ of $\mathcal X$ at $A$ is defined as follows: Choose an \'etale surjection $\ex U\longrightarrow \mathcal X$, and abuse notation a little to denote by $A$ the inverse image of $A$ in $\totb{\ex U}$ and within $\totb{\ex U\times_{\mathcal X}\ex U}$. Apply tropical completion at $A$ to obtain an \'etale proper groupoid with objects parametrized by  $\ex U\tc A$, morphisms parametrized by $(\ex U\times _{\mathcal X}\ex U)\tc A$, and groupoid structure maps the tropical completion at $A$ of the structure maps of the original groupoid representing $\mathcal X$. Then $\mathcal X\tc A$ is the orbifold represented by this groupoid. 

\end{defn}

Again $\mathcal X\tc A$ contains $\mathcal X\rvert_{A}$ as a dense suborbifold, and any map $\mathcal X\longrightarrow \mathcal Y$ sending $A$ to $A'$ induces a canonical map $\mathcal X\tc A\longrightarrow \mathcal Y\tc A$ that restricts to the original map on $\mathcal X\rvert_{A}$. Any (possibly refined)  differential form $\theta$ on $\mathcal X$ defines a differential form $\theta\tc A$ on $\mathcal X\tc A$.

To define the tropical completion of a $K$--category, we again use tropical completion of charts. Here we run into a problem --- $\mathcal K^{st}$ is not determined by these charts, just as the smooth structure on the union of an open half-space with a perpendicular line is not determined by transition data. The problem is that there exist families in $\mathcal K^{st}$ so that some points have no neighborhood sent inside a chart. This technical issue is resolved using extensions of our charts, as families in $\mathcal K^{st}$ are always locally contained in some extended chart.  The following describes the structure we need.

\begin{remark}\label{translation 1}

The following information is sufficient to define an extendable Kuranishi category:

\begin{itemize}
\item  a groupoid in the category of $\C\infty1$ exploded manifolds, and with objects parametrized by a countable disjoint union of  exploded manifolds $\coprod_{i}\ex F_{i}^{\sharp}$, where each $\ex F_{i}$ has fixed dimension,
\item and open subsets $\ex F_{i}\subset \ex F_{i}'\subset \ex F_{i}^{\sharp}$
\end{itemize} 
 satisfying the following conditions.
\begin{enumerate}
\item The morphisms from $\ex F_{i}^{\sharp}$ to itself are given by the action of a finite group $G_{i}$. So, 
\[\ex F_{ii}^{\sharp}=\ex F^{\sharp}_{i}\times_{\ex F^{\sharp}_{i}/G_{i}}\ex F_{i}^{\sharp}\ .\]
where morphisms with source and target $\ex F_{i}^{\sharp}$ are parametrized by $\ex F_{ii}^{\sharp}$.
We require the subsets $\ex F_{i}$ and $\ex F_{i}'$ to be  $G_{i}$--equivariant open subsets of $\ex F_{i}^{\sharp}$.

\item The morphisms with source $\ex F_{i}^{\sharp}$ and target $\ex F_{j}^{\sharp}$ are parametrized by an exploded manifold $\ex F^{\sharp}_{ij}$. For each $i$, we require that $\ex F^{\sharp}_{ij}=\emptyset$ for all but finitely many $j$. The source and target maps, 
\[\begin{tikzcd}\ex F_{i}^{\sharp}&\lar[swap]{\phi_{ij}}\ex F^{\sharp}_{ij}\rar{\phi_{ji}}&\ex F^{\sharp}_{j}\end{tikzcd}\]
must satisfy the following conditions:
\begin{enumerate}
\item\label{extension} If $\bar{\ex F}_{j}'$ indicates the closure of $\ex F_{j}'$ within $\ex F_{j}^{\sharp}$,  $\phi_{ij}(\phi_{ji}^{-1}(\bar{\ex F}'_{j}))\subset \ex F^{\sharp}_{i}$ is a closed subset of $\ex F_{i}^{\sharp}$. Moreover,  $\ex F'_{i}$ contains the closure of $\ex F_{i}\subset \ex F_{i}^{\sharp}$. 

\item\label{G fold cover} If $\dim\ex F_{i}\leq\dim \ex F_{j}$, then $\phi_{ij}$ is a $G_{j}$--fold cover of an open subset of $\ex F_{i}$, and $\phi_{ji}$ a $G_{i}$--fold cover of a exploded submanifold  of $\ex F_{j}$ (locally defined by the transverse vanishing of $\C\infty1$ functions.)
\end{enumerate}
\end{enumerate}

\

Given a $K$--category  $\mathcal K$ with extensions $\mathcal K\exte\mathcal K'\exte\mathcal K^{\sharp}$, we may obtain the above data by setting $\ex F_{i}^{\sharp}=\ex F(\hat f_{i}^{\sharp})$, $\ex F_{i}'=\ex F(\hat f_{i}')$, $\ex F_{i}=\ex F(\hat f_{i})$, and $\ex F^{\sharp}_{ij}=\ex F(\hat f^{\sharp}_{i}\times_{(\mathcal K^{\sharp})^{st}}\hat f^{\sharp}_{j})$. The groupoid in question is the full subcategory of $(\mathcal K^{\sharp})^{st}$ with objects the individual curves in these $\hat f_{i}^{\sharp}$.

\

Given the above data, we define a stack $\mathcal K^{st}$ as follows: A family $\hat f$ in $\mathcal K^{st}$ parametrized by $\ex F(\hat f)$ is
\begin{enumerate}
\item a collection of $G_{i}$--fold covers $X_{i}(\hat f)\subset X^{\sharp}_{i}(\hat f)\longrightarrow \ex F(\hat f)$ of subsets of $\ex F(\hat f)$ so that  every point in $\ex F(\hat f)$ is in the image of some $X_{i}(\hat f)$, and in the interior of the image of some $X^{\sharp}_{i}(\hat f)$\footnote{$X_{i}(\hat f)$ is the fiber product of $\ex F(\hat f)$ with $\ex F(\hat f_{i})=\ex F_{i}$ over $\mathcal K^{st}$.}; 
 \item $G_{i}$--equivariant maps
\[X^{\sharp}_{i}(\hat f)\longrightarrow \ex F^{\sharp}_{i}\]
that are $\C\infty1$ on the interior of $X_{i}^{\sharp}(\hat f)$, so that the inverse image of $\ex F_{i}\subset \ex F_{i}^{\sharp}$ is $X_{i}(\hat f)\subset X_{i}(\hat f)^{\sharp}$;

\item maps $X_{i}^{\sharp}(\hat f)\times_{\ex F(\hat f)}X_{j}^{\sharp}(\hat f)\longrightarrow \ex F_{ij}^{\sharp}$ so that both squares  below are fiber-product diagrams, (in the category of sets)
\[\begin{tikzcd}X^{\sharp}_{i}(\hat f)\rar&\ex F^{\sharp}_{i}
\\  X^{\sharp}_{i}(\hat f)\times_{\ex F(\hat f)}X^{\sharp}_{j}(\hat f)\uar\dar\rar& \ex F^{\sharp}_{ij}\uar\dar
\\ X^{\sharp}_{j}(\hat f)\rar&\ex F^{\sharp}_{j}\end{tikzcd}\]
and so that these maps define a map of groupoids (in the category of sets).
\[\begin{tikzcd}\coprod_{i,j}X^{\sharp}_{i}(\hat f)\times_{\ex F(\hat f)}X^{\sharp}_{j}(\hat f)\dar[shift left]\dar[shift right]
\rar& \coprod_{i,j}\ex F_{ij}^{\sharp}\dar[shift left]\dar[shift right]
\\ \coprod_{i} X^{\sharp}_{i}(\hat f)\rar& \coprod_{i}\ex F_{i}^{\sharp} \end{tikzcd}\]
 
\end{enumerate}

A morphism $\hat f\longrightarrow \hat g$ in $\mathcal K^{st}$ is
\begin{itemize}
\item a $\C\infty1$ map $\ex F(\hat f)\longrightarrow \ex F(\hat g)$,
\item and  $G_{i}$--equivariant maps $X_{i}^{\sharp}(\hat f)\longrightarrow X^{\sharp}_{i}(\hat g)$ so that the following diagrams commute
\[\begin{tikzcd}X^{\sharp}_{i}(\hat f)\rar\ar{dr}&X^{\sharp}_{i}(\hat g)\dar &X_{i}^{\sharp}(\hat f)\times_{\ex F(\hat f)}X_{j}^{\sharp}(\hat f) \ar{dr}\rar &X_{i}^{\sharp}(\hat g)\times_{\ex F(\hat f)}X_{j}^{\sharp}(\hat g)\dar 
\\ &\ex F^{\sharp}_{i} & & \ex F^{\sharp}_{ij}\end{tikzcd}\]
and so that the following  is a  fiber-product diagram.
\[\begin{tikzcd}X^{\sharp}_{i}(\hat f)\rar\dar&X^{\sharp}_{i}(\hat g)\dar
\\\ex F(\hat f)\rar&\ex F(\hat g)\end{tikzcd}\]
\end{itemize}

As defined above, $\mathcal K^{st}$ is a stack. Moreover, $\mathcal K^{st}$ is equivalent to  the original $\mathcal K^{st}$ in the case that $\mathcal K$ was defined as in Definition \ref{K category}.   Of course, when $\mathcal K$ is a subcategory of a nice moduli stack of curves, it is better to think of $\mathcal K^{st}$ as a stack of curves rather than use the above equivalent, but ad-hoc stack. 
\end{remark}

To see that $\mathcal K^{st}$ is a stack, let us verify Definition 4.3 of \cite{stacks}. Note that $\mathcal K^{st}$ has essentially unique pullbacks so is a category fibered in groupoids over the category of $\C\infty1$ exploded manifolds. Isomorphisms in $\mathcal K^{st}$ form a sheaf because morphisms are defined as maps satisfying local conditions. Each descent datum is effective because a descent datum is sufficient to construct the $G_i$--fold covers $X_i^\sharp$, and all other information is locally determined, so further choices form a sheaf.

\begin{example} $\mathcal K^{st}$ is an orbifold in the case that each $\ex F_{i}$ has the same dimension. Then the data above may be thought of as an atlas for $\mathcal K^{st}$, or an \'etale proper Lie groupoid representing $\mathcal K^{st}$. Condition \ref{G fold cover} gives \'etale, and condition \ref{extension} implies properness.  Conversely, given any orbifold, we may obtain the above data from a choice of $3$ nested atlases, with the closure of the charts from a given atlas  contained in the charts from the next atlas. 
\end{example}

\begin{defn}[Tropical completion of a $K$--category] Given an extendable $K$--category $\mathcal K$  (Definition \ref{K category}) and a countable subset $A\subset \totb{\mathcal K}$, define the tropical completion $\mathcal K\tc A$ of $\mathcal K$ at $A$ as follows:

Choose  extensions $\mathcal K^{\sharp}\rexte\mathcal K'\rexte\mathcal K$. Define $\ex F_{i}:=\ex F(\hat f_{i})$,  $\ex F'_{i}:=\ex F(\hat f'_{i})$, $\ex F_{i}^{\sharp}:=\ex F(\hat f_{i}^{\sharp})$, and define $\ex F^{\sharp}_{ij}:=\ex F(\hat f_{i}^{\sharp}\times_{(\mathcal K^{\sharp})^{st}}\hat f_{j}^{\sharp})$ to obtain the data from Remark \ref{translation 1}. Define $\mathcal K\tc A$ to be the $ K$--category with the data
\[\ex F_{i}\tc A\subset \ex F_{i}'\tc A\subset \ex F_{i}^{\sharp}\tc A\]
\[\begin{tikzcd}\ex F^{\sharp}_{i}\tc A&\lar\rar\ex F^{\sharp}_{ij}\tc A&\ex F_{j}^{\sharp}\tc A\end{tikzcd}\]
along with groupoid structure maps corresponding to the tropical completion of the groupoid structures maps from $\mathcal K$. As noted in Remark \ref{translation 1}, this data is sufficient to define an extendable $K$--category $\mathcal K\tc A$.
\end{defn}
%
%There is a functor from $\mathcal K$ to $\mathcal K\tc A$ sending $\hat f$ to a family $\hat f\tc A$. Similarly, given any family $\hat f$ in $\mathcal K^{st}$ and point $p\in\totb{\ex F(\hat f)}$ sent to $A\subset \mathcal K$, there is a corresponding family $\hat f\tc p$ in $(\mathcal K\tc A)^{st}$. We may think of this as a functor from the category of families  $\hat f$ in $\mathcal K^{st}$ along with a choice of point in $\totb{\ex F(\hat f)}$ sent to $A$.  This functor is not a map of stacks (over the  category of exploded manifolds with an extra choice of point in their tropical part) because the family $\ex F(\hat f\tc p)$ parametrizing $\hat f\tc p$ is not $\ex F(\hat f)$, but its tropical completion at $p$.

 Any map $\psi:\mathcal K\longrightarrow \mathcal X$ from $\mathcal K$ to an exploded orbifold or manifold $\mathcal X$ corresponds to a map $\psi\tc A:\mathcal K\tc A\longrightarrow \mathcal X\tc {\psi(A)}$. Any differential form $\theta$ on $\mathcal K$ corresponds to a differential form $\theta\tc A$ on $\mathcal K\tc A$ , any open substack $\mathcal U$ of $\mathcal K^{st}$ corresponds to an open substack $\mathcal U\tc A$ of $(\mathcal K\tc A)^{st}$, any vectorbundle $V$ on $\mathcal U$ corresponds to a vectorbundle $V\tc A$ on $U\tc A$, and any section $s$ of $V$  corresponds to a section $s\tc A$ of $V\tc A$.

\begin{defn}[Tropical completion of a Kuranishi category] The tropical completion of a Kuranishi category (Definition \ref{Kcat}) $\mathcal K$ at a countable subset $A\subset\totb{\mathcal K}$ is the $K$--category $\mathcal K\tc A$ along with the vectorbundles $V_{i}\tc A$ on the open substacks $\mathcal U_{i}\tc A$ and sections $\dbar\hat f_{i}\tc A:\ex F(\hat f_{i}\tc A)\longrightarrow V_{i}\tc A(\hat f\tc A)$. 
\end{defn}

Note that if $\pi:\mathcal K\longrightarrow \ex Z$ is proper (Definition \ref{K proper}), then $\pi\tc p:\mathcal K\tc p\longrightarrow \ex Z\tc{\pi(p)}$ is complete.

\begin{lemma} \label{tropical completion lemma}Given any proper, oriented Kuranishi category $\mathcal K$ and closed form $\theta\in\rof^{*}\mathcal K$, 
\[\int_{[\mathcal K]}\theta=\sum_{p\in\totb{\mathcal K}}\int_{[\mathcal K\tc p]}\theta\tc p\ .\]
Similarly, if  $\mathcal K$ is complete and oriented over $\ex Z$, and $\pi:\mathcal K\longrightarrow \ex A$ is a compatible   map, then given any closed $\theta\in \rof^{*}\mathcal K$ and $p'\in\totb{\ex A}$, 
\[\pi_{!}(\theta)\tc {p'}=\sum_{p\in\totb{\pi}^{-1}(p')}(\pi\tc {p})_{!}(\theta\tc {p})\ \]
where the above equation holds exactly for some choices of construction of $\pi_!$ and $(\pi\tc p)$, and otherwise holds on the level of cohomology.
\end{lemma}

\pf

The construction of $[\mathcal K]$ is compatible with tropical completion. The tropical completions of the functions $\rho_{i}$ used to define $\mathcal K_{\epsilon}$ in Definition \ref{Kdef} are appropriate for defining $(\mathcal K\tc p)_{\epsilon}$.  Tropical completion applied to the metric from Lemma \ref{Kmetric} also gives an appropriate metric on $V$ over $\mathcal K\tc p$. With these choices,  the sheaf $S$ from Definition \ref{Sdef} is compatible with tropical completion --- if $\nu$ is a section of $S(\hat f)$, then $\nu\tc p$ is a section of $S(\hat f\tc p)$.

Any  weighted branched cover $I$ of $\mathcal K$ gives a weighted branched cover $I\tc p$ of $\mathcal K\tc p$ as follows: A connected family $ O$ in $\mathcal K^{st}$ is in $\mathcal O_{I\tc p}$ if there exists a family $O'\in \mathcal K\cap\mathcal O_{I}$ and a map $O\longrightarrow O'\tc p$. Note that given any two such families $O_{i}'$, a connected component of $O_{1}'\times_{\mathcal K^{st}}O_{2}'$ satisfies the same conditions because tropical completion commutes with fiber products. We can therefore define $I\tc p(O)$ as the inverse limit of $I(O')$ over the category of  all such maps $O\longrightarrow O'\tc p$. This $I\tc p$ is a functor because any map $O_{1}\longrightarrow O_{2}$ induces a corresponding map of inverse limits --- simply compose $O_{2}\longrightarrow O'\tc p$ with $O_{1}\longrightarrow O_{2}$. Because $I(O')$ only maps to a finite number of other finite probability spaces, inverse limits are easy: there exists some $O'$ with a map $O\longrightarrow O'\tc p$ so that $I(O')=I\tc p(O)$.  This $I\tc p$  satisfies the requirements of Definition \ref{wb def}. Moreover,  $I\tc p$ is separated if $I$ is, because, given any automorphism $\psi$  of a family $O$ in $\mathcal K\tc p$, there exists a family $O'$ in $\mathcal K$ with an automorphism $\psi'$ so that $O$ is a connected component of $O'\tc p$ and $\psi=\psi'\tc p$. 

Any global section  $\nu$ of $S^{I}$ over $\mathcal K_{\epsilon}$ corresponds to a global section $\nu\tc p$ of $S^{I\tc p}$ over $(\mathcal K\tc p)_{\epsilon}$.  In particular, given any $O\in \mathcal O_{I\tc p}\cap(\mathcal K\tc p)_{\epsilon}$, there exists a family $O'\in \mathcal O_{I}\cap \mathcal K_{\epsilon}$ with a map $O\longrightarrow O'\tc p$ so that $I\tc p(O)=I(O')$. Then, for any $i\in I\tc p(O)$, we have a section $\nu(i)$ of $S(O')$. The tropical completion of this section is a section $\nu(i)\tc p$ of $S(O'\tc p)$, which pulls back under $O\longrightarrow O'\tc p$ to define $\nu\tc p(i)$ as a section of $S(O)$. This defines the section $\nu\tc p$ of $S^{I\tc p}$ over $O$. 

Let us check that $\nu\tc p$ as defined above gives a well-defined global section of $S^{I\tc p}$. Given any map $O_{2}\longrightarrow O$ and map $O_{2}\longrightarrow O'_{2}\tc p$ with $I\tc p(O_{2})=I(O'_{2})$, there exists some $O''_{2}$ with maps so that the following diagram commutes
\[\begin{tikzcd}O_{2}\ar{rr}\dar \ar{dr}&&O\dar
\\ O'_{2}\tc p&O''_{2}\tc p\lar\rar&O'\tc p\end{tikzcd} \]
It follows that the section $\nu\tc p$ of $S^{I\tc p}(O_{2})$ pulled back from $O'_{2}$ coincides with the section pulled back from $O''_{2}$ and the section pulled back from $O'$ via $O$. Therefore, $\nu\tc p$ as defined above is a well-defined global section of $S^{I\tc p}$.

The resulting $[\mathcal K\tc p]$ is similarly related to $[\mathcal K]$. In particular, for any $O\in\mathcal O_{I\tc p}\cap (\mathcal K\tc p)_{\epsilon}$, we can choose an $O'$ with a map $O\longrightarrow O'\tc p$ so that $I\tc p(O)=I(O')$. Then $[\mathcal K\tc p](i)$ is the subfamily of $O$ that is the pullback of $[\mathcal K](i)\tc p\subset O'\tc p$ under the map $O\longrightarrow O'\tc p$.

\begin{claim}\label{local K int} If  $\mathcal K$ is proper and oriented, and $\theta\in\rof^{*}(\mathcal K)$ is any (not necessarily closed) form, then 
\begin{equation}\label{lki}\int_{[\mathcal K]}\theta=\sum_{p\in\totb{\mathcal K}}\int_{[\mathcal K\tc p]}\theta\tc p\end{equation}
when the above method is used to construct the virtual fundamental class $[\mathcal K\tc p]$ of $\mathcal K\tc p$ from $[\mathcal K]$.
\end{claim}

Choose a partition of unity $r:O\longrightarrow \mathbb R$ compatible with $[\mathcal K]$. Then, 
\[\int_{[\mathcal K]}\theta=\sum_{O_{k}\subset O}\sum_{i\in I(O_{k})}\mu(i)\int_{[\mathcal K](i)}r\theta\]
where the sum is over connected components $O_{k}$ of $O$. We may apply tropical completion to $r$ and obtain a partition of unity $r\tc p:O\tc p\longrightarrow \mathbb R$ compatible with $\mathcal K\tc p$, where  we use $p$ in $\totb O$ to mean the inverse image of $p$ under the map $\totb O\longrightarrow \totb {\mathcal K}$. 

For any exploded manifold $\ex B$, the integral of $\theta$ over $\ex B$ is the sum of the integral of $\theta$ over the manifolds $\ex B\rvert_{ p}$ (if $\theta$ has compact support, only a finite number of these integrals will be nonzero.) Furthermore, the integral of $\theta$ over $\ex B\rvert_p$ is equal to the integral of $\theta\tc p$ over $\ex B\tc p$. Applying this to our situation, we get
\[\int_{[\mathcal K](i)}r\theta=\sum_{p\in \totb {\mathcal K}}\int_{[\mathcal K](i)\tc p}r\tc p\theta\tc p\]
The tropical completion of $O_{k}$ may not be connected.  Each connected component $O'$ of $O_{k}\tc p$  is in $\mathcal O_{I\tc p}$, and there is a corresponding map $\iota^{*}:I(O_{k})\longrightarrow I\tc p(O')$ so that $[\mathcal K\tc p](\iota^{*}i)$ is the corresponding collection of connected components of $[\mathcal K](i)\tc p$. Therefore
\[\int_{[\mathcal K](i)\tc p}r\tc p\theta\tc p=\sum_{\iota}\int_{[\mathcal K\tc p](\iota^{*}i)}r\tc p\theta\tc p\]
where the sum is over the different  inclusions $\iota$ of connected components of $ O_{k}\tc p$. 
Noting that each $\iota^{*}$ is surjective and measure preserving then gives 
\[\sum_{O_{k}\subset O}\sum_{i\in I(O_{k})}\mu(i)\int_{[\mathcal K](i)\tc p}r\tc p\theta\tc p=\sum_{O'\subset O\tc p}\sum_{i'\in I(O')}\mu(i')\int_{[\mathcal K\tc p](i')}r\tc p\theta\tc p=\int_{[\mathcal K\tc p]}\theta\tc p \]
which implies  Claim \ref{local K int}.

In the case that $\theta$ is closed, (\ref{lki})  holds regardless of any choices made in the construction of $[\mathcal K]$ and $[\mathcal K\tc p]$ because these integrals do not depend on such choices. It remains to prove the analogous statement for pushforwards:

\begin{claim}\label{tropical pushforward}
 If  $\mathcal K$ is complete and oriented over $\ex Z$, and $\pi:\mathcal K\longrightarrow \ex A$ is a compatible   map,  there exists a construction of $\pi_{!}$ and $(\pi\tc p)_{!}$ so that  given any  $\theta\in \rof^{*}\mathcal K$ and $p'\in\totb{\ex A}$,
\[\pi_{!}(\theta)\tc {p'}=\sum_{p\in\totb{\pi}^{-1}(p')}(\pi\tc {p})_{!}(\theta\tc {p})\] 
\end{claim}

The proof of Claim \ref{tropical pushforward} is very similar to the proof of Claim \ref{local K int} except we must now check that all choices in the construction of $\pi_{!}$ from Definition \ref{pushforward} are compatible with tropical completion. This is indeed the case: the choice of vectorbundle $W\longrightarrow \ex A$, map $x:W\longrightarrow \ex A$ and Thom form $e$ on $W$ all may be tropically completed at $p'$ to give $W\tc p\longrightarrow \ex A\tc p$, $x\tc p:W\tc p\longrightarrow \ex A\tc p$ and $e\tc p$ appropriate for defining $[\pi\tc p]_{!}$ for all $p\in\totb\pi^{-1}p'$. 

 Definition \ref{pushforward} gives that 
\[\pi_{!}(\theta)=\sum_{O_{k}\subset O}\sum _{i\in I(O_{k})} \mu(i)\hat \pi(i)_{!}(r\theta\wedge e)\]
where $\hat \pi(i)$ is defined by the following composition:
\[\begin{tikzcd}\pi(i)^{*}W\dar \rar\ar[bend left]{rr}{\hat \pi(i)}& W\dar \rar{x}& \ex A
\\{} [\mathcal K](i) \rar{\pi(i)}& \ex A\end{tikzcd}\]
Therefore, 
\[\pi_{!}(\theta)\tc {p'}=\sum_{O_{k}\subset O}\sum_{p\in\totb\pi^{-1}(p')}\sum_{i\in I(O_{k})}\mu(i)\hat\pi(i)\tc p(r\tc p\theta\tc p\wedge e\tc p)\ .\]
As in the proof of Claim \ref{local K int}, $O_{k}\tc p$ may have several connected components, $\iota: O'\longrightarrow O_{k}\tc p$, and $[\mathcal K](i)\tc p=\coprod_{\iota}[\mathcal K\tc p](\iota^{*}i)$, and we similarly have $\hat \pi(i)\tc p=\coprod_{\iota}\hat \pi\tc p(\iota^{*}i)$. The same argument as the proof of Claim \ref{local K int} then gives the required result:
\[\pi_{!}(\theta)\tc {p'}=\sum_{p\in\totb\pi^{-1}(p)}(\pi\tc p)_{!}(\theta \tc p)\]

In the case that $\theta$ is closed, both sides of the above equation are independent of all choices, therefore the above equation holds regardless of choices.

\stop

\section{Construction of Gromov--Witten invariants}\label{GW section}

Let us summarise our construction of Gromov--Witten invariants. Start with a complete, basic exploded manifold, $\ex B$, with a $\dbar$--log compatible almost complex structure, $J$, tamed by a symplectic form, $\omega$, as described in definitions 3.1 and 3.5 of \cite{cem} and definitions 3.15 and 4.6 of \cite{iec}. Also assume that the tropical part of $\ex B$ admits a $\mathbb Z$--affine immersion into some  $[0,\infty)^N$ so that Lemma 4.2 and Theorem 6.1 of \cite{cem} imply that the moduli stack of $J$--holomorphic curves is compact when restricted to curves with bounded genus,  energy, and number of punctures or ends.\footnote{An end of a holomorphic curve is a stratum of its domain, $\ex C$, isomorphic to $\et 1{(0,\infty)}$. The smooth part, $\totl{\ex C}$, of $\ex C$ is a nodal curve  with a marked point corresponding to each end of $\ex C$. We can also achieve compactness and construct Gromov--Witten invariants with weaker assumptions when we include more discrete data about the curves at ends. For example holomorphic curves in  $\ex T^{n}$ have tropical parts tropical curves in $\mathbb R^{n}$. Each end of such a curve corresponds to an end of the corresponding tropical curve, which has some derivative $\alpha$ in $\mathbb Z^{n}$. To achieve compactness for the moduli stack of holomorphic curves in $\ex T^{n}$, we have to keep track of these derivatives $\alpha$ at the ends of curves. See Lemma 4.2 of \cite{cem}.} Some examples of $(\ex B,\omega,J)$ satisfying these conditions include:
\begin{itemize}
\item any compact symplectic manifold with a tamed almost complex structure $J$;
\item the explosion\footnote{See section 5 of \cite{iec} for a discussion of the explosion functor, or see section 4 of  \cite{elc} for a discussion of the explosion functor as a base change within log geometry. The tropical part of the explosion of a manifold with a normal crossing divisor is the dual intersection complex with a canonical $\mathbb Z$--affine structure. To ensure that this tropical part immerses  in some $[0,\infty)^N$, we assume that it is  simple, or a union of transversely intersecting, compact, codimension 2, complex submanifolds.} of any compact K\"ahler manifold relative to simple normal crossing divisors;
\item and an exploded manifold corresponding to the singular fiber in any simple normal crossing degeneration of a compact K\"ahler manifold.
\end{itemize}
Let $\Msw_{g,n,\beta}(\ex B)$ be the moduli stack of stable $\C\infty1$ curves in $\ex B$ with genus $g$, $n$ labelled ends, and representing a homology class, $\beta:H^*(\ex B)\longrightarrow \mathbb R$.  As the moduli stack of holomorphic curves within $\Msw_{g,n,\beta}(\ex B)$ is compact, Theorem 7.3 of \cite{evc} provides a complete embedded Kuranishi structure covering the moduli stack of holomorphic curves within $\Msw_{g,n,\beta}(\ex B)$. Next, construct an orienting two-form $\alpha$ on the associated Kuranishi category, using Definition \ref{aKC}, then restrict the associated Kuranishi category to a neighborhood of  the holomorphic curves, so that  $\alpha$ is  orienting.  Call this oriented Kuranishi category $\mathcal K_{g,n,\beta}$. As $\mathcal K_{g,n,\beta}$ is complete,  we can construct a virtual fundamental class, $[\mathcal K_{g,n,\beta}]$, for $\mathcal K_{g,n,\beta}$ using Definition \ref{vdm}.

There are several distinct versions of `evaluation maps' from the moduli stack of curves in exploded manifolds. Let 
\[ev:\Msw_{g,n,\beta}(\ex B)\longrightarrow \ex X\]
be any $\C\infty1$ map to an exploded manifold or orbifold $\ex X$. For example, if $\ex B$ is a smooth manifold, or an exploded manifold with bounded tropical part, (such as any exploded manifold in a connected family also containing a smooth manifold), then each end of a curve is sent to a single point in $\ex B$, so we could take $ev$ as the usual evaluation map at ends, with $\ex X=\ex B^n$. If $\ex B$ has unbounded tropical part, $\totb{\ex B}$, there are different possible `evaluation' maps which we could use. In this case, some ends of holomorphic curves have unbounded image in $\totb{\ex B}$, and evaluation at such punctures lands in a space $2$ real dimensions smaller than $\ex B$. In such an exploded manifold, evaluation at an end lands in an associated exploded manifold, $\End \ex B$, constructed in section 3 of \cite{gfgw}, so we can take $\ex X$ as $(\End \ex B)^n$.

We can then define Gromov--Witten invariants as a map 
\[\rh^*(\ex X)\longrightarrow \mathbb R\]
defined by
\[\theta \mapsto \int_{[\mathcal K_{g,n,\beta}]} ev^*\theta\]
using the integration from Definition \ref{intdef} and the homology  from Definition \ref{refined def}. We can encode finer information about $[\mathcal K_{g,n,\beta}]$ using Definition 5.15 to push forward this virtual fundamental class to define
\[\eta_{g,n,\beta}:=ev_!(1)\in \rhf^*(\ex X)\]
so 
\[\int_{[\mathcal K_{g,n,\beta}]}ev^*\theta=\int_{\ex X}\theta \wedge \eta_{g,n,\beta}\ .\]

 If $2g+n\geq 3$,  there is also a $\C\infty1$ map from $\Msw_{g,n,\beta}(\ex B)$ to  Deligne--Mumford space, $\bar M_{g,n}$, or its explosion, $\expl \bar M_{g,n}$ --- this stabilisation map to $\expl \bar M_{g,n}$ is constructed in section 4 of \cite{evc}, where it is also proved that $\expl \bar M_{g,n}$ represents the moduli stack of stable exploded curves with genus $g$ and $n$ punctures. So, we could also define Gromov--Witten invariants with the $\ex X$ above being  $(\End \ex B)^n\times \expl \bar M_{g,n}$, or a more sophisticated version of this, $\mathcal X_{g,n}(\ex B)$,   described in section 5 of \cite{gfgw}. In the case that $\ex B$ is a symplectic manifold, Gromov--Witten invariants  defined using $\ex X:=\ex B^n\times \bar M_{g,n}$ satisfy Kontsevich and Manin's axioms for Gromov--Witten invariants from \cite{KM}; as explained in \cite{gfgw}, the splitting and genus reduction axioms follow from Theorem 5.2 and Lemma 5.3 of \cite{gfgw}.

We can also define descendant Gromov--Witten invariants incorporating Chern classes of tautological line bundles, $L_i^*$, over $\Msw_{g,n,\beta}(\ex B)$. Let us describe these tautological line bundles. Each of our labeled ends  corresponds to a $\et 1{(0,\infty)}$--bundle, $\ex L_i$, over $\Msw_{g,n,\beta}(\ex B)$, with fiber over a curve $f$ the stratum of the domain of $f$ labeled by $i$; so over a family, $\hat f$,  $\ex L_i(\hat f)\longrightarrow \ex F(\hat f)$ is a union of strata of the domain, $\ex C(\hat f)\longrightarrow \ex F(\hat f)$, of $\hat f$. The moduli stack of $\et 1{(0,\infty)}$--bundles is equivalent to the moduli stack of $\mathbb C^{*}$--bundles or  complex line bundles, so there is an associated complex line bundle $L_i$ on $\Msw_{g,n,\beta}(\ex B)$ --- the relationship between the  fibers of $L_i$  and $\ex L_i$ is the relationship between $\mathbb  C$ and $\et 1{(0,\infty)}\subset \et 1{[0,\infty)}=\expl (\mathbb C, 0)$. We can also think of the fiber of $L_i$ over $f$ as the tangent space of the smooth part\footnote{Every exploded manifold $\ex B$ comes with a natural map $\ex B\longrightarrow \totl{\ex B}$ to its smooth part; the smooth part of an exploded curve $\ex C$ is a nodal curve $\totl{\ex C}$; the smooth part of the explosion of a manifold with normal crossing divisors is the original manifold. } of the domain of $f$ at its $i$th marked point: The smooth part of the domain of $\hat f$, $\totl{\ex C(\hat f)}$, is a family of nodal curves over $\totl{\ex F(\hat f)}$, and the $i$th end corresponds to a  marked point  section $s_i:\totl{\ex F(\hat f)}\longrightarrow \totl{\ex C(\hat f)}$. Our line bundle $L_i$ is the pullback of the vertical tangent bundle of $\totl{\ex C(\hat f)}$ under the composition $\ex F(\hat f_i)\longrightarrow \totl{\ex F(\hat f_i)}\xrightarrow{s_i}\totl{\ex C(\hat f_i)}$. So, the dual line bundles,  $L_i^*$, are analogous to the usual tautological line bundles over the moduli stack of curves. 

With the tautological line bundles $L^*_i$ over $\Msw_{g,n,\beta}(\ex B)$ understood, we can construct their Chern classes  using Remark \ref{chern}, and define $\psi_i\in \rhf^*(\mathcal  K_{g,n,\beta})$ as the first Chern class of the line bundle $L^*_i$ over $\mathcal K$. We can then define the classes
\begin{equation}\label{dgw}\eta_{g,\{a_i\},\beta}:=ev_!\left(\prod_i\psi_i^{a_i}\right)\subset \rhf^*(\ex X)\end{equation}
for non-negative integers $a_1,\dotsc,a_n$, and define Gromov--Witten invariants as the corresponding maps
\begin{equation}\label{dgw2}\begin{tikzcd}\rh^*(\ex X)\ar{rrr}{\theta \mapsto \int_{\ex X}\eta_{g,\{a_i\},\beta}\wedge \theta} &&& \mathbb R
  \end{tikzcd}\ .\end{equation}

In the case that $\totb{\ex B}$ is bounded and $\ex X=\ex B^n$, we can package  Gromov--Witten invariants into correlators without losing important information. When $\totb{\ex B}$ is bounded, $\rh^*(\ex B)$ is isomorphic to $H^*(\ex B)$,\footnote{There is no written proof of this fact, but it can be proved by representing each cohomology class in $\rh^*(\ex B)$  by a closed form on a refinement of $\ex B$;  when $\totb{\ex B}$ is bounded and admits a $\mathbb Z$--affine immersion into $\mathbb R^N$, each refinement of $\ex B$ is cobordant to $\ex B$, and hence has the same cohomology by section 11 of \cite{dre}.} so we can use $H^*$ in place of $\rh^*$ without losing important information about Gromov--Witten invariants. Moreover, unlike $\rh^*$ and $\rhf^*$, K\"unneth's theorem applies to $H^*$, so $H^*(\ex B^n)=H^*(\ex B)^{\otimes n}$. Given classes $\theta_1,\dotsc,\theta_n$ in $H^*(\ex B)$, we can define the following numerical Gromov--Witten invariants,
 \begin{equation}\label{correlator}\langle \tau_{a_1}(\theta_1),\dotsc, \tau_{a_n}(\theta_n)\rangle_{g,n,\beta}:=\int_{\ex B^n}\eta_{g,\{a_i\},\beta}\bigwedge_i\pi_i^*\theta_i=\int_{[\mathcal K_{g,n,\beta}]}\bigwedge_i\psi^{a_i}_i\wedge ( \pi_i\circ ev)^*\theta_i\end{equation}
 where $\pi_i:\ex B^n\longrightarrow \ex B$ is projection onto the $i$th component. These $\tau_{a_i}(\theta_i)$ could also be thought of as cohomology classes on $\ex B$ times the stack of complex line bundles, where our correlator `integrates the pullback' of these classes to $[\mathcal K_{g,n,E}]$ using the evaluation map recording the position of the $i$th puncture and the line bundle $L_i^*$, however our machinery for  differential forms on stacks with infinite isotropy groups is inadequate for removing the above scare quotes. 
 
  Theorem \ref{cobordism theorem} together with Corollary 7.5 of \cite{evc} imply that $\eta_{g,\{a_i\},\beta}$ does not depend on the choices involved in its construction. Theorem 7.3 of \cite{evc} together with Theorem \ref{pullback theorem} provide a kind of invariance in families for the Gromov--Witten invariants (\ref{dgw2}) and (\ref{correlator}).

Let us describe how the correlators (\ref{correlator}) are invariant in families.   Let $\hat{\ex B}\longrightarrow \ex B_0 $ be a connected family of exploded manifolds containing $\ex B$, with a family of $\dbar$--log compatible almost complex structures tamed by a family of taming forms, and suppose that there is a $\mathbb Z$--affine map $\totb{\hat{\ex B}}\longrightarrow [0,\infty)^N$ that is injective when restricted to  each stratum of each fiber.  (For example, the explosion of a simple normal crossing degeneration satisfies this condition.)  In this case, the moduli stack of curves in $\hat{\ex B}$ with bounded genus, energy and number of punctures is proper over $\ex B_0$.

Before describing Gromov--Witten invariants in our family, let us consider how the cohomology of exploded manifolds varies in families.  It is not true that $\rhf^*$ is invariant in families of exploded manifolds, however Section 11 of \cite{dre} proves that $H^*$ is invariant in connected families of exploded manifolds.  More precisely, given any long path\footnote{See Definition 11.1 of \cite{dre}.} $\gamma$ in $\ex B_0$, with inverse image in $\hat{\ex B}$ joining $\ex B$ to a fiber, $\ex B'$,   of $\hat{\ex B}\longrightarrow \ex B_0$, Definition 11.3 of \cite{dre} gives an isomorphism $\Psi_\gamma:H^*(\ex B)\longrightarrow H^*(\ex B')$. This isomorphism depends on the isotopy class of $\gamma$. With $\Psi_{\gamma}$ understood, we can write the invariance of our correlators as follows.

   \begin{equation}\label{invariance}\langle \tau_{a_1}(\theta_1),\dotsc, \tau_{a_n}(\theta_n)\rangle_{g,n,\beta}=\langle \tau_{a_1}(\Psi_\gamma(\theta_1)),\dotsc, \tau_{a_n}(\Psi_\gamma(\theta_n))\rangle_{g,n,(\Psi_\gamma^{-1})^*\beta}\end{equation}

Let us first prove equation (\ref{invariance}) under the assumption that our long path $\gamma$ is contained in a small neighborhood of $\ex B\subset\hat{\ex B}$, and that there exists a  curve $f$ in $\ex B$ representing the homology class $\beta:H^{*}(\ex B)\longrightarrow \mathbb R$. Proposition 5.9 of \cite{evc} implies that $f$ locally extends to a connected family of (not-necessarily-holomorphic) curves, $\hat f$, with $\ex F(\hat f)\longrightarrow \ex B_0$ a submersion. Then,  given any long path $\gamma$ in the image of $\ex F(\hat f)$ joining the image of $f$ with the image of $f'$ the integral of $f^*\theta$ equals the integral of $(f')^*\Psi_\gamma(\theta)$. So, the homology class, $\beta:H^*(\ex B)\longrightarrow \mathbb R$, represented by $f$ locally extends to a map $\hat\beta:H^*(\ex B')\longrightarrow \mathbb R$ for all fibers $\ex B'$ in a neighborhood of $\ex B$ so that for any long path $\gamma$ in this neighborhood, $\hat \beta(\theta)=\hat\beta(\Psi_\gamma(\theta))$; (unlike the case of a family of smooth manifolds, not all homology classes locally extend in this way).  Let $\hat{\ex B}'\longrightarrow \ex B_{0}'$ be such a neighourhood of $\ex B$ where $\hat \beta$ exists, and let $\Msw_{g,n,\hat \beta}(\hat {\ex B}')$ be the moduli stack of stable curves with genus $g$ and $n$ punctures representing $\hat \beta$.

The moduli stack of holomorphic curves in $\Msw_{g,n,\hat \beta}(\hat {\ex B}')$ is complete over $\ex B_0'$, so using Theorem 7.3 of \cite{evc} we can construct an embedded Kuranishi structure with an associated Kuranishi category, $\hat{\mathcal K}_{g,n,\hat \beta}$, oriented and complete over $\ex B_0$. Then, Definition \ref{vdm} gives us a virtual fundamental class $[\hat{\mathcal K}_{g,n,\hat \beta}]$.

In this case, evaluation at punctures gives an evaluation map 
\[\hat{ev}:\Msw_{g,n,\hat \beta}(\hat{\ex B}')\longrightarrow \hat{\ex X}\]
where $\hat{\ex X}$ is the $n$th fiber product of $\hat {\ex B}'$ over $\ex B_0'$. The line bundle $L_i$ still make sense over $\Msw_{g,n,\hat \beta}(\hat {\ex B}')$, so we can define  $\psi_i\in \rhf^*(\hat{\mathcal K}_{g,n,\hat \beta})$ as the first Chern class of $L_i^*$ using Remark \ref{chern}. Because $\hat {ev}:\hat{\mathcal K}_{g,n,\hat \beta}\longrightarrow \hat{\ex X}$ is proper and relatively oriented,  we can define 
\[\hat \eta_{g,\{a_i\},\hat \beta}:=\hat{ev}_!\left(\prod_i\psi_i^{a_i}\right)\subset \rhf^*(\hat{\ex X})\]
 using Definition \ref{pushforward}. Then, Theorem \ref{pullback theorem} implies that  $\eta_{g,\{a_i\},\beta}$ is the pullback of $\hat \eta_{g,\{a_i\}, \hat \beta}$ under the corresponding inclusion $\ex B^n\subset \hat {\ex X}$, and if $\gamma$ is a long path in $\ex B_{0}'$ joining $\ex B$ to $\ex B'$, $\eta_{g,\{a_{i}\},(\Psi_{\gamma}^{-1})^{*}\beta}$ is the pullback of $\hat\eta_{g,\{a_{i}\},\hat \beta}$ under the inclusion $(\ex B')^{n}\subset \hat{\ex X}$.   In this case, equation (\ref{invariance}) holds  because  $\Psi_\gamma(\theta_i)$ is defined by extending $\theta_i$ to a closed form $\hat \theta_i$ over the long path $\gamma$, and $\Psi_\gamma(\theta_i)$ is the restriction of $\hat \theta_i$ to $\ex B'$; so, the integral of $\hat \eta_{g,\{a_i\},\hat \beta}\bigwedge\hat \pi_i^*\hat\theta_i $ over fibers is constant, and in particular, its integral over $\ex B^{n}$ equals its integral over $(\ex B')^{n}$. 
 
 So far, we have proved that equation (\ref{invariance}) holds for  long paths in $\ex B_{0}'$  starting at a fiber containing a  curve representing $\beta$. Now we argue that  equation (\ref{invariance}) always holds. As with usual paths, we can reparametrize a long path $\gamma$ into $\gamma_{1}$ followed by $\gamma_{2}$ so that $\gamma_{1}$ starts at the same point as $\gamma$, and ends at a chosen point in the middle of $\gamma$, and $\gamma_{2}$ starts at this chosen point, and ends where $\gamma$ ends. As usual, whenever  $\gamma_{1}$ followed by $\gamma_{2}$ is isotopic to a reparametrization of $\gamma$,  $\Psi_{\gamma}=\Psi_{\gamma_{2}}\circ\Psi_{\gamma_{1}}$, and reversed paths induce inverse isomorphisms. Let $S$ be the set of points in the domain of $\gamma$ for which equation (\ref{invariance}) holds --- so if $\gamma_{1}$ ends at a point in $S$, equation (\ref{invariance}) holds for $\gamma_{1}$. Let $\gamma_{1}$ end at a point in the closure of $S$, and let $\ex B'$ be the fiber over this point. If our Gromov--Witten invariant (\ref{correlator}) is nonzero,  then  $\ex B'$ must contain a holomorphic curve representing $(\Psi_{\gamma_{1}}^{-1})^{*}\beta$, because the moduli stack of holomorphic curves with bounded energy, genus, and number of ends is proper over $\ex B_{0}$. Then, the above  argument  shows that equation (\ref{invariance}) holds with $\beta$ replaced by $(\Psi_{\gamma_{1}}^{-1})^{*}\beta$, and $\gamma$ replaced by any sufficiently small long path starting at $\ex B'$. As equation (\ref{invariance}) holds for a path if and only if it holds for the reversed path, we can reparametrize $\gamma_{1}$ using two long paths for which equation (\ref{invariance}) holds, so equation (\ref{invariance}) holds for $\gamma_{1}$. Similarly, equation (\ref{invariance}) holds for any long path that can be reparametrized as $\gamma_{1}$ followed by a sufficiently small long path. It follows that the set $S$ where equation (\ref{invariance}) holds is both open and closed, and therefore includes the entire domain of $\gamma$. We have therefore shown that equation (\ref{invariance}) holds so long as our Gromov--Witten invariant is nonzero. If it is zero, and equation (\ref{invariance}) failed to hold, then our Gromov--Witten invariant at the other end of $\gamma$ would be nonzero, and therefore equation (\ref{invariance}) would hold for the reversed long path, and therefore must hold for $\gamma$ itself.

\bibliographystyle{plain}
\bibliography{ref.bib}
 \end{document}